\documentclass[12pt,leqno]{amsart} 
\theoremstyle{plain} 

\newtheorem{global-theorem}{Theorem}
\newtheorem{theorem}{Theorem}[section]
\newtheorem{lemma}[theorem]{Lemma}

\newtheorem{hypothesis}[theorem]{Hypothesis}
\newtheorem{corollary}[theorem]{Corollary}

\newtheorem{proposition}[theorem]{Proposition}

\newtheorem{prop-def}[theorem]{Proposition-Definition}
\newtheorem{lemma-def}[theorem]{Lemma-Definition}

\newcommand{\eop}{\ \hfill $\Box$}

\newcommand{\revision}{\null}
\numberwithin{equation}{section}

\newcommand{\cc}{{\mathbb C}}
\newcommand{\pp}{{\mathbb P}}
\newcommand{\rr}{{\mathbb R}}
\newcommand{\qq}{{\mathbb Q}}
\newcommand{\zz}{{\mathbb Z}}
\newcommand{\ff}{{\mathbb F}}
\newcommand{\Gm}{{\mathbb G}_{\rm m}}

\newcommand{\Tt}{{\mathcal T}}
\newcommand{\Dd}{{\mathcal D}}
\newcommand{\Ff}{{\mathcal F}}
\newcommand{\Oo}{{\mathcal O}}
\newcommand{\Uu}{{\mathcal U}}

\newcommand{\Vv}{{\mathcal V}}
\newcommand{\mcS}{{\mathcal S}}

\newcommand{\Aa}{{\mathcal A}}
\newcommand{\Bb}{{\mathcal B}}

\newcommand{\Hhh}{{\mathfrak H}}

\begin{document}

\author[K. Corlette]{Kevin Corlette}
\address{University of Chicago
\\ 5734 University Avenue\\
Chicago, IL 60637, USA}
\email{kevin@math.uchicago.edu}
\urladdr{http://math.uchicago.edu/$\sim$kevin/}

\author[C. Simpson]{Carlos Simpson}
\address{CNRS, Laboratoire J. A. Dieudonn\'e, UMR 6621
\\ Universit\'e de Nice-Sophia Antipolis\\
06108 Nice, Cedex 2, France}
\email{carlos@math.unice.fr}
\urladdr{http://math.unice.fr/$\sim$carlos/} 

\title[Rank two representations]{On the classification of rank two representations of quasiprojective
fundamental groups}

\begin{abstract}
Suppose $X$ is a smooth quasiprojective variety over $\cc$ and $\rho : \pi _1(X,x) \rightarrow SL(2,\cc )$ 
is a Zariski-dense representation with quasiunipotent monodromy at infinity. Then $\rho$ factors through
a map $X\rightarrow Y$ with $Y$ either a DM-curve or a Shimura modular stack.  
\end{abstract}


\maketitle


\section{Introduction} \label{sec-introduction}

Let $X$ be a connected 
smooth quasiprojective variety over $\cc$ with basepoint $x$. We look at 
representations $\rho : \pi _1(X,x)\rightarrow
SL(2,{\mathbb C})$. We assume throughout that {\em the monodromy at infinity is quasi-unipotent}.
If $X\subset \overline{X}$ is a normal-crossings compactification with complementary divisor 
$D=\sum D_i$, and if $\gamma _i$ are loops going around the components $D_i$, this condition 
means that the $\rho (\gamma _i)$ are quasi-unipotent, in other words their eigenvalues are roots of unity.

A representation $\rho$ is  {\em Zariski-dense} if the Zariski-closure of its image is the whole group $SL(2,\cc )$.
A reductive representation of rank two is either Zariski-dense, or else it becomes  reducible upon pullback to a 
finite unramified covering of $X$.  We will classify representations
$\rho$ which are Zariski-dense and have quasi-unipotent monodromy at infinity. See \cite{Beauville} \cite{Arapura}
\cite{Delzant2} \cite{Dimca} \cite{FalkYuzvinsky} for a similar classification
in the reducible case. 

The geometry of the fundamental group of an algebraic variety has been studied from many
different angles \cite{AllcockCarlsonToledo}
\cite{AmorosBurgerEtAl} \cite{Arapura} \cite{ArapuraBresslerRamachandran} \cite{Campana94}
\cite{Catanese} \cite{Donaldson} \cite{Gromov} \cite{JostYau} \cite{KapovichMillson} \cite{Kollar}
\cite{Panov} \cite{Reznikov} \cite{RobbTeicher} \cite{Toledo} \cite{ViehwegZuo}. 
The methods we will use here are based on the theory of harmonic mappings, both to symmetric spaces and 
combinatorial complexes
\cite{canonical} \cite{rigid} \cite{Delzant} \cite{DonaldsonApp} \cite{Eyssidieux} 
\cite{GromovSchoen}  \cite{JostZuo1} \cite{JostZuo2} \cite{Katzarkov} \cite{KatzarkovRamachandran} \cite{Klingler}
\cite{TMochizuki2} \cite{NapierRamachandran3} 
\cite{hbls} \cite{ubiquity} \cite{Siu} \cite{Zuo}.

Our classification is obtained by looking at the interplay between different properties of $\rho$.
The main property is {\em factorization}: we say that $\rho$ {\em factors through a map}
$f:X\rightarrow Y$ if it is isomorphic to the pullback of a representation of $\pi _1(Y,f(x))$.
This notion can be extended in a couple of ways, for example $\rho$ {\em projectively factors}
through $f$ if the projected representation into $PSL(2,\cc )$ factors through $f$. The other extension is
that it is convenient (and basically --almost essential---to look at the notion of factorization through maps to
{\em Deligne-Mumford stacks} $Y$ rather than just varieties. In a certain sense this takes the place of 
complicated statements involving coverings of $X$. It even subsumes the notion of projective factorization,
because projective factorization is equivalent to factorization through a new
DM-stack obtained by putting a stack structure with group $\zz / 2$ (the center of $SL(2,\cc )$)
over the generic point of $Y$. 

One of the main cases of factorization we shall be concerned with is factorization through a curve.
A smooth one-dimensional DM-stack will  be called a {\em DM-curve}. 
Recall that an {\em orbicurve} is a DM-curve whose generic stabilizer is trivial. An orbicurve
is given by the data of a smooth curve together with a collection of marked points assigned integer
($\geq 2$) weights. Factorization through a DM-curve is equivalent to projective factorization through
an orbicurve (Corollary \ref{factorequivalence}). 

The other case we need to consider arises when the representation is motivic, in fact comes from a family of 
abelian varieties. The families of abelian varieties whose monodromy representations break up into rank two pieces
are given by maps to certain Shimura varieties or stacks. These Shimura varieties are closely analogous to Hilbert 
modular varieties. However, Hilbert modular varieties parametrize abelian varieties with real multiplication, while
in general we need to look at abelian varieties with multiplication by a totally imaginary extension of a totally 
real field. The condition that the tautological representation goes into $SL(2)$ basically says that the universal
covering of the Shimura variety is a polydisk i.e. a product of one-dimensional disks. We work without level structure
and call these things
{\em polydisk Shimura modular DM-stacks}. A classical example is the case of Shimura curves. The precise construction will be 
reviewed in \S \ref{sec-shimura} below. If $H$ is a polydisk Shimura DM-stack then $\pi _1(H)$ has a tautological
representation into $SL(2,L)$ for a totally imaginary extension $L$ of a totally real field, and this gives a
collection of tautological representations into $SL(2,\cc )$ indexed by the embeddings $\sigma : L\rightarrow \cc$. 

Classifying our rank two Zariski-dense representations
quasi-unipotent at infinity, will consist then of showing that any such representation factors through a map
$f:X\rightarrow Y$, with $Y$ being either a DM-curve, or else a polydisk Shimura modular DM-stack. We consider as
``known'' the representations on these target stacks $Y$.
There may be some overlap between these two cases, but one of our basic tasks is to have properties
which determine which case of the classification we will want to prove for a given
representation. 

Since we are looking at representations on quasiprojective varieties, we define {\em rigidity} in a way
which takes into account the monodromy at infinity. Fix a normal crossings compactification of 
$X$. For each component $D_i$ of the divisor at infinity, we have a well-defined conjugacy class of elements of
$\pi _1(X,x)$ corresponding to a loop $\gamma _i$ going around that component. Thus for a given representation
$\rho$ this gives a conjugacy class $C_i$ in the target group. We are assuming that these monodromy elements
are quasi-unipotent, so $C_i$ is a quasi-unipotent conjugacy class. We can define an affine variety  
$R(X,x,SL(2), \{ \overline{C}_i \} )$ of representations such that the monodromies $\rho (\gamma _i)$ are contained
in the closures of the $C_i$. Let $M(X,SL(2), \{ \overline{C}_i \} )$ denote its universal categorical
quotient by the conjugation action. We say that $\rho$ is {\em rigid} if it represents an isolated point
in the moduli space $M(X,SL(2), \{ \overline{C}_i \} )$ obtained by looking at its own conjugacy classes. 
In the case of a Zariski-dense representation, this form of rigidity means that there is no
non-isotrivial family of representations all having the same conjugacy classes at infinity, going through
$\rho$ (Lemma \ref{rigidequiv}). 

A property which plays a similar role but which is easier to state is {\em integrality}. Say that a representation
$\rho$ is {\em integral} if it is conjugate, in $SL(2,\cc )$, to a representation 
$\rho : \pi _1(X,x) \rightarrow SL(2,A)$ for $A$ a ring of algebraic integers. For Zariski-dense
representations, this is equivalent to 
asking that the traces $Tr (\rho (\gamma ))$ be algebraic integers for all $\gamma \in \pi _1(X,x)$. 

Say that {\em $\rho$ comes from a complex variation of Hodge structure} 
if there is a structure of complex variation of Hodge structure on the corresponding
local system $V$.

The main relationship between all of these notions is the following first result. 

\begin{global-theorem}
\label{full-global}
Suppose $\rho : \pi _1(X,x) \rightarrow SL(2,\cc )$ is a representation with quasi-unipotent monodromy at infinity,
such that $\rho$ does not
projectively factor through an orbicurve, or equivalently $\rho$ doesn't factor through a map to a DM-curve. 
Then $\rho$ is rigid and integral. Rigidity implies that $\rho$ comes from 
a complex variation of Hodge structure. 
\end{global-theorem}

This is already known in the case when $X$ is projective 
from \cite{ubiquity} for rigidity, and for integrality Gromov-Schoen
\cite{GromovSchoen}, and \cite{lefschetz}, the latter of which was
designed to support the original dormant version of this paper. The variation of Hodge structure
follows from \cite{canonical}, see \cite{hbls}. In the present, 
we extend the result to the quasi-projective case
for representations with quasi-unipotent monodromy at infinity. 
The various statements in Theorem 1 appear as
Theorems \ref{nonfactorsrigid}, 
\ref{nonfactorsintegral} and
\ref{rigidvhs} below. 

The underlying argument for both rigidity and integrality comes from Theorem 
\ref{mainfactorization}
about harmonic maps to Bruhat-Tits trees. This strategy is perhaps worth commenting on. It has its origins in the
work of Bass and Serre \cite{Bass} \cite{Serre} \cite{Serre2}, Culler-Shalen \cite{CullerShalen} and 
Gromov-Schoen \cite{GromovSchoen}. 

It would certainly have
been possible to treat the rigidity question using harmonic maps to symmetric spaces \cite{canonical} 
\cite{DiederichOhsawa} \cite{DonaldsonApp} \cite{EellsSampson} 
\cite{TMochizuki2}. 
For integrality, though, it is
necessary to use the theory of  harmonic maps to Bruhat-Tits trees \cite{GromovSchoen}. 
Furthermore, there is a sort of analogy between the two notions: 
integrality means that a representation into $SL(2,\qq _p)$ goes
into a compact subgroup, whereas rigidity may be thought of as saying that a representation into $SL(2,\cc (t))$ goes
into a compact subgroup, much as in \cite{CullerShalen}. 
So, we thought it would be interesting to use harmonic maps to trees to treat both cases at once.
This reduces the volume of material about  harmonic maps. On the other hand it introduces an additional difficulty,
because the Bruhat-Tits tree for $SL(2,\cc (t))$ is not locally compact. We overcome this by making a reduction to the case
of representations in $SL(2, \ff _q (t))$ where $\ff _q$ is a finite field. This reduction is fairly standard but
it requires a finiteness theorem for the number of possible maps to a hyperbolic DM-curve
(Proposition \ref{finiteDMfactors}), see Theorem \ref{nonfactorsrigid}. Delzant appears to have independently found
this proof of the fact that nonfactorization implies rigidity. Some pieces of his proof, including the finiteness statement, 
appear in \cite{Delzant2}
but that paper is oriented towards proving a similar statement for representations in a solvable group. 
The classification for solvable representations in \cite{Beauville} and \cite{Delzant2} complements the present paper because we restrict here 
to Zariski dense representations, which are irreducible over any finite covering.

Suppose now that $\rho$ does not factor through a curve. From the above results we obtain
that $\rho$ is rigid, integral, and comes from a complex variation of Hodge structure. 
Let $A$ be a ring of algebraic integers such that $\rho$ is defined over $A$.
Let $\sigma : A\rightarrow {\mathbb C}$ be an embedding (not necessarily the identity one).
Let $\rho ^{\sigma}$ be the composed representation
$$
\pi _1(X,x)\stackrel{\rho _A}{\longrightarrow}
SL(2,A)\stackrel{\sigma}{\longrightarrow}
SL(2,{\mathbb C}).
$$
Then $\rho ^{\sigma}$ is rigid too. Hence $\rho ^{\sigma}$ comes from a complex variation of Hodge structure
for every $\sigma$. We will use this data to construct a factorization of $\rho$ through a polydisk Shimura modular DM-stack $H$.
As before this means that there is a map $f:X\rightarrow H$ such that
$\rho$ is the pullback $f^{\ast}$ of one of the tautological rank two local systems on $H$. 
The main idea is that because the rank is $2$, the Hodge types can be chosen to be $(1,0)$ and $(0,1)$; then, by integrality we get
a family of abelian varieties. 
The construction of the map to $H$ is straightforward but it has to take into account the notion of polarization
and the special structure of abelian varieties whose Hodge structures split into rank two pieces over a totally imaginary field. 
This gives the statement of our main classification result, see \S \ref{sec-classif}:

\begin{global-theorem}
\label{classification-global}
Suppose $X$ is a smooth quasiprojective variety
and $\rho : \pi _1(X,x) \rightarrow SL(2,\cc )$ is a Zariski-dense representation. 
Suppose that the monodromy transformations around components of the divisor at infinity are quasi-unipotent.
Then either $\rho$ comes from a map $f:X\rightarrow Y$ to a DM-curve $Y$, or else $\rho$ comes from 
pullback of one of the tautological representations by a map $f:X\rightarrow H$ to a polydisk Shimura modular DM-stack $H$.
\end{global-theorem}

The two cases described in this theorem can overlap: there can be rigid local systems on an orbicurve.
For any rank such things are classified by Katz's algorithm \cite{Katz}, but in the rank two case
they can be seen explicitly as hypergeometric systems \S \ref{sec-hypergeometric}. 

This paper is a project which we have been entertaining since around 1990. It was motivated by
Gromov's paper \cite{Gromov}, and spurred on by a lecture by R. Schoen in Chicago about 
\cite{GromovSchoen}. In the projective case
the paper \cite{lefschetz} was done in the context of this program,
in order to obtain the proof of Theorem \ref{full-global} and thus
Theorem \ref{classification-global}; in particular \cite{lefschetz} should be considered as
an integral part of the present project. 

An important element was M. Larsen's explanation of how to get information on the field of definition
of the representation. His argument, first reported in \cite[Lemma 4.8]{hbls}, plays a crucial role in 
starting off \S \ref{sec-improvements} below. 

As we have taken such a long time to write up the classification result,
the ambient state of technology has evolved in the meantime \cite{Biquard} \cite{JostZuo2} \cite{TMochizuki2} \cite{Panov} which makes it reasonable
to give statements for quasiprojective varieties. Another recent advance is that the notion of DM-stack has
become standard and well-understood, for example Behrend-Noohi \cite{BehrendNoohi} have classified DM-curves
in a way which is very useful for our considerations.

One of the reasons for getting back to this project is that recently
there have been some explicit constructions of rank two local systems on quasiprojective varieties,
which in some cases can be seen as coming from projective varieties by passing to a finite ramified cover.
The examples we know of are those of Boalch \cite{Boalch} and Panov \cite{Panov}. 
It would be interesting explicitly to determine the factorizations for these examples, but we don't treat that 
question here. That question seems to be answered in some cases by a recent paper of Ben Hamed and Gavrilov 
\cite{BenHamedGavrilov} which gives an explicit geometric origin for solutions of Painlev\'e VI 
equations. Kontsevich has a number of conjectures about local systems on curves \cite{Kontsevich}.

\section{Local systems on Deligne-Mumford stacks}
\label{sec-dms}

In order to obtain optimal statements of the type we shall consider throughout this paper,
it is convenient to consider the notion of local system over a Deligne-Mumford stack.

Recall that a {\em Deligne-Mumford stack} is a $1$-stack $X$ on the site of schemes over $\cc$
with the etale topology, such that there exists a surjective morphism of stacks $f:Z\rightarrow X$ from
a scheme $Z$ to $X$, with $f$ ``representable and etale''. These conditions mean that if $Y$ is any scheme
then $Z\times _XY\rightarrow Y$ is an etale morphism of schemes. We refer the reader to 
\cite{AbramovichVistoli}
\cite{BehrendNoohi} \cite{Campana04} \cite{DeligneMumford}
\cite{LaumonMB} \cite{Noohi}
for general references about this notion. Intuitively,  
a Deligne-Mumford stack is an object $X$ which looks like the
quotient of an algebraic variety $Z$ by an ``equivalence relation'' $R\rightarrow Z\times Z$
such that the projections from $R$ to the two factors are etale
(in particular quasi-finite). Technically speaking $R:= Z\times _XZ$ 
is also provided with a multiplication making it into the morphism-object for a groupoid whose
object-object is $Z$. In practice, in the cases we shall consider in this paper,
$X$ will often be the quotient of a smooth variety $Z$ by an action of a finite group $G$. 
This is the case when $X$ is a Shimura modular stack, and almost always the case when $X$ is 
an orbicurve or DM-curve. In the Shimura modular case, $X$  has many Galois coverings by Shimura modular varieties obtained
by imposing some level structure. In the one dimensional case, aside from a small number of degenerate situations,
an orbicurve has a Galois covering which is a smooth curve (see Lemmas \ref{nonspherical} and 
\ref{DM-curve-quotient} below). 

If $X$ is a DM-stack then we obtain its {\em coarse moduli space} $X^{\rm coarse}$ which is the universal algebraic space with a
map from $X$. This exists, and the map $X\rightarrow X^{\rm coarse}$ is finite, by Keel-Mori \cite{KeelMori}.

\begin{theorem}
\label{localquotient}
If $X$ is a DM-stack, there is a Zariski open covering of $X^{\rm coarse}$ such that the pullbacks of the open sets  
are quotient stacks by finite group actions. 
\end{theorem}  
{\em Proof:} This is stated in
\cite{ToenThesis} Proposition 1.17, with for proof a reference to \cite{Vistoli} 2.8 which should be taken in light of
\cite{KeelMori}. 
A sketch of proof is given in \cite{AbramovichVistoli} Lemma 2.2.3. 
\eop

Because of this theorem, the reader may without danger imagine that
the words ``DM-stack'' and so forth, basically mean varieties modulo finite group actions.  The {\em irreducible components}
of a DM-stack will by definition be the irreducible components of the coarse moduli space, which pull back to closed substacks
of $X$. 

Local properties of $X$ are defined
by requiring the same properties for the etale covering scheme $Z\rightarrow X$ occuring in the above definition. In particular,
a DM-stack $X$ is {\em smooth} if for one (or equivalently any) surjective etale map from a scheme $Z\rightarrow X$,
the scheme $Z$ is smooth. We restrict our attention in this paper to smooth DM-stacks. A {\em point} is a morphism
$x: Spec (\cc )\rightarrow X$, which we denote abusively by $x\in X$. A point can have ``automorphisms'',
which is the phenomenon new to stacks. 

The notion of Deligne-Mumford stack is a generalization and transfer to the algebraic category, of the classical 
notion of ``orbifold'' or $V$-manifold \cite{Kawasaki} \cite{Satake1} \cite{Satake2}. 
This notion adapts to the algebraic category:
an {\em orbifold} is a smooth DM-stack such that the general point of any irreducible component has
trivial automorphism group. The associated complex analytic stack in this case is exactly a complex analytic $V$-manifold. 
Then there is a dense Zariski open subset of the coarse moduli space $X^{\rm coarse}$ over which the projection is
an isomorphism, thus an orbifold has a dense open substack which is an algebraic space.  

If $X=Z/G$ is a quotient stack by the action of a finite group (with the quotient being irreducible, say), 
then $X$ is an orbifold if and only if $G$ acts faithfully. The coarse moduli space is the classical
quotient space of the action.  Orthogonal to this case is the important example
of the quotient of a single point by a trivial action of a finite group $G$: this is a DM-stack denoted $BG$. 
All of its points are isomorphic and they all have automorphism group isomorphic to $G$. 

These two examples permit us to form all possible smooth DM-stacks, as shown by the following structure result. 

\begin{proposition}
\label{gerb}
Suppose $X$ is a smooth DM-stack. There is a  universal orbifold $X^{\rm orb}$ (i.e. a smooth DM
stack with trivial generic stabilizers) with a map $\phi : X\rightarrow X^{\rm orb}$
such that locally in the etale topology over the base, the map $\phi$ is isomorphic to the projection of a product of the 
base orbifold with a $BG$ (this is usually called a {\em gerb}). 
\end{proposition}
{\em Proof:} (See \cite{BehrendNoohi}). Using Theorem \ref{localquotient}, over an open set where the stack has the form $Z/G$ for a smooth irreducible
variety $G$, let $H$ be the kernel of the action. It is a normal subgroup (since it is the kernel of the map from $G$ to the
automorphisms of $Z$). In this case, $X^{\rm orb}$ is the quotient stack of $Z$ by $G/H$, and $X$ is a gerb over $X^{\rm orb}$
with fiber $BH$ (use $Z$ itself for the etale neighborhood in question: we have $X \times _{X^{\rm orb}}Z \cong Z\times BH$). 
This satisfies a universal property so these local constructions glue together to give $X^{\rm orb}$ in the global case. 
\eop

\noindent
{\em Remark:} The coarse moduli spaces of $X$ and $X^{\rm orb}$ coincide. 

\subsection{Fundamental group and local systems}

If $X$ is a DM-stack and $x\in X$ is a point then we obtain the {\em fundamental group} $\pi _1(X,x)$. 
See Noohi \cite{Noohi} for
a general discussion.
In the quotient case $X = Z / G$, the fundamental group may be viewed more simply as an extension
$$
1\rightarrow \pi _1(Z,z) \rightarrow \pi _1(X,x) \rightarrow G \rightarrow 1.
$$
The group in the middle may be defined as the set of paths in $Z$ starting at $z$ and going to any
preimage of $x$. Paths are composed by first translating by an appropriate element of $g$, then juxtaposing paths.

We also have the notion of {\em local system} over $X$, which is by now a standard notion,  
see  \cite{LaszloOlsson} for example, or \cite{ToenThesis} for the case of $\Dd$-modules. 
The case of local systems on the moduli stack
of hyperelliptic curves was mentionned in a talk by R. Hain, see \cite[p. 12]{Hain}. A local system is
a collection $L_{Z,f}$ of local systems over schemes $Z$ for every section $f:Z\rightarrow X$, together with 
functoriality maps: if 
$$
f:Z\rightarrow X, \;\;\; f':Z'\rightarrow X
$$
are two maps, and if $g:Z'\rightarrow Z$ is a map together with a natural transformation $\eta$ from $f\circ g$ to $f'$ then
we get a map from $g^{\ast} (L_{Z,f})$ to $L_{Z',f'}$. These satisfy some natural axioms. We obtain a tensor category of
local systems over $X$. 

If $L$ is a local system of rank $r$ over $X$
and if $x\in X$ is a basepoint then the fiber $L_x$ is a $\cc$-vector space of dimension $r$, and the monodromy is an action
of $\pi _1(X,x)$ on it. If $X$ is smooth and irreducible with fixed basepoint $x$ then the category of local systems over $X$
is equivalent to the category of representations of $\pi _1(X,x)$. 

If $X=Z/G$ is a quotient stack of a variety $Z$, 
a local system on $X$ may be seen as a $G$-equivariant
local system on $Z$, i.e. a pair $(V_Z, \alpha )$  where $V_Z$ is a local system on $Z$ and
$\alpha$ is an action associating to each $g\in G$ an isomorphism $\alpha (g) : g^{\ast }V_Z\cong V_Z$.
The action is required to satisfy a cocycle condition: if $g,h\in G$ then $\alpha (g) \circ g^{\ast}
\alpha (h) = \alpha (gh)$. 

The notion of local system on a DM stack 
beams back down to the world of varieties in the following way: if $Y$ is a variety
and $f:Y\rightarrow X$ is a map to a Deligne-Mumford stack $X$ then for any local system $V$ on $X$
the pullback $f^{\ast}V$ is a local system on $Y$. 

Something new here is that a map $f$ can have 
{\em automorphisms}. Roughly speaking an automorphism of $f$ is a section of the sheaf of stabilizer groups
pulled back over $f$. More generally we should speak of {\em isomorphisms} between maps $f,g:Y\rightarrow X$.
We have the following functoriality: if $a:f\Rightarrow g$ is an isomorphism from $f$ to $g$ then
we obtain an isomorphism of pullback local systems 
$a^{\ast}V : f^{\ast} V\cong g^{\ast} V$. This satisfies some usual functoriality and associativity
identities. This phenomenon will appear in our statements when we want to say that a map is unique.

\subsection{DM-curves and orbicurves}

We are particularly interested in the case of objects of dimension $1$. A {\em DM-curve} is a smooth DM-stack of dimension $1$.
An {\em orbicurve} is an orbifold of dimension $1$, so an orbifold is a DM-curve with trivial generic stabilizer.
In this case Proposition \ref{gerb} says that a DM-curve is always a gerb over a canonical orbicurve. 
If $X$ is an orbicurve then its coarse moduli space $X^{\rm coarse}$ is again a smooth curve (this is special to the case of
orbicurves: in general the coarse moduli space of a smooth orbifold will have finite quotient singularities). 

An etale covering is 
a finite etale map; over the complex numbers this is the same thing as a finite topological
covering space. For DM-curves, etaleness is measured in terms of a local chart for the stack. 
In practical terms this comes down to saying that the ramification index, taking correctly into account the
orbifold indices, is $1$.

The data of an orbicurve is determined by the smooth curve $X^{\rm coarse}$, together with a finite set of points $P_j$
and an integer $n_j\geq 2$ attached to
each point. A nonconstant morphism from a connected smooth curve to $X$ is any nonconstant map 
$f:Z\rightarrow X^{\rm coarse}$ such that
for any $z\in Z$ lying over some $P_j$ the ramification order of $f$ at $z$ is divisible by $n_j$. The map is etale
when the ramification at $P_j$ is exactly $n_j$. 

We denote an orbicurve by the notation $(Y,n_1, \ldots , n_k)$ with $Y = X^{\rm coarse}$ a smooth curve and 
$n_i$ the sequence of integers organized in decreasing
order (left out of this notation is the choice of points $P_1,\ldots , P_k\in X$).

The paper of Behrend and Noohi \cite{BehrendNoohi} is a very complete description of the possibilities for DM-curves and orbicurves,
and the reader is referred there. Recall  here the general outlines of their classification which generalizes the classification of
curves. An orbicurve is either {\em spherical}, {\em elliptic}, or {\em hyperbolic}.  The spherical orbicurves are
$\pp ^1$, the {\em drops} $(\pp ^1, a)$ and {\em footballs} $(\pp ^1, a,b)$ as well as the finite list of cases
$$
(\pp ^1, 2, 2, 2), \;\; 
(\pp ^1, 2, 3, 3), \;\; 
(\pp ^1, 2, 3, 4), \;\; 
(\pp ^1, 2, 3, 5).
$$
The universal coverings of these orbicurves are either $\pp ^1$, drops, or footballs with relatively prime indices.

\begin{lemma}
\label{nonspherical}
An orbicurve which is not spherical as listed above, has an etale covering which is a regular curve different from $\pp ^ 1$. 
In particular the
associated complex-analytic orbicurve has a contractible universal covering. 
\end{lemma}
{\em Proof:} See \cite{BehrendNoohi}.
\eop

Now we can define the {\em elliptic} orbicurves to be those which have a covering by either ${\mathbb A}^ 1$ or 
${\mathbb G}_m$ or an elliptic curve (thus, they are those whose universal covering is $\cc$); and the 
{\em hyperbolic} orbicurves are the remaining ones---those whose universal covering is a disc. 

Note also that an elliptic or hyperbolic orbicurve has infinite fundamental group except for the elliptic case
where the covering is ${\mathbb A}^ 1$. In all elliptic and hyperbolic cases, the homotopy type is a 
$K(\pi _1, 1)$.

We say that a DM-curve $X$ is spherical, elliptic or hyperbolic according to the type of the associated orbicurve
$X^{\rm orb}$.

\begin{lemma}
\label{zdhyperbolic}
Suppose $X$ is a DM-curve and $\rho : \pi _1(X,x) \rightarrow SL(2,K)$ is a Zariski-dense representation for any
infinite field $K$. Then $X$ is hyperbolic.
\end{lemma}
{\em Proof:}
Spherical orbicurves have finite fundamental groups, and 
elliptic orbicurves have fundamental groups which are virtually abelian. 
Thus the same holds for spherical and elliptic DM-curves.
Neither of these types of groups can have Zariski-dense representations to $SL(2,K)$
for an infinite field $K$.
\eop

In view of this lemma, for factorization questions we will mainly be looking at hyperbolic orbicurves.

\begin{lemma}
\label{DM-structure}
Suppose $X$ is a DM-curve with basepoint $x$. 
Suppose the underlying orbifold $Y:=X^{\rm orb}$ is either elliptic or hyperbolic. Then
$\pi _2(Y) = 0$ and the gerb $p:X \rightarrow Y$ (locally a product with $BG$) 
is completely determined by the induced extension
$$
1 \rightarrow G \rightarrow \pi _1(X, x) \rightarrow \pi _1(Y, p (x)) \rightarrow 1.
$$
Thus we may think of $X$ as being given by the data of an orbicurve $Y$ plus an extension of 
the fundamental group by a group $G$. If $Z$ is a scheme then a submersive map from $Z$ to $X$ is the
same thing as a submersive map from $Z$ to $Y$ plus a lifting of the map on fundamental groups from
$\pi _1(Z,z)$ to $\pi _1(X,x)$. 
\end{lemma}
{\em Proof:} This is basically the same as Behrend-Noohi \cite{BehrendNoohi} Proposition 4.7.
 
If $X$ is already an orbicurve then
the generic stabilizer group $G$ is trivial, $BG = {\rm Spec} (\cc )$ and $Y=X$.

In general, let $Z$ be the coarse moduli space for $X$. It is a normal (hence smooth) curve.
Let $Y$ be the orbicurve obtained by setting, for each point $z\in Z$, the orbifold structure at $z$
to be the ramification degree of the map $X\rightarrow Z$ over $z$ (this is different from $1$ for only
a finite number of points). There is a unique factorization to a map $p:X\rightarrow Y$ (note that there
is no question of natural transformations here because $p$ is submersive and
the generic stabilizers on $Y$ are trivial). The map $p$ is etale by the choice of
orbifold indices. An etale covering is a fibration in the etale topology, and the fiber over a
point (say a general point) is a one-point DM-stack, hence of the form $BG$ for a finite group $G$.
We obtain a long exact sequence in homotopy. If $X$ is elliptic or hyperbolic, then by definition the same
is true of $Y$. Thus  
the universal covering of $Y$ is $\cc$ or a disk \cite{BehrendNoohi}, so $\pi _2(Y)=0$.
In particular, the long exact homotopy sequence of homotopy groups gives an extension as in the statement of the lemma.
The fact that this extension uniquely determines the fibration $p:X\rightarrow Y$, plus the statement
about maps to $X$, will be left as
exercises in nonabelian cohomology \cite{Breen}. 
\eop

We will see below that the notion of DM-curve gives a handy way of dealing with the difference between
$SL(2)$ and $PSL(2)$ for factorization statements. In particular, since we are mainly using the notion
of DM-curve for that purpose, it mostly suffices to consider the case where the group $G$ in Lemma
\ref{DM-structure} is of order two. 

Putting together the above information we obtain the following. 

\begin{lemma}
\label{DM-curve-quotient}
If $X$ is an elliptic or hyperbolic DM-curve, then there is an etale Galois covering
$Z\rightarrow X$ with $Z$ a smooth curve; in particular $X$ is a quotient of $Z$ by the 
Galois group $G$---note that the action of $G$ on $Z$ might not be faithful; this is how we get
nontrivial generic stabilizers. As described above we have an extension
$$
1\rightarrow \pi _1(Z) \rightarrow \pi _1(X)\rightarrow G \rightarrow 1.
$$
\end{lemma}
{\em Proof:}
Use Lemma \ref{DM-structure} to express $X$ as a gerb over an orbicurve $Y$. From the  long exact
sequence we see that $Y$ has infinite fundamental group, so we are not in the degenerate cases 
of Lemma \ref{nonspherical}. Thus there is an etale covering map $V\rightarrow Y$ from a smooth curve to $Y$.
Pulling back, we get a gerb $U:=X\times _YV$ over $V$. The gerb corresponds to an extension
$$
1\rightarrow H\rightarrow \pi _1(U) \rightarrow \pi _1(V)\rightarrow 1.
$$
Let $U'\rightarrow U$ denote the covering corresponding to the map $\pi _1(U)\rightarrow {\rm Aut(H)}$.
This is again a gerb corresponding to a central extension
$$
1\rightarrow Z(H)\rightarrow \pi _1(U') \rightarrow \pi _1(V')\rightarrow 1.
$$
The central extension is classified by an element of $H^2(V',Z(H))$. There is a finite etale covering space
$V''\rightarrow V'$ such that the pullback of the element is zero, that is the extension splits.
This gives an etale covering morphism from $V''$ back down to $X$. Galois completion yields the curve $Z$.
\eop

Due to this lemma, and the fact that the other stacks we consider (Shimura modular stacks) are by
construction analytic quotients and are seen to be algebraic quotients by looking at moduli spaces for
objects with level structure, it is possible to restrict our attention to Deligne-Mumford stacks which are quotients.

\subsection{Finiteness}

We prove a finiteness statement which is the key to our method of proof of the implication
``nonfactoring implies rigid'' in \S \ref{sec-rigidity} below. It has also independantly appeared
in Delzant's preprint \cite{Delzant2} with a more conceptual proof. 

\begin{lemma}
\label{indexbound}
Given $g$ and $b$, there is an integer $N$ giving the following bound. For any 
hyperbolic orbicurve $X$, and any regular curve $Y$ with a nonconstant map $f:Y\rightarrow X$,
assume that the projective compactification
$\overline{Y}$ of $Y$ has genus $\leq g$, and the number of points in $\overline{Y}-Y$ 
is $\leq b$. Then for any point $P$ in the image of $f$, the orbifold index of $X$ at $P$ is $\leq N$.
\end{lemma}
{\em Proof:}
We may assume that $X$ is compact. Indeed, if not then we could add in orbifold structures at the points
at infinity of $\overline{X} - X^{\rm coarse}$. Note that we are not asking in the hypothesis that the
map $f$ be surjective, only that the point $P$ in question be in the image. If we fix an orbifold index
of $\geq 7$ at these additional points then the new orbifold will also be hyperbolic
(see \cite{BehrendNoohi}), so we may assume $X$ is compact. 

Fix $f:Y\rightarrow X$, and points $Q\in Y$ and $P=f(Q)\in X$. Let $n$ be the orbifold index of
$X$ at $P$. \revision Denote by $\overline{f}$ the extension of $f$ to $\overline{Y} \rightarrow \overline{X}$
(where $\overline{X}$ denotes the usual curve which is the projective completion of $X^{\rm coarse}$).

The ramification degree of $f$ at $Q$ is divisible by $n$,
so the contribution to the Hurwitz formula from $Q$ is at least $n-1$. If $g(\overline{X}) \geq 1$ we have 
$$
2g(\overline{Y}) - 2 \geq d (2g(\overline{X} - 2) + (n-1) \geq n-1
$$
which gives the bound in question. 

The remaining problem is to treat the case $\overline{X} = \pp ^1$.  Since $X$ is hyperbolic, there are at least
two other orbifold points $P_1$ and $P_2$ besides $P_0:=P$. Let $n_1$ and $n_2$ denote their orbifold indices.
Let $d$ denote the degree of the map $\overline{f}$. If $Q'$ is a point in $\overline{f}^{-1}(P_i)$
then either $Q'\in Y$, in which case the ramification degree $r_{Q'}$ is divisible by $n_i$ (if $i=0$ we are
denoting $n_0:=n$); or else
$Q' \in \overline{Y}-Y$. There are at most $d/n_i$ points in the first case, and at most $b$ points in the last case. 

The Hurwitz contribution for the fiber over $P_i$ is bounded:
$$
\sum _{\overline{f}(Q')=P_i} (r_{Q'}-1) = d - \# (\overline{f}^{-1}(P_i)) 
\geq d - b - d/n_i .
$$
Therefore the Hurwitz formula gives
$$
2g(\overline{Y}) - 2 \geq -2d + 3d -3b - d (\frac{1}{n_0} + \frac{1}{n_1} + \frac{1}{n_2}).
$$
For $n_0 = n$ this gives a bound 
$$
\frac{1}{n} \geq \frac{2 - 2g(\overline{Y}) -3b}{d} + (1 - \frac{1}{n_1} - \frac{1}{n_2}).
$$

In the case that at least one of $n_1$ or $n_2$ are $\geq 3$ we obtain 
$$
(1 - \frac{1}{n_1} - \frac{1}{n_2}) \geq \frac{1}{6}
$$
so 
$$
\frac{1}{n} \geq \frac{2 - 2g(\overline{Y}) -3b}{d} + \frac{1}{6}.
$$
If $n \leq 6$ then we  have our bound so we can assume that $n\geq 7$. Thus
$$
\frac{3b + 2g(\overline{Y}) - 2}{d} \geq \frac{1}{42},
$$
or 
$$
d \leq 42 (3b + 2g(\overline{Y}) - 2).
$$
This gives a bound for $d$ hence for $n\leq d$.

Assume now that $n_1=n_2=2$. In this case there must be a fourth orbifold point, because the orbifolds with
indices $(2,2,a)$ are spherical for any $a$. Adding in the fourth point to the above calculations the bound
becomes 
$$
\frac{1}{n} \geq \frac{2 - 2g(\overline{Y}) -4b}{d} + (2 - \frac{1}{n_1} - \frac{1}{n_2}- \frac{1}{n_3}).
$$
Here 
$$
(2 - \frac{1}{n_1} - \frac{1}{n_2}- \frac{1}{n_3}) \geq \frac{1}{2}
$$
so
$$
\frac{1}{n} \geq \frac{2 - 2g(\overline{Y}) -4b}{d} + \frac{1}{2}.
$$
We can assume that $n\geq 3$,  thus
$$
\frac{4b + 2g(\overline{Y}) - 2}{d} \geq \frac{1}{6},
$$
or 
$$
d \leq 6 (4b + 2g(\overline{Y}) - 2).
$$
Again this gives a bound for $n\leq d$. 
\eop

We now get our finiteness result, cf \cite[Theorem 2]{Delzant2}. 

\begin{proposition}
\label{finiteDMfactors}
Suppose $X$ is a connected smooth quasiprojective variety. The set of isomorphism classes of pairs 
$(Y,f)$ where $Y$ is an hyperbolic orbicurve, and $f:X\rightarrow Y$ is
a surjective algebraic morphism, is finite.
\end{proposition}
{\em Proof:}
Let $\overline{X}$ be a compactification such that the complementary divisor has normal crossings. 
We can find a projective embedding of $\overline{X}$ and choose a connected family of complete intersection
curves in $\overline{X}$ passing through any point. We may assume that through any point there is a 
smooth curve in the family. Let $g$ be the genus of the curves in the family, and let $b$ be the maximal
number of intersection points of a curve with the divisor $D:=\overline{X} - X$. We may also assume that for
any point $Q\in X$ there is a curve $\overline{C}$ in the family passing through $Q$ and intersecting the divisor $D$
transversally. Let $C:= \overline{C}\cap X$, and note that $\pi _1(C)\rightarrow \pi _1(X)$ is surjective. 

Suppose $f:X\rightarrow Y$ is a surjective morphism to a hyperbolic orbicurve, and suppose $P$ is
a point in the base. Choose a point $Q\in X$ with $f(Q)=P$ and choose a curve $C$ in the above family
passing through $Q$. 
The map $f|_C$ is nonconstant since the map on fundamental groups is nontrivial because
$\pi _1(C)\rightarrow \pi _1(X)$ is surjective and 
$\pi _1(X)\rightarrow \pi _1(Y)$ has image of finite index in an infinite group. 
The point $P$ is in the image of $f|_C$, so Lemma {indexbound} provides a bound for the orbifold index of
$Y$ at $P$ in terms of the bounds $g$ and $b$ for the family of curves $C$. Thus we obtain a bound $N$
depending only on $X$, for the orbifold indices of points in $Y$. Assume $N\geq 7$. 

Note also that the genus of $\overline{Y}$ is bounded by the genus $g$ of the curves $\overline{C}$ in the family.
The number of points at infinity (i.e. the cardinality of $\overline{Y}-Y^{\rm coarse}$) is bounded by
the bound $b$ for the curves $C$. Finally, an argument similar to that of Lemma \ref{indexbound}
provides a bound for the number of orbifold points in $Y$. Thus there are only finitely many possible
combinatorial types for the target orbicurve $Y$. 

For the remainder of the argument, we can fix a combinatorial type for a hyperbolic orbicurve $Y$.
There is an etale Galois covering $Z\rightarrow Y$ such that $Z$ is a curve (without orbifold structure).
We may also assume that the projective completion $\overline{Z}$ has genus $g_Z\geq 2$. Let $G$ denote the
Galois group of $Z/Y$. The data of $g_Z$ and the group may be determined by the combinatorial type of $Y$
up to finitely many choices, so we can view them as being fixed.  

Now if $f:X\rightarrow Y$ is any map, the pullback $X\times _YZ$ is a
Galois covering of $X$ with group $G$. There are only finitely many possibilities for this covering,
so we may view it as being a fixed covering $V\rightarrow X$ with Galois action of $G$. 
The map $f$ corresponds to a $G$-equivariant map from $V$ to $Z\subset \overline{Z}$ where
$\overline{Z}$ is a variable Riemann surface of genus $g_Z$. Once the map $V\rightarrow \overline{Z}$
is fixed, the image open set $Z$ and the action of $G$ on $Z$ are determined. 

We are now reduced to a classical statement  \cite{KobayashiOchiai}, \cite[p. 137]{Lang}: 
given a quasiprojective variety $V$ there are only finitely
many isomorphism classes of submersive maps to a compact Riemann surface of genus $g_Z\geq 2$. 
\eop

\section{Factorization}
\label{sec-factor}

Superrigidity theory says that under certain circumstances, varieties with the same topology have to be
isomorphic. The notion of factorization, which extends this idea to maps which are not isomorphisms,
is by now classical too 
\cite{Beauville} \cite{Campana94} \cite{Catanese} \cite{rigid} \cite{Delzant} \cite{Eyssidieux}
\cite{Gromov} \cite{GromovSchoen}
\cite{Katzarkov} \cite{KatzarkovRamachandran} \cite{Klingler} \cite{Kollar} 
\cite{ubiquity} \cite{Siu} \cite{Zuo}.

Suppose $X$ is a smooth quasiprojective variety and
$\rho :\pi _1(X,x) \rightarrow G$ is a representation into some group $G$. A {\em
factorization} is a triple $(Y,f,\psi )$ where $Y$ is a DM-stack, $f:X\rightarrow Y$ is an algebraic map,
and $\psi : \pi _1(Y,f(x))\rightarrow G$ is a representation such that 
$\rho = \psi \circ \pi _1(f)$. In this case we say that {\em $\rho$ factors through $f:X\rightarrow Y$}
or just {\em through $Y$}. If $\pi _1(f)$ is not surjective then
the representation $\psi $ might not be unique. 

We say that {\em $\rho$ projectively factors through $f:X\rightarrow Y$} if the composed representation
$$
\tilde{\rho} : \pi _1(X)\rightarrow G/Z(G)
$$
to the quotient of $G$ by its center, factors through $f$. This terminology is intended for the case
$G=SL(2,K)$ and $G/Z(G)=PSL(2,K)$.

For example we say that a representation $\rho$ {\em factors through a DM-curve} if there
exists a map $f:X\rightarrow Y$ to a smooth Deligne-Mumford curve and a factorization of $\rho$ through $f$.

Recall that an  orbicurve is a DM-curve such that the generic stabilizer is trivial.
The following lemmas say that it is equivalent to speak  of projective factorization through
maps to orbicurves, or actual factorization through maps to DM-curves. In the remainder of the paper
we will use the language of factorization through DM-curves since this is easier to handle notationally
speaking. 

\begin{lemma}
\label{projfactorsDMfactors}
Suppose $X$ is a smooth quasiprojective variety and $\rho : \pi _1(X,x)\rightarrow SL(2,K)$ a representation
with $K$ a field. Suppose $\rho$ projectively factors through 
a map $f: X\rightarrow Y$ to an orbicurve, i.e. there is a factorization $(Y,f,\psi )$ for
the representation $\tilde{\rho}$ into $PSL(2,K)$. 
Then there is a structure of DM-curve $Y'$ mapping to $Y$ and a lifting 
$f':X\rightarrow Y'$ such that $\rho$ factors through $f'$. 
\end{lemma} 
{\em Proof:}
Recall from \ref{zdhyperbolic} that $Y$ must be hyperbolic.
Use $\psi $ to pull back the central extension
$$
1 \rightarrow \{ \pm 1 \} \rightarrow SL(2,K) \rightarrow PSL(2,K) \rightarrow 1.
$$
We obtain a central extension of $\pi _1(Y)$ by a group of order $2$. By Lemma \ref{DM-structure} this corresponds to a
structure of DM-curve $Y'$ over $Y$ and the representation $\psi$ lifts to a representation
$\psi ' : \pi _1(Y')\rightarrow SL(2,K)$. Note that the square 
$$
\begin{array}{ccc}
\pi _1(Y') & \rightarrow & \pi _1(Y) \\
\downarrow & & \downarrow \\
SL(2,K) & \rightarrow & PSL(2,K)
\end{array}
$$
is cartesian. The projective factorization gives a map into the fiber product of the lower right angle of this diagram,
thus it gives a map $\pi _1(X)\rightarrow \pi _1(Y')$. Thus by \ref{DM-structure} 
the map $f:X\rightarrow Y$ lifts to a map $f':X\rightarrow Y'$
giving a factorization of the representation $\rho$. 
\eop

In the above lemma, the structure of DM curve isn't unique. 

\begin{lemma}
\label{DMfactorsprojfactors}
Suppose $X$ is a smooth quasiprojective variety and $\rho : \pi _1(X,x)\rightarrow SL(2,K)$ is a Zariski-dense representation
with $K$ a field. 
If $\rho$ factors through a map
to a DM-curve $f:X\rightarrow Y'$ then letting $Y$ be the underlying orbicurve as defined above, $\rho$ projectively
factors through $f:X\rightarrow Y$. 
\end{lemma} 
{\em Proof:}
Suppose $(Y',f',\psi ')$ is a factorization for $\rho$. Note that $Y'$ is hyperbolic so it is classified by the
extension 
$$
1\rightarrow G \rightarrow \pi _1(Y')\rightarrow \pi _1(Y) \rightarrow 1.
$$
The group $G$ is finite. The subgroup of $SL(2,K)$ which normalizes the image $\psi ' (G)$
is an algebraic subgroup containing the image of $\rho$. By the Zariski-density hypothesis, this group is the whole
of $SL(2,K)$, that is $\psi '(G)$ is a finite normal subgroup, hence contained in the center. Thus the
projection $\tilde{\rho}$ into $PSL(2,K)$ is trivial on $G$ so we get a factorizing representation
$\psi : \pi _1(Y) \rightarrow PSL(2,K)$.
\eop

\begin{corollary}
\label{factorequivalence}
Suppose $X$ is a smooth quasiprojective variety and $\rho : \pi _1(X,x)\rightarrow SL(2,K)$ a Zariski-dense representation
with $K$ a field. Then $\rho$ factors through a map to a DM-curve, if and only if it projectively factors through
a map to an orbicurve.
\end{corollary}
\eop

\begin{lemma}
\label{zdhyperbolicfactor}
If $\rho : \pi _1(X,x )\rightarrow SL(2,K)$ is a Zariski-dense representation which factors through
a map $f:X\rightarrow Y$ to a DM-curve, then $Y$ is a hyperbolic DM-curve. 
\end{lemma}
{\em Proof:}
This follows from Lemma \ref{zdhyperbolic} above. 
\eop

\subsection{Invariance under open sets and finite coverings}

\begin{lemma} 
\label{extendfactorization}
Suppose  $U\subset X$ is a Zariski open set in a smooth quasiprojective variety.
Suppose $K$ is an infinite field and $\rho : \pi _1(X,x)\rightarrow SL(2,K)$ is a Zariski-dense
representation such that the restriction $\rho |_U$
factors through a DM-curve. 
Then $\rho$ itself factors through a DM-curve.
\end{lemma}
{\em Proof:}
It is easier to use Corollary \ref{factorequivalence} to replace factorization through a DM-curve,
by projective factorization through an orbicurve in the statement. Thus we look at the representation
$\tilde{\rho }$ projected into $PSL(2,K)$. 

Let $f:U\rightarrow Y$ be the factorization map to an orbicurve $Y$.
We can restrict the open set $U$, so we can suppose that $Y$ is a smooth curve. 
Let $\overline{Y}$ denote the smooth projective compactification of $Y$. 
For $P\in \overline{Y}-Y$ consider the monodromy transformation for $\tilde{\rho} _Y$ around $P$.
If it is of infinite order, then leave out the point $P$.
If it is of some finite order $n$, then add an orbifold point of order $n$ at $P$.
This gives us an orbicurve partial completion $Y'$ of $Y$ with a representation $\tilde{\rho}_{Y'}$ restricting
to $\tilde{\rho}_Y$ on $Y$. 

Now we claim that the map $f$ extends to a map $f':X\rightarrow Y'$. 
For this, note that there is a birational modification $p:X''\rightarrow X$ containing $U\subset X'$
as the complement of a divisor with normal crossings,
and a map $f'' : X'' \rightarrow \overline{Y}$. On the other hand, for any component of the 
complementary divisor in $X''$, the monodromy around that component is trivial for the representation $\tilde{\rho}$.
This implies that any divisor component maps to a point of $Y'$ (i.e. maps to a point in $\overline{Y}$ where
the monodromy is not of infinite order). Thus we obtain a map $X'' \rightarrow (Y')^{\rm coarse}$.

By looking at the order of the monodromy around any complementary divisor component,
we get a lifting to a map $X'' \rightarrow Y'$. Suppose $P\in (Y')^{\rm coarse}$ is an orbifold point with index $n$.
The inverse image $(f'')^{-1}(P)=\sum _i r_iB_i$ is a divisor in $X''$. If $B_i$ is a component with multiplicity
$r_i$ then $r_i$ is divisible by $n$, because $\tilde{\rho}|Y$ \revision has monodromy of order $n$ around $P$ and
the pullback $f^{\ast}(\tilde{\rho})Y)$ extends to $X$ so it extends across the component $B_i$. Now each $r_i$ is
divisible by $n$ so the divisor $(f'')^{-1}(P)$ is divisible by $n$, which says that $X''\rightarrow  (Y')^{\rm coarse}$
lifts to a map into $Y'$ locally at $P$.

Finally, we can go to a model $X^{(k)}$ birational over $X''$
and with the property that $X^{(k)}$ \revision is obtained from $X=X^{(0)}$ by a sequence of blow-ups 
$\varphi _{i+1}:X^{(i+1)}\rightarrow X^{(i)}$ with smooth centers $V^{(i+1)}$.
Working by downward induction on $i$, we descend the map $X^{(k)}\rightarrow Y'$ to a map $X\rightarrow Y'$
as follows. Suppose we are given $f_{i+1}: X^{(i+1)}\rightarrow Y'$. For any point on the center of the blow-up
$P\in V^{(i+1)}\subset X^{(i)}$ the inverse image $\varphi _{i+1}^{-1}(P)$ is a connected variety such that
$$
f_{i+1}^{\ast}(\tilde{\rho}_{Y'})|_{\varphi _{i+1}^{-1}(P)} = \rho |_{\varphi _{i+1}^{-1}(P)}
$$
\revision
is trivial. Therefore $f_{i+1}$ maps $\varphi _{i+1}^{-1}(P)$ to a single point in $Y'$. Thus $f_{i+1}$ descends to a map
$$
X^{(i)}\rightarrow (Y')^{\rm coarse},
$$
and looking again at the monodromy around orbifold points as in the previous paragraph we get that it lifts to
$$
X^{(i)}\rightarrow Y'.
$$
Continuing by downward induction on $i$ we eventually get to a map $f':X\rightarrow Y'$ as required. 
This constitutes a factorization of the representation $\rho$, indeed
$(f')^{\ast} (\tilde{\rho}_{Y'})= \rho $ because this is true over the Zariski open dense set $U$. 
\eop

\begin{lemma} 
\label{descendfactorization}
Suppose  $p:Z\rightarrow X$ is a quasi-finite map of smooth quasiprojective varieties. 
Suppose $K$ is an infinite field and $\rho : \pi _1(X,x)\rightarrow SL(2,K)$ is a Zariski-dense
representation such that the pullback $\rho |_Z$
factors through a DM-curve. 
Then $\rho$ itself factors through a DM-curve.
\end{lemma}
{\em Proof:}
This is a sort of baby case of the Shafarevich factorization, see \cite{Kollar} as well as \cite{Campana94} 
\cite{Eyssidieux} \cite{Katzarkov}. We do it explicitly here. 

For constructing the factorization, we may pass to a Zariski open set of $X$ and then invoke
Lemma \ref{extendfactorization} above. 
In particular, we may assume that $p$ is a finite map. 

Let $f:Z \rightarrow Y$ be the factorization map. By \ref{zdhyperbolicfactor}, $Y$ is an hyperbolic DM-curve.
Therefore there exists a finite covering $Y'\rightarrow Y$ by a smooth curve. Replacing $Z$ by
its fiber product with $Y'$, we find that we may in fact assume that $\rho |_Z$ factors through
a smooth curve. Thus, take now this situation and assume that $f:Z\rightarrow Y$ is the factorization
map to a smooth curve. We may also take the Stein factorization so we can assume that the general fiber
of $f$ is connected. Remove from $X$ the closure of the image of all of the non-smooth or reducible fibers. 
Thus we may assume that all of the fibers of $f$ are irreducible and smooth.  

We can assume that $Z$ is Galois over $X$ with group $G$. 
We claim that $G$ permutes the fibers of $f$. This is because a fiber is characterized as a maximal 
connected smooth curve such that the restriction of $\rho$ is trivial. Therefore the action of $G$
extends to an action on $Y$ compatible with the map $f$. Let $V:= Y/G$ be the DM-curve quotient.
We obtain a factorization map $X\rightarrow V$ for $\rho$. 

The fact that the action of $G$ on $Z$ preserves $\rho$ (because $\rho$ is a pullback from $X$)
allows us to define an action on $\rho _Y$ covering the action on $Y$, so $\rho _Y$ descends
to a representation $\rho _V$ on $V$ whose pullback to $X$ is $\rho$. 
\eop

\subsection{Invariance under field extension}

Suppose $f:X\rightarrow Y$ is a map. Let $K$ denote the kernel of the induced map
$\pi _1(X,x) \rightarrow \pi _1(Y,f(x))$. A representation $\rho : \pi _1(X,x)\rightarrow G$ 
factors through $f$ if and only if the restriction of $\rho$ to $K$ is trivial. This immediately
gives the following statement.

\begin{lemma}
Suppose $\alpha : G\hookrightarrow H$ is an injection of groups. A representation 
$\rho : \pi _1(X,x)\rightarrow G$ factors through a map $f:X\rightarrow Y$, if and only if the 
composed representation $\alpha \circ \rho : \pi _1(X,x)\rightarrow H$ factors through $f$.
\end{lemma}
\eop

A corollary is independence of factorization with respect to change of coefficient field.

\begin{corollary}
\label{invariancefieldextension}
Suppose $K\subset L$ is a field extension. Then a representation 
$\rho : \pi _1(X,x) \rightarrow SL(r,K)$ factors (resp. projectively factors) through a map
$f:X\rightarrow Y$ if and only if the representation obtained by extension of coefficients
$\rho _L : \pi _1(X,x)\rightarrow SL(r,L)$ factors (resp. projectively factors) through $f$.
\end{corollary}
\eop

\section{Construction of a pluriharmonic mapping to the Bruhat-Tits tree}
\label{sec-harmonic-initial}

An element $A\in PSL(2,K)$ is said to be {\em unipotent} if its representatives in $SL(2,K)$ are
conjugate to upper triangular matrices with $\pm 1$ on the diagonal. An element is {\em quasi-unipotent}
if there is some $n$ such that $A^n$ is unipotent. 
We consider the following situation, stated as a hypothesis for future reference. 

\begin{hypothesis}
\label{unboundedZD}
Fix a quasiprojective variety $X$ with compactification $\overline{X}$ whose complementary divisor 
$D$ is a union of smooth irreducible components $D_i$ intersecting at normal crossings. Let
$x$ denote a basepoint in $X$.

Let $\Oo _K$ be a complete local ring with fraction field $K$ and residue field $F$;
we assume that $F$ is finite. Suppose we have a representation $\rho : \pi _1(X,x)\rightarrow PSL(2,K)$ 
such that the monodromy transformations around the $D_i$ are quasi-unipotent. We also assume that the
monodromy transformations around the $D_i$ are not the identity (i.e. that $X$ is a maximal open
subset of $\overline{X}$ of
definition for the representation $\rho$). We assume that the Zariski closure of the image of
$\rho$ is all of $PSL(2,K)$, and we assume that the image is not contained in any compact subgroup. 
\end{hypothesis}

In the next section, we will prove Theorem \ref{mainfactorization} saying that 
under the above hypothesis
the representation $\rho$ factors through a DM-curve. In general Hypothesis \ref{unboundedZD} 
will be in effect throughout this
section and the next. 

The basic technique is to let $PSL(2,K)$ act on the {\em Bruhat-Tits tree} $\Tt$, and to choose an equivariant pluriharmonic
map $\widetilde{X}\rightarrow \Tt$ via the theory of Gromov and Schoen \cite{GromovSchoen}. 
The first step, not totally trivial in the quasiprojective
case, is to choose an initial map having finite energy on the fundamental region ${\rm Reg}(\widetilde{X}/X)$.
By \cite{GromovSchoen} this deforms to a pluriharmonic map which will then be studied in \S \ref{sec-harmonic-properties}. 

This theory has already been treated by Jost and Zuo in \cite{JostZuo2}. They treat the
case of harmonic maps to buildings of arbitrary rank. In the case of trees which was treated specifically
in Gromov-Schoen \cite{GromovSchoen}, and which is by now pretty well-known, 
many parts of the argument are easier to understand so it seems worthwhile to recall them here. 
Also, we feel that the discussion of the choice of initial map on \cite[p. 12]{JostZuo2} needs some amplification.

Looking at harmonic mappings to trees or more generally buildings, mirrors much of the study of harmonic bundles which
has been extended to the quasiprojective case \cite{Biquard} \cite{actesToulouse} \cite{JostZuo1} 
\cite{Li} \cite{LiNarasimhan} \cite{LiNarasimhan2} \cite{TMochizuki} \cite{TMochizuki2} \cite{NiRen} \cite{hbnc}. 
Other related techniques are also available \cite{Delzant2} \cite{Gromov} \cite{NapierRamachandran}
 \cite{NapierRamachandran3} \cite{NapierRamachandran2}.

\subsection{The Bruhat-Tits tree}

Recall the construction of the {\em Bruhat-Tits tree} associated to $K$, extensively studied by Bass and Serre 
\cite{Bass} \cite{Serre} \cite{Serre2}. 
There is a 
uniformizing parameter $t\in \Oo _K$ which generates the maximal ideal.
The residue field $F:= \Oo _K /(t)$ is assumed to be finite, with $q=p^k$ elements.
A {\em lattice} is a finitely generated $\Oo _K$-submodule $M\subset K^2$ 
which spans $K^2$ over $K$. It follows that $M\cong \Oo _K^2$.

Say that two lattices $M_1$ and $M_2$ are {\em equivalent} if one is a scalar multiple of the other
by an element of $K^{\ast}$. The vertices of the {\em Bruhat-Tits tree} are defined to be the equivalence
classes of lattices. Two lattices are {\em adjacent} if after choosing appropriate representatives
up to equivalence, one is contained in the
other and if the quotient is an $\Oo _K$-module of length $1$ isomorphic to $F$. The  edges of the Bruhat-Tits tree
are defined as connecting lattices
which are adjacent in the above sense. See \cite{Serre2}. There are $q+1$ edges at each vertex.

Metrize the tree by assigning a linear metric of length $1$ to each edge. 

The group $SL(2,K)$ acts isometrically on $\Tt$ via its action on $K^2$, thus $g\in SL(2,K)$ sends
a lattice $M$ to $gM$. 
The center $\{ \pm 1\}$ acts trivially on $\Tt$, so the action factors through the group
$PSL(2,K)$. 

Choosing a basis for a lattice $M$ gives an isomorphism between 
$PSL(2,\Oo _K)$ and the 
subgroup of elements which fix $M$, in particular the latter group is compact.  

Conversely any compact subgroup  $H \subset PSL(2,K)$ fixes a lattice. Indeed, the action on the tree is continuous,
so the orbit of a point under the action of $H$ is compact. The convex subset spanned by this orbit is $H$-invariant
and compact, so the set of endpoints is finite and the barycenter of this finite set is a fixed point for $H$.

Recall that an {\em end} of $\Tt$ is an equivalence class of injective isometric paths $\rr _{\geq 0}
\rightarrow \Tt$, under the equivalence relation that two paths are said to be equivalent if
they remain a bounded distance apart, or equivalently if their image subsets coincide after a finite time.
Recall also that a parabolic subgroup of $PSL(2,K)$ is a subgroup conjugate to the group of  upper-triangular
matrices. Such a subgroup is uniquely determined by its corresponding rank one subspace $L\subset K^2$.

\begin{lemma}
\label{ends}
The subgroup of elements fixing an end is a parabolic subgroup, and a parabolic subgroup fixes a unique end. In this way
the set of ends of $\Tt$ is naturally in one-to-one correspondence with the set of parabolic subgroups
of $PSL(2,K)$. 
\end{lemma}
{\em Proof:}
See  \cite{Bass} \cite{Serre} \cite{Serre2}. 
\eop

\subsection{Quasi-unipotent matrices}

Let $\overline{K}$ denote the algebraic closure of $K$. A matrix $A\in PSL(2,K)$ is quasi-unipotent if and only if
the image in $PSL(2,\overline{K})$ is conjugate to an upper triangular matrix whose diagonal entries are roots of unity.
The diagonal entries are $\alpha , \alpha ^{-1}$ with $\alpha ^n = 1$ for some $n$. 

Let $Z(A)\subset PSL(2,K)$ be the centralizer of $A$, that is the group of matrices in $PSL(2,K)$ which commute with $A$. 
It is an algebraic group and its extension of scalars to $\overline{K}$ is the same as the centralizer of $A$ in $PSL(2,\overline{K})$.
The centralizer will play an important role in our construction below. We divide into three cases according to the 
type of $Z(A)$ up to conjugacy. 

\begin{lemma}
\label{quasiunip-classif}
Suppose $A\in PSL(2,K)$ is a nontrivial quasiunipotent matrix. Then one of three possibilities holds:
\newline
---{\em (unipotent)} the matrix $A$ is conjugate to an upper triangular matrix with $\pm 1$ on the diagonal,
and the centralizer $Z(A)$ is conjugate to the unipotent group of all such matrices;
\newline
---{\em (split torus)} the matrix $A$ is conjugate to a diagonal matrix with entries $\alpha ^{\pm 1}$ where $\alpha \neq \pm 1$
is a nontrivial root of unity, and the centralizer $Z(A)$ is conjugate to the split torus $\Gm$ of diagonal matrices; or
\newline
---{\em (non-split torus)} the matrix $A$ has eigenvalues $\alpha ^{\pm 1}$ which are not in $K$, and the centralizer $Z(A)$ is
a non-split torus which becomes conjugate in $PSL(2,\overline{K})$ to the group of diagonal matrices.
\end{lemma}
{\em Proof:}
Over $\overline{K}$ the matrix can be put in Jordan normal form. If it has distinct eigenvalues, then they are either in $K$
which gives the split torus case, or not in $K$ giving the non-split torus case. If both eigenvalues are the same, that is to
say $\pm 1$, and if the matrix is nontrivial in $PSL(2,\overline{K})$ then it is conjugate to 
$\left( \begin{array}{cc}1 & 1 \\ 0 & 1 \end{array}\right)$. In this case, there is an unique eigenvector, and the dimension
of ${\rm ker}(A-1)$ is invariant under field extension. Thus 
$A$ has an eigenvector over $K$ so it can be conjugated to a unipotent upper-triangular matrix in $PSL(2,K)$. This gives the unipotent case.
\eop

\begin{lemma}
\label{samecentralizers}
Suppose $A$ and $B$ are commuting quasi-unipotent matrices, both nontrivial in $PSL(2,K)$. Then $Z(A)=Z(B)$. 
\end{lemma}
{\em Proof:}
It suffices to do this in $PSL(2,\overline{K})$. Thus, we may assume that $A$ is either diagonal with distinct eigenvalues, 
or $\pm 1$ times an upper
triangular matrix. Respectively we see that $B$ has to be diagonal or $\pm 1$ times an upper
triangular matrix. Again, respectively in these cases the centralizers in $PSL(2,\overline{K})$ are either the group of diagonal matrices,
or the group of projectively upper triangular matrices. The centralizers of $A$ and $B$ in $PSL(2,\overline{K})$ are the same in either case.
Formation of the centralizer commutes with field extension, so this shows that $Z(A)=Z(B)$ in $PSL(2,K)$.
\eop

\subsection{Local study at infinity}

Recall that part of our assumption is that the monodromy around all components of the divisor $D$ is nontrivial in
$PSL(2,K)$ Thus, the connected components of $D$ form well-defined subsets depending on $\rho$. 
The choice of initial map is a local problem near connected components of $D$. In view of this,
we will for simplicity assume that $D$ is connected for most of the remainder of this section.

Let $D_j$ denote the irreducible components of $D$, which we are assuming are smooth. Let $N_j$ denote a tubular neighborhood of
$D_j$ in $\overline{X}$ and let $N_j^{\ast}$ denote the complement of $D$, that is $N_j^{\ast}=N_j\cap X$. 
Let $D_{ij}:= D_i\cap D_j$ etc.\ and let $N_{ij}:= N_i\cap N_j$ which we assume is a tubular neighborhood of $D_{ij}$ (with polydisc transversal section).
Let $N:= \bigcup _jN_j$ be the tubular neighborhood of $D$ and $N^{\ast} := N\cap X$. 

For now, we assume having chosen a basepoint $x\in N^{\ast}$. Let $\widetilde{N}^{\ast}$ denote the resulting universal covering
of $N^{\ast}$. 
Let ${\rm Reg}(\widetilde{N}^{\ast} / N^{\ast} ) \subset \widetilde{N}^{\ast}$ denote a closed fundamental region for  $N^{\ast}$,
as will be chosen more specifically further along (Lemma \ref{regionproperties}).
Let
$\pi _1(N^{\ast})$ denote the fundamental group of $N^{\ast}$ based at $x$, also viewed as the group of deck transformations of $\widetilde{N}^{\ast}$.

Let $\gamma _j$ be a path going around the divisor component $D_j$. Recall that we are assuming that $\rho (\gamma _j)$
is quasiunipotent and nontrivial in $PSL(2,K)$. We assume that the basepoint $x$ is near $D$ and that the paths
used to define the $\gamma _j$ stay within the tubular neighborhood $N$ of $D$.

\begin{lemma}
\label{samecentralizer}
Assuming Hypothesis \ref{unboundedZD} as well as connectedness of $D$ and with the above notations,
let $Z(\rho (\gamma _j))$ be the centralizer in $PSL(2,K)$ of the matrix $\rho (\gamma _j)$. 
Then the $Z(\rho (\gamma _j))$ are all the same group for different $j$. 
\end{lemma}
{\em Proof:}
We may assume that the basepoint $x$ is near one of the components $D_0$. Write 
$$
\gamma _j = \eta \xi _j \eta ^{-1}
$$
where $\xi _j$ is a path going around $D_j$ from a nearby basepoint $x_j$, and $\eta $ is a path going from $x$ to $x_j$. 
Decompose $\eta$ into paths
$$
\eta = \eta _{1} \cdots \eta _{k}
$$
where for any $a=1,\ldots , k$, the piece $\eta _{a}$ is a path staying in $N_{i_a}$. Let $\zeta _a$ be a loop going from the endpoint of $\eta _a$
around the divisor component $D_{i_a}$ and let $\zeta '_a$ be a loop going around the divisor component $D_{i_{a+1}}$. Note that the
endpoint of $\eta _a$ is in $N_{i_ai_{a+1}}$ so we can look at both of these loops. 

Commutativity of the fundamental group of the product of punctured polydiscs says that $\zeta _a$ and $\zeta '_a$ commute.
This step corresponds to the phrase on lines 4-7 of page 12 of \cite{JostZuo2} and in a certain sense the rest of our discussion above and below consists
of filling in the surrounding details. 

For any $i$ the loop around the divisor $D_i$ is in the center of $\pi _1(N_i^{\ast})$, as can be seen by restricting to a Zariski open
set of $D_i$ over which the normal bundle is trivialized, and noting that the fundamental group of the open set surjects to the original one. 
This implies that
$$
\zeta '_a = \eta _{a+1} \zeta _{a+1} \eta _{a+1}^{-1}.
$$
Set 
$$
\tau _a := \eta _1 \cdots \eta _a \zeta _a \eta _a ^{-1} \cdots \eta _1^{-1}.
$$
Combining with commutativity of $\zeta _a$ and $\zeta ' _a$ from the previous paragraph we conclude that $\tau _a$ and $\tau _{a+1}$ commute. 
By our hypothesis on $D$, all of the $\rho (\tau _a)$ are nontrivial quasiunipotent matrices in $PSL(2,K)$. By 
Lemma \ref{samecentralizers} we find
$$
Z(\rho (\tau _a)) = Z(\rho (\tau _{a+1})).
$$
It follows by induction on $a$ that the centralizers $Z(\rho (\tau _a))$ are the same for all $a$. 
We can extend some of our notation to the boundary cases $a=0$ and $a=k$ as is left to the reader. With this, 
$\tau _0 = \gamma _0$ whereas $\tau _k = \gamma _j$, and we get
$$
Z(\rho (\gamma _0)) = Z(\rho (\gamma _j)).
$$
This proves the lemma.
\eop

In view of this lemma, let $Z\subset PSL(2,K)$ denote the centralizer. We can assume by an appropriate change of
basis that it is either the group of unipotent
upper triangular matrices,  the group of diagonal matrices, or a non-split torus, see Lemma \ref{quasiunip-classif}.
Let $NZ\subset PSL(2,K)$ denote its normalizer. In the unipotent case, $NZ$ is the parabolic subgroup $P$ of upper-triangular matrices;
in the split torus case, $NZ/Z$ is the Weyl group of order two transposing the two standard basis elements of $K^2$,
and in the non-split torus case, $NZ/Z$ has order $1$ or $2$ by mapping it into the Weyl group over the algebraic closure.

\begin{lemma}
\label{toNZ}. 
Let $NZ$ denote the normalizer of the subgroup $Z$ which is the centralizer of any $\rho (\gamma _j)$. 
Then 
the image by $\rho$ of the fundamental group  $\pi _1(N^{\ast})$ lies inside $NZ$.
Furthermore, the representation gives a map of exact sequences 
$$
\begin{array}{ccccccc}
1 \rightarrow & \pi _1(N^{\ast} \times _{D} \widetilde{D})  & \rightarrow & \pi _1(N^{\ast})& \rightarrow & \pi _1(D ) &\rightarrow 1 \\
 & \downarrow & & \downarrow & & \downarrow & \\
1 \rightarrow &Z  & \rightarrow & NZ & \rightarrow & NZ/Z &\rightarrow 1 .
\end{array}
$$
\end{lemma}
{\em  Proof:}
Let $D_0$ be the irreducible component of $D$ near the basepoint $x$, and let $\gamma _0$ be the loop going from $x$ around $D_0$.
We claim that for any $\eta \in \pi _1(N^{\ast}, x)$ then
$$
Z(\rho (\eta \gamma _0 \eta ^{-1})) = Z(\rho (\gamma _0)).
$$
The proof of this claim is exactly the same as the proof of Lemma \ref{samecentralizer} above: decompose 
$$
\eta = \eta _{1} \cdots \eta _{k},
$$
fix the loops $\zeta _a$ as previously and set 
$$
\tau _a := \eta _1 \cdots \eta _a \zeta _a \eta _a ^{-1} \cdots \eta _1^{-1}.
$$
This time we choose $\zeta _k := \gamma _0$ which gives for endpoints $\tau _0 = \gamma _0$ and $\tau _k= \eta \gamma _0 \eta ^{-1}$.
As before $\tau _a$ and $\tau _{a+1}$ commute, so $Z(\rho (\tau _a)) = Z(\rho (\tau _{a+1}))$ by Lemma \ref{samecentralizers}. 
Thus $Z(\rho (\tau _k)) = Z(\rho (\tau _0))$ which gives the claim.

On the other hand, the centralizers of conjugate matrices are conjugate, that is
$$
Z(\rho (\eta \gamma _0 \eta ^{-1})) = \rho (\eta  ) Z (\rho (\gamma _0))\rho ( \eta ) ^{-1}.
$$
Putting this together with the claim and the notation $Z = Z(\rho (\gamma _j))$ for any $j$, we find that 
$$
\rho (\eta  )Z\rho ( \eta )^{-1} =  Z.
$$
Thus $\rho (\eta ) \in NZ$ is in the normalizer of $Z$. 

We finish by showing the map of exact sequences. The fact that the image of $\rho$ on $\pi _1(N^{\ast})$ is contained in
the normalizer of $Z$ implies that the normal subgroup of $\pi _1(N^{\ast})$ generated by the $\gamma _j$
maps into $Z$. However, this normal subgroup is the kernel of the map 
$$
\pi _1(N^{\ast}) \rightarrow \pi _1(N) = \pi _1(D),
$$
as may be seen by noting that $N$ is obtained homotopically from $N^{\ast}$ by adding some $2$-cells whose boundaries are the loops $\gamma _j$,
and then adding only cells of dimensions $\geq 3$. Thus the kernel of the upper exact sequence maps into $Z$,
so $\rho$ induces the map of exact sequences as claimed. 
\eop

\begin{corollary}
\label{allcommute}
The matrices $\rho (\gamma _j)$ all commute with each other.
\end{corollary}
{\em Proof:}
They are all inside the group $Z$, which in the case of $PSL(2,K)$ is commutative.
\eop

{\em Remark}---
\label{sigmaremark}
Let $\Sigma$ denote the incidence simplicial complex of the divisor $D$. This has appeared in the works of Stepanov \cite{Stepanov} and Thuillier \cite{Thuillier}.
There is a vertex for each divisor component $D_j$, an edge for each irreducible component of $D_{ij}$, and higher simplices for the multiple intersections.
We have a map $N\rightarrow \Sigma$
inducing $\pi _1(N^{\ast}) \rightarrow \pi _1(\Sigma )$. In Lemma \ref{toNZ}, we could have replaced $D$ by $\Sigma$, that is to say
the representation sends the kernel of $\pi _1(N^{\ast})\rightarrow \pi _1(\Sigma )$ into the center $Z$. We don't need this for our
proof, and the end theorem of the paper shows that the situation we are considering here doesn't really arrive
in practice; however this remark might be interesting from a geometrical point of view in related situations such as \cite{Delzant2}.   

\subsection{Local construction of an initial map---generalities}

Fix a complete finite-volume K\"ahler metric of Poincar\'e type on $X$, see \cite{SiuYau}, \cite{JostYau}, 
\cite{CattaniKaplanSchmid}, \cite{JostZuo2}.
Our first problem will be to choose an initial finite-energy equivariant map to the Bruhat-Tits tree under the action of $PSL(2,K)$.

Composing our representation $\rho$
with the inclusion $\pi _1(N^{\ast})\rightarrow \pi _1(X)$ we obtain an action denoted again by $\rho$ of $\pi _1(N^{\ast})$ on the
tree $\Tt$. 
We would like to construct an initial map $\phi : \widetilde{N}^{\ast}\rightarrow \Tt$, which is $\rho$-equivariant and whose restriction to the
fundamental region ${\rm Reg}(\widetilde{N}^{\ast} / N^{\ast} )$ has finite energy.

We assume that we have chosen projections $p_j:N_j\rightarrow D_j$ such that on $N_{ij}$ the different projections commute and give
a projection $p_{ij}$ to $D_{ij}$, and similarly for higher intersections. 
In order to insure the finite-energy condition, we shall try to construct a map which is constant in the direction of the fibers of $N_j\rightarrow D_j$.
This would mean in particular that it is constant on the fibers of $N_{ij}\rightarrow D_{ij}$ and similarly for multiple intersections. 
The discussion by Jost and Zuo in \cite[pp 8-12]{JostZuo2} serves to explain why this condition yields a finite-energy map. To be more precise
about this, say that two points in $\widetilde{N}^{\ast}$ (or equally well, $N^{\ast}$ or $N$) 
are {\em projection-equivalent} if they are joined by a path which projects into 
$N^{\ast}$ to a path which stays in any fiber of any projection $p_j$ which it meets. We refer to \cite{JostZuo2} for the fact that any 
Lipschitz and piecewise smooth map $\phi : \widetilde{N}^{\ast}\rightarrow \Tt$ which sends projection-equivalent points to the same point in $\Tt$,
has finite energy on the fundamental region.

If we take the quotient of $N$ by the equivalence relation of projection-equivalence we obtain a space which is homotopy-equivalent to $D$. Thus,
by modifying the projections $p_j$ near the multiple intersections, we can define a map $\sigma : N \rightarrow D$ with the property that
all fibers of the projections $p_j$ go to single points in $D_j$. The restriction $\sigma |_D$ will  no longer be the identity,
for example points in $N_{ij}\cap D$ have to get mapped to $D_{ij}$. However, $\sigma |_D$ will be homotopic to the identity. 
The map $\sigma$ induces
$$
\sigma _{\ast} : \pi _1(N^{\ast})\rightarrow \pi _1(D ),\;\;
\widetilde{\sigma} : \widetilde{N}^{\ast} \rightarrow \widetilde{D},
$$
where $\widetilde{D}\rightarrow D$ is the 
universal cover corresponding to  $\sigma (x)$ as basepoint for $D$. 

We assume that the choice of fundamental region 
${\rm Reg}(\widetilde{N}^{\ast} / N^{\ast} )$ is
made compatibly with a choice of ${\rm Reg}(\widetilde{D}/ D )$ via the map $\sigma$. 
More precisely, proceed as follows. Let $N^+$ denote the real blow-up of the divisor $D$ in $N$.
It consists of adding a boundary to $N^{\ast}$, the boundary being homeomorphic to the boundary of a tubular
neighborhood of $D$. The inclusion of an open dense subset $N^{\ast}\hookrightarrow N^+$ is a homotopy equivalence. 
We assume that our map extends to a map 
$$
\sigma ^+ : N^+\rightarrow D
$$
and the homotopy equivalence between $N^{\ast}$ and $N^+$ commutes with $\sigma ^+$. 
Let $\widetilde{N}^+$ denote the universal cover of $N^+$ corresponding to the same basepoint $x$,
so again we have an open dense inclusion $\widetilde{N}^{\ast} \hookrightarrow \widetilde{N}^+$ on which the homotopy equivalence lifts. 

Choose a fundamental region
$$
{\rm Reg}(\widetilde{D}/ D ) \subset \widetilde{D}.
$$
Now $\sigma ^+$ induces $\widetilde{\sigma} ^+: \widetilde{N}^{+} \rightarrow \widetilde{D}$.
Consider the region 
$$
Q^+:= (\widetilde{\sigma}^+) ^{-1}({\rm Reg}(\widetilde{D}/ D )) \subset \widetilde{N}^{+},
$$
and let $Q^{\ast} = Q^+\cap \widetilde{N}^{\ast}$.

Note that $\pi _1(N^+)= \pi _1(N^{\ast})$ and $\sigma _{\ast} = (\sigma ^+)_{\ast}$ maps this group to $\pi _1(D)$.
The group $\ker (\sigma _{\ast})\subset \pi _1(N^{\ast})$ acts freely on $Q^+$, and the quotient is
$$
Q^+/\ker (\sigma _{\ast}) = N^{+} \times _D {\rm Reg}(\widetilde{D}/ D ).
$$
Let ${\rm Reg}(\widetilde{N}^{+} / N^{+} )\subset Q^+$ be a fundamental region for the action of $\ker (\sigma _{\ast})$.

In general $\ker (\sigma _{\ast})$ will not be finitely generated. However, $N^{+} \times _D {\rm Reg}(\widetilde{D}/ D )$ is
homeomorphic to a finite simplicial complex, so its fundamental group is finitely generated. Thus, $Q^+$ will usually be disconnected,
a disjoint union of coverings of $N^{+} \times _D {\rm Reg}(\widetilde{D}/ D )$ corresponding to finitedly generated groups. 

The fact that the quotient is a finite complex means we can choose a compact fundamental region in $Q^+$, indeed in one of the
connected components. Thus, ${\rm Reg}(\widetilde{N}^{+} / N^{+} )$ is compact and is a finite simplicial complex. 

Now let ${\rm Reg}(\widetilde{N}^{\ast} / N^{\ast} ) := {\rm Reg}(\widetilde{N}^{+} / N^{+} ) \cap N^{\ast}$. 
Note that it is a fundamental region in $Q^{\ast}$ for the action of $\ker (\sigma _{\ast})$, surjecting to the quotient
$$
Q^{\ast}/\ker (\sigma _{\ast}) = N^{\ast} \times _D {\rm Reg}(\widetilde{D}/ D ).
$$

\begin{lemma}
\label{regionproperties}
The region ${\rm Reg}(\widetilde{N}^{\ast} / N^{\ast} )$ constructed in this way is fundamental, that is any point of
$\widetilde{N}^{\ast}$ is a $\pi _1(N^{\ast})$-translate of a point in the region; also the subset 
$\Aa \subset \pi _1(N^{\ast})$ of elements $\alpha$ such that 
$$
{\rm Reg}(\widetilde{N}^{\ast} / N^{\ast} ) \cap \alpha \cdot {\rm Reg}(\widetilde{N}^{\ast} / N^{\ast} ) \neq \emptyset
$$
is finite. 
\end{lemma}
{\em Proof:}
It suffices to do this for ${\rm Reg}(\widetilde{N}^{+} / N^{+} )$ because
$$
{\rm Reg}(\widetilde{N}^{\ast} / N^{\ast} ) = {\rm Reg}(\widetilde{N}^{+} / N^{+} ) \times _{N^+}N^{\ast},
$$
and $\widetilde{N}^{\ast} = \widetilde{N}^{\ast}\times _{N^+}N^{\ast}$ too. 

The map $\sigma _{\ast} : \pi _1(N^{\ast})\rightarrow \pi _1(D)$ is surjective. This implies that any point in
$\widetilde{N}^{+}$ is a translate of a point in $Q^+$. On the other hand, any point in $Q^+$ is a translate of a point in 
${\rm Reg}(\widetilde{N}^{+} / N^{+} )$ by choice of the latter. 

The set of points $\Aa$ is a subset of the
set $\Aa ^+$ defined in the same way for ${\rm Reg}(\widetilde{N}^{+} / N^{+} )$ so for the second statement it suffices to
prove that $\Aa ^+$ is finite. 
For this, note that $\widetilde{N}^{\ast}$ is the universal covering of a finite complex, and 
${\rm Reg}(\widetilde{N}^{+} / N^{+} )$ is a subset bounded with respect to any $\pi _1(N^+)$-invariant metric. 
There are only finitely many group elements corresponding to paths of bounded length, so $\Aa ^+$ is finite. 
\eop

By construction we now have a map between fundamental regions denoted
$$
\sigma _{\rm Reg} : {\rm Reg}(\widetilde{N}^{\ast} / N^{\ast} )\rightarrow {\rm Reg}(\widetilde{D}/ D ).
$$

\begin{lemma}
\label{inorderto}
In order to define an equivariant map $\phi$ on $\widetilde{N}^{\ast}$ with finite energy on the fundamental region, 
it suffices to construct a Lipschitz piecewise $C^1$ map 
$$
\phi _{\rm Reg} : {\rm Reg}(\widetilde{N}^{\ast} / N^{\ast} )\rightarrow \Tt
$$
sending projection-equivalent points to the same point in $\Tt$, and
with the property that for any $x,y\in {\rm Reg}(\widetilde{N}^{\ast} / N^{\ast} )$ and any $\alpha \in \pi _1(N^{\ast})$ such that
$\alpha x = y$, we should have
$$
\phi _{\rm Reg}(y) = \rho (\alpha )\phi _{\rm Reg}(x).
$$
\end{lemma}
{\em Proof:}
Suppose given $\phi _{\rm Reg}$. Then define $\phi : \widetilde{N}^{\ast} \rightarrow \Tt$ as follows. For $u\in \widetilde{N}^{\ast}$
choose a group element $g\in \pi _1(N^{\ast})$ such that $gu\in {\rm Reg}(\widetilde{N}^{\ast} / N^{\ast} )$. Put
$$
\phi (u) := \rho (g) ^{-1} \phi _{\rm Reg} (gu).
$$
If $g'$ is another element with $g'u\in {\rm Reg}(\widetilde{N}^{\ast} / N^{\ast} )$
then putting $\alpha = g' g^{-1}$, $x=gu$ and $y=g'u$ we have $\alpha x = y$ so this is data such as in the hypothesis of the lemma. In this case,
$$
\rho (g) ^{-1} \phi _{\rm Reg} (gu) = \rho (g') ^{-1} \rho (\alpha ) ^{-1} \phi _{\rm Reg} (x) =
\rho (g') ^{-1}\phi _{\rm Reg}(\alpha x) = \rho (g') ^{-1}\phi _{\rm Reg}(g'u),
$$
which shows that $\phi$ is well-defined. If $\phi _{\rm Reg}$ is Lipschitz piecwise $C^1$ then the same is true of $\phi$. 
As stated above, by \cite{JostZuo2} if $\phi _{\rm Reg}$ sends projection-equivalent
points to the same point in $\Tt$ then it will have finite energy. 
\eop

Our strategy will be to use this lemma and try to construct a map which in a certain sense factors through the projection $\sigma$. In the most delicate case
of unipotent monodromy, this isn't exactly possible, however it will be possible to construct a map which factors over the fundamental region.

According to Lemma \ref{samecentralizers}, the centralizer $Z$ is an invariant. Since the centralizers are different in the three cases of Lemma 
\ref{quasiunip-classif}, it means that the matrices $\rho (\gamma _j)$ all fall into the same case of that classification.
Thus, for our construction of an initial map we can treat separately these three cases, as we shall do in the next subsections in
increasing order of difficulty: the non-split torus, the split torus, and the unipotent case.

\subsection{The non-split torus case}

In both torus cases, the $\rho (\gamma _j)$ are commuting matrices of finite order. In particular, there exists a matrix $A$ which
also has finite order, and such that for each $j$ there is $m_j$ with $\rho (\gamma _j) = A^{m_j}$. In fact, $A$ can be taken as an appropriate product
of powers of the $\rho (\gamma _j)$. 

Let $R\subset \Tt$ be the set of fixed points of $A$. The action of $A$ fixes the ends of $R$. This implies that $A$ is contained in the parabolic
subgroup associated to any end. However, since the eigenvalue of $A$ is not defined over $K$, this is impossible: the parabolic subgroups
are conjugate to groups of upper-triangular matrices in particular the eigenvalues of elements of a parabolic subgroup are defined over $K$.
Thus, $R$ has no ends, in other words it is compact. 

Now, the group $Z(A)$ fixes $R$: if $\beta \in Z(A)$ and $x\in R$ then
$$
A(\beta x) = \beta (Ax) = \beta x, \;\;\; \Rightarrow \beta x \in R.
$$
It follows that $Z(A)$ fixes the barycenter of $R$. Let $R'\subset \Tt$ be the set of fixed points of 
$Z$. By what we have just said, it is nonempty. Note also that $A\in Z$ so $R'\subset R$, in particular $R'$
is compact too. 

If $\beta \in NZ$ and $x\in R'$ then for any $\xi \in Z$ we have
$$
\xi (\beta x) = \beta (\beta ^{-1}\xi \beta x) = \beta x \;\;\; \mbox{since} \;\; \beta ^{-1}\xi \beta \in Z,
$$
thus $\beta x\in R'$. In other words, $NZ$ sends $R'$ into itself. Again, this implies that the barycenter of $R'$  
is a fixed point of $NZ$. 

By Lemma \ref{toNZ}, the representation $\rho$ sends $\pi _1(N^{\ast})$ into $NZ$, thus $\Tt$ admits a fixed point for the 
action of $\pi _1(N^{\ast})$. We can take as initial finite-energy map the constant map sending all of $N^{\ast}$ to this
fixed point.

\subsection{The split torus case}

As above, there is a matrix $A\in Z$ of finite order such that the $\rho (\gamma _j)$ are powers of $A$.
In the split case, $A$ may be conjugated to a diagonal matrix with roots of unity on the diagonal. The
subset $R\subset \Tt$ of points fixed by $A$ is a real line connecting two ends of $\Tt$. The two ends correspond
to the parabolic subgroups of upper and lower triangular matrices. By direct calculation, the centralizer
$Z$ which is the diagonal torus $\Gm \subset PSL(2,K)$ preserves $R$. The normalizer $NZ$ is the extension of
$Z$ by the Weyl group of transpositions, and again we see directly that $NZ$ preserves $R$, with the 
nontrivial irreducible component acting by interchanging the two ends, reversing the orientation of $R$. 
By Lemma \ref{toNZ}, the action of $\pi _1(N^{\ast})$ via $\rho$ goes through $NZ$, hence it preserves $R\subset \Tt$. 

The morphism $NZ \rightarrow Aut(R)$ sends $A$ to the identity. Hence it sends the normal subgroup generated
by the $\rho (\gamma _j)$ to the identity. This normal subgroup is the kernel of the map $\pi _1(N^{\ast})\rightarrow \pi _1(N)$.
This map is the same as $\sigma _{\ast} : \pi _1(N^{\ast})\rightarrow \pi _1(D)$. 
Therefore the action of $\pi _1(N^{\ast})$ on $R\subset \Tt$ factors through an action of $\pi _1(D)$.

Now $D$ is a usual compact variety and we can choose an equivariant map $\widetilde{D}\rightarrow R$.
Composing with the projection $\sigma$ this gives an equivariant map $\widetilde{N}^{\ast} \rightarrow R$ which
sends projection-equivalent points to the same point in $R$, in particular it has finite energy on the
fundamental region.

\subsection{The unipotent case}

We would like to mimic the preceding construction. However, $Z$ is now the unipotent
subgroup of upper triangular matrices with $\pm 1$ on the diagonal. 
The fixed point sets of finitely generated subgroups of $Z$ are nonempty: they are contractible rays. The action
of $NZ$ will not preserve these rays. Basically, a diagonal matrix with nontrivial valuations of its diagonal elements will
change the fixed ray. If we think of $Z$ as isomorphic to the additive group of the field $K$, then any finitely generated
subgroup goes into a subgroup of the form $t^q\Oo _K$ where $t$ is the uniformizing parameter. The ray of fixed points depends on $q$ and
conjugating with a diagonal matrix can change $q$. 

We get around this problem by noting that we really only need to choose a map on the fundamental region, equivariant with respect
to the set of group elements which serve to identify different points of this region. This set is finite, so we can start far enough out
along the fixed ray so that the finite set of translations involved doesn't take us out of the set of fixed points. 
We obtain a map from ${\rm Reg}(\widetilde{N}^{\ast} / N^{\ast} )$ to $\Tt$ which factors through
${\rm Reg}(\widetilde{D} / D)$, even though the resulting equivariant map on the whole universal covering 
$\widetilde{N}^{\ast}$ will not factor through $\widetilde{D}$. The different copies of ${\rm Reg}(\widetilde{N}^{\ast} / N^{\ast} )$
which make up $\widetilde{N}^{\ast}$ are related by paths which, even if they project to the same path in $\pi _1(D)$, act differently on $\Tt$
and this serves to create a divergence in the equivariant map. 

We now explain the
argument more carefully.

\begin{corollary}
\label{maponregion}
Suppose we can construct a Lipschitz piecewise $C^1$ map $\varphi : {\rm Reg}(\widetilde{D}/ D )\rightarrow \Tt$ with the property that for any 
$x,y\in {\rm Reg}(\widetilde{N}^{\ast} / N^{\ast} )$ and any $\alpha \in \pi _1(N^{\ast}/N)$ with 
$\alpha x = y$ then 
$$
\varphi (\sigma _{\rm Reg}y ) = \rho (\alpha ) \varphi (\sigma _{\rm Reg}x).
$$
With this, posing $\phi _{\rm Reg}(y):= \varphi (\sigma _{\rm Reg}(y))$ leads to an equivariant map $\phi : \widetilde{N}^{\ast} \rightarrow \Tt$
which has finite energy on the fundamental region. 
\end{corollary}
{\em Proof:}
Use Lemma \ref{inorderto}. The criterion of that lemma is exactly the same as the hypothesis here. 
Note that if $\phi _{\rm Reg}$ factors through $\sigma _{\rm Reg}$ then it automatically sends projection-equivalent points to 
the same point, since $\sigma$ has this property. 
\eop

Define the subsets of group elements corresponding to boundary identifications of the fundamental regions:
$$
\Aa := \{ \alpha \in \pi _1(N^{\ast}), \; \exists x,y\in {\rm Reg}(\widetilde{N}^{\ast} / N^{\ast} ), \; \alpha x = y \} ,
$$
and
$$
\Bb := \{ \beta \in \pi _1(D), \; \exists x,y\in {\rm Reg}(\widetilde{D}/ D ), \; \beta x = y \} .
$$
Both of these subsets are finite, and we have a map $\Aa \rightarrow \Bb$. 

Let $\Aa \times _{\pi _1(D )}\Aa$ be the set of pairs $(\alpha , \beta )\in \Aa \times \Aa$ with 
 $\sigma _{\ast}(\alpha ) =\sigma _{\ast} (\beta )$ in $\pi _1(D )$. 

\begin{lemma}
\label{ABsurj}
The map $\Aa \rightarrow \Bb$ is surjective, in other words $\Bb$ is the image of $\Aa$ in $\pi _1(D)$ and
it is the quotient of $\Aa$ by the equivalence relation $\Aa \times _{\pi _1(D )}\Aa$.
\end{lemma}
{\em Proof:}
Look at the choice of ${\rm Reg}(\widetilde{N}^{\ast} / N^{\ast} )$ described above Lemma \ref{regionproperties}. 
We first claim that the map 
$\sigma _{\rm Reg} : {\rm Reg}(\widetilde{N}^{\ast} / N^{\ast} )\rightarrow {\rm Reg}(\widetilde{D}/ D )$
is surjective. To see this, note that $N^{\ast} \rightarrow D$ is surjective, so 
$$
N^{\ast} \times _D{\rm Reg}(\widetilde{D}/ D )\rightarrow {\rm Reg}(\widetilde{D}/ D )
$$ 
is surjective. On the other hand, any point $N^{\ast} \times _D{\rm Reg}(\widetilde{D}/ D )$ lifts
to a point of ${\rm Reg}(\widetilde{N}^{\ast} / N^{\ast} )$ by construction of this latter.  This proves the claim.

Now we continue with the proof of the lemma. 
Suppose $\beta \in \Bb$ with $u,v\in {\rm Reg}(\widetilde{D}/ D )$ such that $\beta u = v$. Choose $x\in {\rm Reg}(\widetilde{N}^{\ast} / N^{\ast} )$
such that $\sigma _{\rm Reg}(x)=u$. Choose $\beta ' \in \pi _1(N^{\ast})$ mapping to $\beta$, and put $y':= \beta ' x$. We have
$\widetilde{\sigma}(y')=\beta u = v$. Thus,
$$
y'\in \widetilde{\sigma}^{-1}({\rm Reg}(\widetilde{D}/ D )) = Q^{\ast}.
$$
By choice of ${\rm Reg}(\widetilde{N}^{\ast} / N^{\ast} )$ as a fundamental region for the action of $\ker (\sigma _{\ast})$ on 
$Q^{\ast}$, there is a group element $g\in \ker (\sigma _{\ast})$ such that $y:= g\cdot y' \in {\rm Reg}(\widetilde{N}^{\ast} / N^{\ast} )$.
Now, putting $\alpha := g\beta '$ we get $\alpha x = g\beta ' x = g y' = y$, so $\alpha \in \Aa$ and it maps to $\beta$. 
\eop

\begin{lemma}
\label{maponregion2}
Suppose we can construct a map $\varphi : {\rm Reg}(\widetilde{D}/ D )\rightarrow \Tt$ with the property that for any 
$u,v\in {\rm Reg}(\widetilde{D}/ D )$ and any $\alpha \in \Aa$ such that 
$\sigma _{\ast}(\alpha ) u = v$, then $\varphi (v) = \rho (\alpha ) \varphi (u)$. Then $\varphi$ satisfies the criterion of
Corollary \ref{maponregion}. 
\end{lemma}
{\em Proof:}
Suppose $x,y\in {\rm Reg}(\widetilde{N}^{\ast} / N^{\ast} )$ and $\alpha \in \pi _1(N^{\ast}/N)$ with 
$\alpha x = y$. Note by definition $\alpha \in \Aa$. Put $u:= \sigma _{\rm Reg}(x)$ and $v := \sigma _{\rm Reg}(y)$. We have 
$\sigma _{\ast}(\alpha ) u = v$. By the hypothesis of this lemma, $\varphi (v) = \rho (\alpha ) \varphi (u)$, and in view of
the definition of $u$ and $v$ that gives exactly the condition of Corollary \ref{maponregion}. 
\eop

Note that $D$ is obtained from ${\rm Reg}(\widetilde{D}/ D )$ by glueing together points $u,v$ whenever
there is $\beta \in \Bb$ such that $\beta u = v$. 

\begin{corollary}
\label{maponregion3}
In order to apply Corollary \ref{maponregion} to construct a finite-energy initial map, it suffices to construct
a map $\varphi : {\rm Reg}(\widetilde{D}/ D )\rightarrow \Tt$ with the property that for any two points
$u,v\in {\rm Reg}(\widetilde{D}/ D )$ which are equivalent by a group element $\beta \in \Bb \subset \pi _1(D )$,
then for any lifting $\alpha \in \Aa$ with $\sigma _{\ast}(\alpha )=\beta$ we should have 
$\varphi (v) = \rho (\alpha ) \varphi (u)$.
\end{corollary}
{\em Proof:}
Given that the condition here makes sense by Lemma \ref{ABsurj}, the statement is the same as that of Lemma \ref{maponregion2}.
\eop

We have to look carefully at the standard procedure for creating an initial equivariant map into a contractible 
space, because $\pi _1(D )$ doesn't actually act naturally on $\Tt$.

\begin{proposition}
\label{maponregion4}
Suppose we have a nested sequence of subsets 
$$
R_{-1} \subset R_0\subset R_1 \subset \cdots \subset R_d \subset \Tt
$$
where $d={\rm dim}_{\rr}(D )$, such that each $R_k$ is contractible,
such that $\rho (\Aa \times R_k) \subset R_{k+1}$ and such that  
$$
\rho (\alpha )t = \rho (\beta )t,  \; \; \forall \, t\in R_k, \; \forall \, (\alpha , \beta )\in \Aa \times _{\pi _1(D )}\Aa .
$$
Then we can construct an initial map $\varphi$ as required by the previous Corollary \ref{maponregion3}.
\end{proposition}
{\em Proof:}
Let ${\rm Reg}_k(\widetilde{D}/ D )$ be the $k$-skeleton of the fundamental region. Define inductively maps
$\varphi _k : {\rm Reg}_k(\widetilde{D}/ D )\rightarrow R_k$ as follows. For $k=-1$ the domain is empty and the map
is trivially defined. For any $k$, assuming the map is defined on the $k-1$-skeleton, then choose a subset of $k$-simplices 
of ${\rm Reg}_k(\widetilde{D}/ D )$ mapping isomorphically to the set of $k$-simplices of $D$. 
By contractibility of $R_{k-1}$ we can map these simplices into $R_{k-1}$ arbitrarily, given the map $\varphi _{k-1}$ which
we already know on their boundaries. The remaining 
$k$-simplices of ${\rm Reg}_k(\widetilde{D}/ D )$ are mapped into $\Tt$ in a unique and well-defined way
using the group translation condition for $\varphi$. Note that
any $k$-simplex of ${\rm Reg}_k(\widetilde{D}/ D )$ is related to one in the original set by a unique group element
because the group $\pi _1(D )$ acts freely on $\widetilde{D}$, and furthermore this group element is
in $\Bb$ by the definition of the latter. It lifts to an element of $\Aa$ by Lemma \ref{ABsurj}.  Thus, with our condition that
$\rho (\Aa \times R_{k-1})\subset R_k$ we get that these maps send all of ${\rm Reg}_k(\widetilde{D}/ D )$
into $R_k$. This completes the inductive construction of $\varphi$.
\eop

We have now reduced to a purely algebraic problem. Recall that 
the centralizer $Z$ is the subgroup of unipotent upper-triangular matrices, its normalizer is $NZ=P$ the group of upper-triangular
matrices in $PSL(2,K)$, and $NZ/Z \cong \Gm$ is the diagonal torus. 

\begin{lemma}
\label{fixedrays}
Let $P$ be the group of upper-triangular matrices in $PSL(2,K)$, with projection $\xi : P \rightarrow \Gm$ to the
diagonal torus. Suppose we have a finite subset $A\subset P$. For any $d$
there is a sequence of contractible real rays $R_{-1} \subset \cdots \subset R_d \subset \Tt$, such that
$A\cdot R_{k}\subset R_{k+1}$ and 
$$
\forall \alpha , \beta \in A, \xi (\alpha ) = \xi (\beta ) \Rightarrow \forall t\in R_k, \alpha \cdot t = \beta \cdot t.
$$
\end{lemma}
{\em Proof:}
Let 
$$
D:= \{ (\alpha \beta ^{-1}) \in Z, \;\;\; \mbox{for} \;\; \alpha , \beta \in A \;\; \mbox{with} \;\; \xi (\alpha ) = \xi (\beta ) \} .
$$
This is a finite subset of $Z$. We can rewrite the conditions of the lemma by saying that we want rays $R_k$ for $-1\leq k \leq d$ such that
$R_k$ are fixed by elements of $D$, and $A\cdot R_k \subset R_{k+1}$. 

If we consider the standard basis $\{ {\bf e}_1, {\bf e}_2\}$ for $K^2$, it corresponds to two ends of $\Tt$ which are the two parabolic subgroups of upper
and lower triangular matrices. The real line joining these two ends consists of the segments joining together all adjacent lattices of the form 
$M_q := \langle t^q {\bf e}_1 , {\bf e}_2\}$ where $t\in K$ is the uniformizing parameter. The group of points of $Z$ is just $K$, and an element $a\in K$
acts on the basis elements by ${\bf e}_1\mapsto {\bf e}_1$, ${\bf e}_1\mapsto a{\bf e}_1 + {\bf e}_2$. 
There is an upper bound on the valuation of elements of the finite subset $D\subset K$. This leads to $q'$ such that the lattice $M_q$ is fixed by $D$ for any 
$q\leq q'$. We may consider any $q_k \leq q_0$ and define $R_k$ to be the ray of segments  joining the lattices $M_q$ for $q\leq q_k$.
With this choice we get the first condition that $R_k$ consists of fixed points of $D$. 

To get the second condition, notice that the diagonal group $\Gm \cong K^{\ast}$ acts on our real line by translations through the valuation map 
$\nu : K^{\ast} \rightarrow \zz$. An element $b\in K^{\ast}$ sends $M_q$ to $M_{q+\nu (b)}$. Again, since $\xi (A)\subset K^{\ast}$ is a finite subset,
there is a bound on the valuations of elements. Thus, there is $n$ such that $\nu (b)\leq n$ for $b\in \xi (A)$. If we choose $q_{k+1}= q_k + n$ for example, then 
the second condition $A\cdot R_k \subset R_{k+1}$ will hold. 

To complete the proof, just choose $q_{-1}\leq q'-(d+1)n$, then $q_k = q_{-1} + (k+1)n \leq q'$ for every $k=-1,\ldots , d$ and we obtain the required choice of rays.
\eop

This suffices for our construction of $\varphi$.

\begin{corollary}
\label{maponregion5}
In the unipotent case, 
putting $A := \rho (\Aa )\subset P\subset PSL(2,K)$ and  $R_k:= R^{d-k}$ we obtain regions 
satisfying the criterion of Proposition \ref{maponregion4}, hence going back to Corollary 
\ref{maponregion} we get the desired finite-energy initial map.
\end{corollary}
{\em Proof:}
By the map of exact sequences in Lemma \ref{toNZ}, the points of $\Aa \times _{\pi _1(D)}\Aa$ map to $A\times _{\Gm}A$. 
Thus, the subsets provided by Lemma \ref{fixedrays} serve to satisfy the criterion of Proposition \ref{maponregion4}.
\eop

{\em Remark:}---
In the above construction for the unipotent case,  we could have replaced the divisor $D$ by its incidence complex $\Sigma$, cf the remark on p. \pageref{sigmaremark}. 
There is a projection $D\rightarrow \Sigma$ and composing gives a projection $N^{\ast} \rightarrow \Sigma$. We could have constructed the initial 
map in a way which factors through this projection over the fundamental region. This is special to the unipotent case: in the
split torus case it looks like we need to consider the projection to $D$. 

\subsection{The pluriharmonic map}

We sum up the above constructions and apply Gromov-Schoen to deduce the existence of an equivariant harmonic map
of finite energy on the fundamental region.  In this subsection, we no longer assume that $D$ is connected.

Our now arbitrary basepoint $x\in X$ corresponds to a universal covering $f:\widetilde{X}\rightarrow X$. 

Let $D_{(i)}$ denote a connected component of the divisor $D$. It may be a union of several
irreducible components. Let $N_{(i)}$ denote a tubular neighborhood of $D_{(i)}$ in
$\overline{X}$ and let $N_{(i)}^{\ast} := N_{(i)}\cap X$ denote the ``punctured'' tubular 
neighborhood. Define the fundamental group of $N_{(i)}^{\ast}$ inside $\pi _1(X,x)$ 
by choosing a basepoint $x_{(i)}$  in $N_{(i)}^{\ast}$ and joining it to $x$ by a path.
Denote the subgroup of paths which consist of going from $x$ to $x_{(i)}$ along our path, then arbitrarily
inside $N_{(i)}^{\ast}$, and again along our path back to $x$, by $\pi _1(N_{(i)}^{\ast}) \subset \pi _1(X,x)$.

Let $\widetilde{N}_{(i)}^{\ast}$ denote the universal covering determined by the basepoint $x_{(i)}$.
The inverse image $f^{-1}N^{\ast}_{(i)}$ inside $\widetilde{X}$ is a disjoint union of copies of
$\widetilde{N}_{(i)}^{\ast}$. In order to define a $\pi _1(X,x)$-equivariant map starting from this disjoint union,
it is equivalent to defining a $\pi _1(N^{\ast}_{(i)}, x_{(i)})$-equivariant map starting from 
$\widetilde{N}_{(i)}^{\ast}$. 

Our preceding discussion where we focalized on a single connected component of $D$, allows us to choose 
a  $\pi _1(N^{\ast}_{(i)}, x_{(i)})$-equivariant map $\widetilde{N}_{(i)}^{\ast}\rightarrow \Tt$ which has finite energy on a fundamental region.
This gives a $\pi _1(X,x)$-equivariant map $\widetilde{X}\supset f^{-1}N^{\ast}_{(i)} \rightarrow \Tt$.
Do this for each connected component $D_{(i)}$. 
By contractibility of $\Tt$ we can extend these maps to a $\pi _1(X,x)$-equivariant map $\phi : \widetilde{X}\rightarrow \Tt$. The procedure
was recalled in detail in the proof of Proposition \ref{maponregion4} above. The complement of the union of translates of the
fundamental regions for the $\widetilde{N}^{\ast}_{(i)}$ is relatively compact
in the fundamental region of $\widetilde{X}$, so the map $\phi$ has finite energy on the fundamental region.

The action of $\pi _1(X)$ on $\Tt$ is reductive in the sense of \cite{JostZuo2}, 
because of our hypothesis that the Zariski closure of the image of
$\rho$ is $SL(2,K)$.

By the general theory of Gromov-Schoen \cite{GromovSchoen}
our initial map can be replaced by \revision a harmonic $\rho$-equivariant map 
$\Phi : \widetilde{X}\rightarrow \Tt$. 

The harmonic map is in fact pluriharmonic, by Gromov-Schoen \cite{GromovSchoen} 
Theorem 7.3.

The hypothesis that the image of $\rho$ is not contained in a compact subgroup, implies that
the pluriharmonic map is not constant, because if it were constant then its image point would be
fixed by $\rho$ and the subgroup of elements fixing a point of $\Tt$ is compact. We get:

\begin{proposition}
\label{pluriharmonicexists}
Under Hypothesis \ref{unboundedZD}, there exists a nonconstant pluriharmonic $\rho$-equivariant map from the universal
cover of $X$ to the Bruhat-Tits tree $\Tt$. With respect to the Poincar\'e-like 
complete finite-volume metric on $X$, the pluriharmonic map has finite energy.
\end{proposition}
\eop

\subsection{Fixed points of normal subgroups}

We close out this section with a study of the following situation. Suppose a finitely presented group
$\Gamma$ acts on $\Tt$ by a representation
$\rho : \Gamma \rightarrow PSL(2,K)$, and suppose $\Upsilon \subset \Gamma$ is a normal subgroup.
Let $R\subset \Tt$ be the set of fixed points of $\Upsilon$, that is the
set of points $y\in \Tt$ such that $uy = y$ for all $u\in \Upsilon$.
Note that $R$ is convex in the sense that if $y,z\in R$ then the unique path joining $y$ to $z$ 
is contained in $R$. Thus $R$ is a subtree of $\Tt$. 

The fact that $\Upsilon$ is a normal subgroup implies that $ \Gamma R \subset R$, indeed 
$$
\forall u\in \Upsilon \;\; \forall \gamma \in \Gamma \;\; \forall y\in R, \,\,  
u(\gamma y) = \gamma ((\gamma ^{-1} u \gamma ) y) = \gamma y.
$$ 
In particular we get an action of 
the quotient $\Gamma / \Upsilon$ on $R$. 

The subgroup $\Upsilon$ fixes the set of ends of $R$ (which are also ends of $\Tt$).

If $R$ has at least three distinct ends, this implies that the image of $\Upsilon$ in 
$\Tt$ is contained in the intersection of at  least three distinct normal subgroups. These normal subgroups
correspond to  distinct one-dimensional subspaces $L_i\subset K^2$ preserved by $\Upsilon$. 
Up to conjugation we may assume
that three of these subspaces are the two coordinate lines plus the diagonal. The subgroup of 
elements preserving these subspaces is the center $\{\pm 1\} \subset SL(2,K)$, that is the trivial subgroup of
$PSL(2,K)$. Thus, if $R$ has at
least three ends then the image of $\Upsilon$ is trivial.

On the other hand, the set of ends of $R$ remains invariant, although maybe acted upon, by the action of 
$\Gamma$. A subset consisting of one or two ends corresponds to a subgroup (either a parabolic subgroup or
a torus), so in either of these cases the image of $\Gamma$ cannot be Zariski-dense in $PSL(2,K)$.

Finally if $R$ has no ends, thus it  is compact, then $\Gamma$ fixes its barycenter so $\Gamma$
is contained in a compact subgroup. We can sum this up in the following lemma. 

\begin{lemma}
\label{normalsubgroup}
Suppose $\rho : \Gamma \rightarrow PSL(2,K)$ is a representation, and suppose $\Upsilon \subset \Gamma$
is a normal subgroup acting with a nonempty set of fixed points. Suppose that the image $\rho (\Gamma )$
is Zariski-dense and not contained in a compact subgroup of $PSL(2,K)$. Then the image of $\Upsilon$ 
is trivial. 
\end{lemma}
\eop

It will also be useful to have the invariance of our hypothesis under generically finite coverings. 

\begin{lemma}
\label{hypothesisinvariant}
Suppose $(X,\overline{X}, \rho )$ satisfies Hypothesis \ref{unboundedZD}.Then for any
generically surjective map $Z\rightarrow X$ and any smooth normal crossings compactification $Z\subset \overline{Z}$,
if we let $\rho _Z$ denote the pullback of $\rho$ to $Z$ and $Y$ the maximal open subset of $\overline{Z}$ on which
$\rho _Z$ extends to $\rho _Y$, then $(Y, \overline{Z},\rho _Y)$ satisfies Hypothesis 
\ref{unboundedZD}.
\end{lemma}
{\em Proof:}
Note first that $Y$ is still the complement of a divisor with normal crossings, which is the divisor $\overline{Z}-Z$ with
possibly some components removed. 

The local monodromy transformations on $Y$ come from those of $X$. If we have a small loop around a component of the
divisor at infinity in $Y$ then it goes to a loop in one of the neighborhoods of a multiple intersection point on $D$.
In such a neighborhood, the monodromy is an abelian representation whose generators are quasi-unipotent. 
A commuting
family of quasi-unipotent matrices can be simultaneously upper-triangularized, so any matrix in the
group they generate is also quasi-unipotent. This shows that the monodromy of $\rho_Y$ is
quasi-unipotent at infinity. 

The image of $\pi _1(Z)$ is of finite index in $\pi _1(X)$, call this subgroup $\Upsilon $. 
Furthermore $\rho _Y(\pi _1(Y)) = \rho _Z(\pi _1(Z)) = \rho (\Upsilon )$.
Let $\Upsilon '$ be the intersection of the conjugates of $\Upsilon$. It is a normal subgroup.  
Let $H$ be the Zariski closure of $\rho (\Upsilon ')$. It is a normal subgroup of $SL(2,K)$ since the image of $\rho$ is
Zariski dense,  so $H$ is either the whole
group or the center. If it were the center then the image of $\pi _1(X)$ in $PSL(2,K)$ would be finite, contradicting
Zariski density of $\rho$. Thus $H= PSL(2,K)$ which shows that $\rho _Y$ is Zariski dense. 

Suppose now that $\rho _Y(\pi _1(Y))$ were contained in a compact subgroup. Using the same group $\Upsilon '$ as in
the previous paragraph, we would get  $\rho (\Upsilon ')$ contained in a compact subgroup. By Lemma \ref{normalsubgroup}
the image $\rho (\Upsilon ')$ would have to be contained in the center, again contradicting Zariski density of $\rho$. 
This completes the verification of the conditions of Hypothesis \ref{unboundedZD}.
\eop

\section{Properties of the pluriharmonic map}
\label{sec-harmonic-properties}

In this section, we look further at the properties of the pluriharmonic mapping constructed above.
The case of trees was discussed explicitly in Gromov-Schoen \cite{GromovSchoen}, when $X$ is compact.
For the quasiprojective case, 
our situation is subsumed under the more general case of maps to buildings
considered by Jost and Zuo \cite{JostZuo2}. Since the case of maps to trees is considerably simpler in that there
is only a single spectral form to consider, it seems opportune to look at the details here. 

We keep Hypothesis \ref{unboundedZD} and look at the pluriharmonic map given by Proposition 
\ref{pluriharmonicexists}. Our goal is to obtain the factorization result Theorem \ref{mainfactorization} at the end. 

\subsection{The spectral form}

The pluriharmonic map yields spectral data which is a two-valued one-form $\theta$ on $X$.
This is easy to describe in the case of a map to a tree: the edges of $\Tt$  have a linear structure
and also metrics so we can identify them with intervals in $\rr$. This identification is well defined up
to translation and inversion, so the differential $df$ is well defined up to $\pm 1$. The differential is 
pluriharmonic, so
it is the real part of a holomorphic differential form $\theta$. It turns out, due to the regularity statements
of Gromov and Schoen \cite{GromovSchoen} which are local so they work equally well in the quasiprojective case,
that the symmetric square $Sym^2(\theta )$ is a well-defined and holomorphic section of $Sym^2(\Omega ^1_X)$.

Near a smooth point of the divisor at infinity, choose a coordinate system $(z_1,\ldots , z_n)$ such that
$D$ is given by $z_1=0$. Then \revision we can write
$$
Sym^2(\theta ) = \sum _{i\leq j} a_{ij}(z)dz_idz_j
$$
with the functions $a_{ij}(z)$ holomorphic on the complement of $z_1=0$. 

At any point $z\in X$, $Sym^2(\theta )$ is the square of a one-form $\theta (z)$ well-defined up to sign.
Furthermore, the energy is the integral of the square norm of $\theta (z)$, which in turn is the same as
$| Sym^2(\theta ) |$. We conclude that the integral (over our coordinate neighborhood) is finite: 
$$
\int \left| \sum _{i\leq j} a_{ij}(z)dz_idz_j \right| d\mu <\infty .
$$
Up to a bounded change in the metric on the base, we may assume that the metric is the product of the Poincar\'e
metric in $z_1$ with the euclidean metric on the $z_i$ for $i>1$. In particular the $dz_idz_j$ are
orthogonal and we get
$$
\int \left| a_{ij}(z)dz_idz_j \right| d\mu <\infty
$$
for all $i\leq j$. Recall that $|dz_i|= 1$ for $i>1$ and if we set $r:=|z_1|$
$$
| dz_1| = r |\log r|   
$$
and
$$
d\mu = \frac{dz_1\wedge d\overline{z}_1\wedge \cdots \wedge dz_n\wedge d\overline{z}_n}{r^2|\log  r|^2}.
$$

For showing meromorphy of $a_{ij}$ along $D$ and calculating their maximum order of poles, it suffices to
restrict to curves given by $z_i= p_i$ fixed for $i>1$. This is because if the integral is finite then it is finite
for almost all curves, and we have a Hartogs principle for separately almost everywhere meromorphic functions
\cite{Shiffman} \cite{Shiffman2} \cite{Siciak}.
The slices in the other directions are automatically almost everywhere holomorphic, indeed for any $z_1\neq 0$. 

Thus we are reduced to considering a function $a_{ij}^p(z_1):= a_{ij}(z_1,p_2,\ldots , p_n)$ of one variable,
holomorphic on a punctured neighborhood of the origin,
with the hypothesis that 
$$
\int \left| a_{ij}^p(u)\right| (r|\log  r| )^{m-2} dud\overline{u} < \infty
$$
where $m= 2$ if $i=j=1$, or $m=1$ if $i=1<j$ and $m=0$ if $i,j>1$. Here $r= |u|$. 

We are in a situation similar to Bell's removable singularities theorem, \cite[p. 687]{Bell}.

An argument using Cauchy's theorem on discs not meeting the origin shows that $a_{ij}^p$ is 
meromorphic at the origin.
To get the maximum order of poles, suppose 
$$
a_{ij}^p = b u^{-n} + \ldots 
$$
is the Laurent expansion with first coefficient $b \neq 0$. Then for $u$ near the origin 
$$
|a_{ij}^p (u)| \sim r^{-n}
$$
so we have
$$
\int _{0}^{1} r^{m-n-1} |\log r|^{m-2}dr < \infty .
$$
Making the change of variables $u=-\log r$ we get
$$
\int _{0}^{\infty} e^{(n-m)u}u^{m-2}du < \infty .
$$
We conclude that in the case $m=0$ we get $n\leq 0$, whereas for $m\geq 1$ we get $n<m$. 
Thus for $m=0,1$ we have $n\leq 0$ and for $m=2$ we have $n\leq 1$. In view of the definition of $m$ this
means that the $a^p_{ij}$ are holomorphic on $\overline{X}$ if either of $i$ or $j$ is different from $1$,
whereas $a^p_{11}$ has a pole of order $\leq 1$.

Going back to all of $\overline{X}$, we conclude by Hartogs \cite{Shiffman} \cite{Shiffman2} \cite{Siciak} 
that the $a_{ij}$ are meromorphic along smooth components of the divisor.
By the more classical Hartogs, they are meromorphic everywhere. Their orders of poles are bounded as described
previously. 

We can say what this bound means in invariant terms. Define the {\em half-logarithmic symmetric square}
to be the subsheaf denoted, for lack of better notation, 
$$
Sym ^2(\Omega ^1_{\overline{X}}; \frac{1}{2}\log D) \subset 
j_{\ast}(Sym ^2(\Omega ^1_X))
$$
as the subsheaf which is generated along smooth points of 
$D_i$ by $Sym ^2(\Omega ^1_{\overline{X}})$ plus $\frac{dz^2}{z}$ where $z$ is a coordinate of $D_i$. 
One can check that the subsheaf $Sym ^2(\Omega ^1_{\overline{X}} ; \frac{1}{2}\log D)$ 
is independent of the choice of coordinate. This defines the subsheaf at points of
$\overline{X}$ which are not normal crossing points of $D$. This defines a unique reflexive coherent sheaf on 
$\overline{X}$. By choosing local coordinates at a normal crossing point and noting that the generators of
$\Ff$ have a simple monomial form, we can see that $\Ff$ is actually locally free.  

Our estimate can be stated as follows: 

\begin{lemma}
\label{halfsymm}
With the above notations, for the finite energy pluriharmonic map the differential 
$Sym^2(\theta )$ extends to a holomorphic section of the bundle 
$Sym ^2(\Omega ^1_{\overline{X}}; \frac{1}{2}\log D) $ over $\overline{X}$
defined in the previous paragraph. 
\end{lemma}
\eop

\begin{lemma}
\label{pullbackF}
Suppose $p:\overline{Z}\rightarrow \overline{X}$ is a morphism from a smooth projective variety to 
$\overline{X}$, and suppose 
we have a divisor $D_Z\subset \overline{Z}$ which
is assumed to have normal crossings. Suppose that the inverse image of $D$ is contained in $D_Z$.
Thus if $Z:= \overline{Z}-D_Z$ then $p:Z\rightarrow X$. 
Define as above the half-logarithmic symmetric squares for $(\overline{X}, D)$ and
$(\overline{Z}, D_Z)$. Then the natural map 
$$
p^{\ast}Sym ^2(\Omega ^1_{X})\rightarrow Sym^2(\Omega ^1_{Z})
$$
extends to a map 
$$
p^{\ast} Sym ^2(\Omega ^1_{\overline{X}}; \frac{1}{2}\log D) 
\rightarrow Sym ^2(\Omega ^1_{\overline{Z}}; \frac{1}{2}\log D_Z) .
$$
\end{lemma}
{\em Proof:}
The sheaves on both sides of the desired map are vector bundles. 
Thus, to obtain the map we just need to construct it in codimension
$1$. In particular we can look at near a point on a smooth point $P\in D_Z$, where we have coordinates
$z_1,\ldots , z_n$ such that the divisor $D_Z$ is given by $z_1=0$ (with $P$ corresponding to $z_i=0$). 
The image $Q=p(P)$ in $X$ is a point which we may assume is on $D$, because if it were not then the
statement would be immediate. Thus we
have coordinates $x_1,\ldots , x_m$ on $X$ again centered at $Q$ and with $D$ given by $x_1\cdots x_k =0$. 
Write 
$$
x_i \circ p = u_i(z)z_1^{a_i},
$$
with $u_i$ nonvanishing along $D_Z$ at the point $P$. The inverse image of $D$ is given by the equation
$$
0=(x_1\cdots x_k)\circ p = u_1(z)\cdots u_k(z)z_1^{a_1+\ldots +a_k}.
$$
By hypothesis, the zero set of this equation is contained in $D_Z$, which near $P$ is just the
smooth divisor $z_1=0$. This implies that the zero sets of the $u_i(z)$ don't contain any other
irreducible components for $i=1,\ldots , k$. Therefore, $u_i(P)\neq 0$ for $i=1,\ldots ,k$. 
This implies that $a_i >0$ for $i=1,\ldots , k$ because our coordinates $x_i$ are centered at $Q$
so $x_i\circ p (P)=0$.
Of course we can't say anything for the $u_{k+1}\ldots u_n$ but these won't enter into the subsquent argument.

Now
$$
Sym ^2(\Omega ^1_{\overline{X}}; \frac{1}{2}\log D)
$$
is locally free, with generators near $Q$ of the form
$$
\frac{dx_i^2}{x_i}, \;\; 1\leq i\leq k,
$$
$$
dx_j^2,\;\; j>k,
$$
and 
$$
dx_idx_j, 1\leq i < j \leq n.
$$

Note that the pullback of $dx_i$ is
$$
p^{\ast}dx_i = a_iu_i(z)z_1^{a_i-1}dz_1 + \sum \frac{\partial u_i}{\partial z_j}z_1^{a_i}dz_j.
$$
This is a section of $\Omega ^1_{\overline{Z}}$, so the pullbacks of the generators $dx_j^2$ and
$dx_idx_j$  are in $Sym^2\Omega ^1_{\overline{Z}}$. We just have to check the pullback 
of the fractional generator for $i=1,\ldots , k$. We have
$$
p^{\ast} \left( \frac{dx_i^2}{x_i} \right) = z_1^{2a_i-3}a_i^2u_i(z)^2dz_1^2 + z_1^{2a_i -2}A
$$
where $A$ is a section of $Sym^2\Omega ^1_{\overline{Z}}$. The fact that $a_i \geq 1$ 
means that the term with $A$ is a section of $Sym^2\Omega ^1_{\overline{Z}}$. 
Similarly, the first term is a holomorphic multiple of $dz_1^2/z_1$. Thus,
$p^{\ast}(dx_i^2/x_i)$ is a section of $Sym ^2(\Omega ^1_{\overline{Z}}; \frac{1}{2}\log D_Z)$
as required.
\eop

\begin{corollary}
\label{pullbackPhi}
Suppose $\Phi$ is a finite energy pluriharmonic map on $X$ corresponding to a representation satisfying
Hypothesis \ref{unboundedZD}. Suppose $p:\overline{Z}\rightarrow \overline{X}$ is as in Lemma \ref{pullbackF}. Then
the pullback $\Phi \circ p$ is a pluriharmonic map on $Z$ satisfying the condition of Lemma \ref{halfsymm}, that
the symmetric square of its differential is a section of $Sym ^2(\Omega ^1_{\overline{Z}}; \frac{1}{2}\log D_Z)$. 
\end{corollary}
{\em Proof:}
Combine Lemmas \ref{halfsymm} and \ref{pullbackF}.
\eop

\subsection{Uniqueness}

The pluriharmonic map itself may not be unique. However, its differential is unique, and different pluriharmonic
maps are related by a geodesic homotopy which preserves the differential. 

Suppose $\Phi _0$ and $\Phi _1$ are two different pluriharmonic maps of finite energy equivariant for 
$\rho$. Define the {\em geodesic homotopy} $\{ \Phi _t \} _{t\in [0,1]}$ as follows. For any point $x\in \widetilde{X}$,
let $\Phi _t(x)$ denote the geodesic path joining $\Phi _0(x)$ and $\Phi _1(x)$.

\begin{lemma}
\label{geod}
The geodesic homotopy is a continuous, piecewise differentiable map
$$
\Phi : [0,1] \times \widetilde{X} \rightarrow \Tt .
$$
There is a dense open set $U\in \widetilde{X}$ such that on $[0,1]\times U$ the map $\Phi$ is 
harmonic in the second variable and linear in the first variable. On each connected component of
$U$ there is a choice of 
chart in $\Tt$ containing the images of all $\Phi _t(x)$, and in this chart the differential with respect
to the second variable satisfies 
$$
d\Phi _t (x) = t \, d\Phi _1(x) + (1-t) d\Phi _0(x).
$$
\end{lemma}
{\em Proof:}
Consider the function $G(x,y,t)$ which to points $x,y\in \Tt$ and $t\in [0,1]$ associates the
point in $\Tt$ on the geodesic from $x$ to $y$, with coordinate $t$. Thus $G(x,y,0)=x$ and 
$G(x,y,1)=y$. A geometric picture shows that it is Lipschitz, hence continuous. Furthermore if $x$
and $y$ are not vertices of $\Tt$ then near $(x,y)$ the function is bilinear onto the flat of $\Tt$ joining $x$
to $y$. In terms of a real coordinate $r$ on this flat, $G$ has the
form 
$$
rG(x,y,t) = (1-t)r(x) + tr(y).
$$
Our function $\Phi$ is defined by 
$$
\Phi (t,x) := G(\Phi _0(x),\Phi _1(x),t).
$$
The open set $U$ is the set of points $x$ such that $\Phi _0(x)$ and $\Phi _1(x)$ are not vertices of $\Tt$.
This is the complement of a one dimensional real analytic subvariety, in particular it is dense. 
On $U$, the map $\Phi$ is linear in $t$ and harmonic in $x$ because it is a linear combination of harmonic functions.
In terms of a real coordinate $r$ on the flat joining $\Phi _0(x)$ to $\Phi _1(x)$, we have
$$
r \Phi _t(x) = (1-t)r\Phi _0(x) + t r\Phi _1(x).
$$
This gives the formula for the differential.
\eop

\begin{corollary}
\label{energyconvex}
The energy of $\Phi _t$ satisfies
$$
\int _X |d\Phi _t|^2 = e_0 + e_1 t + e_2 t^2
$$
with 
$$
e_2 = \int _X |d\Phi _0 - d\Phi _1|^2.
$$
Here the $d\Phi _0$ and $d\Phi _1$ are taken using the charts of $\Tt$ as in the previous lemma, well defined
on connected components of $U$. 
\end{corollary}
{\em Proof:}
Plug in the formula 
$$
d\Phi _t (x) = t \, d\Phi _1(x) + (1-t) d\Phi _0(x)
$$
to get
$$
|d\Phi _t|^2 = |t \, d\Phi _1(x) + (1-t) d\Phi _0(x)|^2
$$
$$
=t^2 |d\Phi _1(x)|^2 + (1-t)^2|d\Phi _0(x)|^2 + 2t(1-t)d\Phi _0(x)d\Phi _1(x)
$$
$$
= E_0(x) + tE_1(x) + t^2 \left( |d\Phi _1(x)|^2 + |d\Phi _0(x)|^2 -2d\Phi _0(x)d\Phi _1(x)\right) .
$$
This gives
$$
\int _X |d\Phi _t|^2 = e_0 + e_1t + t^2\int _X |d\Phi _0 - d\Phi _1|^2.
$$
\eop

We get the following uniqueness statement for the differential of the pluriharmonic map. While the methods
do not seem to be related at all, there is a uniqueness statement for harmonic maps which has a somewhat
similar form in the recent preprint \cite[Lemma 3.1]{GelanderKarlssonMargulis}.

\begin{theorem}
\label{uniquenessdiff}
If $\Phi _0$ and $\Phi _1$ are two pluriharmonic maps of finite energy equivariant for the same representation
$\rho$, then with appropriate determinations of charts of $\Tt$ we have $d\Phi _0 = d\Phi _1$. The continuous maps
$\Phi _t$ in the geodesic homotopy relating $\Phi _0$ and $\Phi _1$
are in fact differentiable and pluriharmonic and their differentials are all the same.
\end{theorem}
{\em Proof:}
The condition that $\Phi _0$ and $\Phi _1$ are pluriharmonic implies that $e_2=0$ in the previous corollary.
This is because the pluriharmonic condition means that the maps are critical points of the energy functional,
so the quadratic polynomial $e_0 + e_1 t + e_2 t^2$ has critical points at $t=0$ and $t=1$. This implies that 
$e_2=e_1=0$. 
In particular, with the determinations of charts as in the lemma, we have $d\Phi _0 = d\Phi _1$,
hence $Sym ^2(\theta _0)=Sym ^2(\theta _1)$.
Over the open set $U$ then, the differentials of $\Phi _t$ are all the same. 
\eop

We don't have a fully worked-out example, but it seems clear that the maps themselves need not be unique,
rather there can be the latitude \revision for a small interval of linear translations within $\Tt$. This is related
to the fact that critical points of $\theta$ need not go to vertices of $\Tt$.

\begin{lemma}
\label{pullbackFE}
Suppose $f:Y\rightarrow X$ is a generically surjective 
map, suppose $\rho : \pi _1(X,x)\rightarrow SL(2,K)$ is a representation
satisfying \revision Hypothesis \ref{unboundedZD}, and suppose $\Phi : \widetilde{X}\rightarrow \Tt$ 
is a pluriharmonic map of finite energy equivariant for $\rho$. Then
$\Phi \circ \widetilde{f}$ is a pluriharmonic map of finite energy equivariant for $f^{\ast}(\rho )$.
\end{lemma}
{\em Proof:}
The hypothesis of \revision quasi-unipotence at infinity is satisfied by $f^{\ast}(\rho )$, and the condition that
$f$ be generically surjective implies that $f^{\ast}(\rho )$ is still Zariski dense.
Thus \revision from Proposition \ref{pluriharmonicexists}
we get the existence of an equivariant pluriharmonic map of finite energy 
$\Phi ': \widetilde{Z}\rightarrow \Tt$. Choose a generically finite
covering $Z\rightarrow Y$ such that the pullbacks of the
differentials $Sym ^2(\theta )$ for $\Phi$ and $Sym ^2(\theta ')$ for $\Phi '$ become single-valued on $Z$. 

We may assume that $Z$ is smooth and has a normal crossings compactification.

Let $\Phi _0$ and $\Phi _1$ denote respectively the pullbacks of $\Phi$ and $\Phi '$ to $Z$.
We claim that they are of finite energy. Let $Sym^2(\theta _0)$ and $Sym ^2(\theta _1)$ denote the squares
of their differentials. By our assumption on $Z$ these are both squares of differential forms $\alpha _0$ and
$\alpha _1$ respectively on $Z$. By Lemma \ref{pullbackPhi},
$$
\alpha _0^2, \alpha _1^2 \in H^0(Sym ^2(\Omega ^1_{\overline{Z}}; \frac{1}{2}\log D_Z)).
$$ 
This implies by a local calculation that 
$\alpha _0,\, \alpha _1$ extend to holomorphic forms on the
compactification of $Z$. In turn this means that $\Phi _0$ and $\Phi _1$ have finite energy. 

Now, Theorem \ref{uniquenessdiff} says that the squares of their differentials
are the same, in other words $Sym^2(\theta _0) = Sym ^2(\theta _1)$. 
This implies that $Sym ^2(\theta ') = f^{\ast}Sym ^2(\theta )$ on $Y$. As $|Sym ^2(\theta ')|$
is integrable on $Y$, the same holds for  $|f^{\ast}Sym ^2(\theta )|$ which means that $\Phi \circ  \widetilde{f}$ 
has finite energy.
\eop

\begin{corollary}
\label{thetasplitting}
There is a finite quasiprojective ramified covering $Z\rightarrow X$ such that
over $Z$, the pullback of $\Phi$ is a finite-energy pluriharmonic map,
and the determination of a direction of the edges of the tree can be made and the 
differential becomes a holomorphic section
$\alpha \in H^0(\overline{Z},\Omega ^1_{\overline{Z}})$. The 
pluriharmonic map is constant on the leaves of the foliation
defined by the real one-form $\Re \alpha$. 
\end{corollary}
\eop

The fact that the pluriharmonic map is not constant implies that the differential is not identically zero,
so $Sym^2(\theta )$ is not identically zero, and on the covering $Z$ the form $\alpha$ is nonzero.

\subsection{Geometry of the pluriharmonic map}
\label{sec-geom}

We now discuss in a preliminary way some of the real geometry of this picture. This subsection is optional. 

\begin{lemma}
In the situation of Corollary \ref{thetasplitting}, the periods of the real one-form $\Re \alpha$
are integers. The associated class $[\Re \alpha ]\in H^1(\overline{Z},\zz )$ corresponds to the map
$$
\varphi : \overline{Z} \rightarrow S^1 = \Tt '/ SL(2,K)'
$$
where $\Tt '$ is the subtree of $\Tt$ consisting of edges whose interior is touched by $\Phi$,
and we denote by $SL(2,K)'\subset SL(2,K)$ the subgroup of transformations preserving $\Tt '$ and preserving
the orientations of all edges of 
$\Tt '$. This  restricts to a map on the open set $Z$ and 
the representation $\rho$ on $\pi _1(Z)$ factors through the subgroup $SL(2,K)'$. 
\end{lemma}
{\em Proof:}
The choice of a determination of $\alpha$ corresponds to a choice of direction for the edges of $\Tt '$,
so $\rho :\pi _1(Z)\rightarrow SL(2,K)'$. Now we can form the quotient metric space $S^1 = \Tt '/ SL(2,K)'$
and the pluriharmonic map descends to give a map $\varphi$ well defined on $Z$. The periods of 
$\Re \alpha $ correspond to signed lengths of images of loops by $\varphi$ so they are clearly integers.
In particular the periods of $\Re \alpha$ on $\overline{Z}$ are integers so the Albanese map projects to 
a map $\varphi$ defined on the compactification $\overline{Z}$.
\eop

In the situation of the lemma, the map $\varphi$ factors as the composition of the Albanese map of $\overline{Z}$ with 
the linear map whose differential is $\Re \alpha$:
$$
Z\rightarrow Alb(\overline{Z})\rightarrow S^1.
$$
We can also do a ``real Stein factorization'' of the map $\varphi$ through a map
$$
Z\rightarrow G \rightarrow S^1
$$
where $G$ is a graph. Note in the case when $Z$ is noncompact, the graph $G$ may
be non-Hausdorff, which is what leads to monodromy around the points at infinity.

The graph $G$ is a version of the {\em Reeb graph} of a Morse function \cite{Reeb}. This notion came to our attention in
a talk of B. Falcidieno. It has been cited recently in the papers
\cite{CasallesEtAl}, \cite{ColeEtAl}, and \cite{EntovPolterovich}. 

There is a distinguished point on $S^1$ corresponding to the image of the vertices of
$\Tt '$. However, the vertices of $G$ might not go to this point. In the case where no vertex of $G$ goes to 
the distinguished point, then one can compose with small rotations of $S^1$ which then lift, and get 
a family of translations of the pluriharmonic map $\Phi$. This is why we said above that it looked like the
map was not necessarily unique. To get an example, one would have to get an example where the graph $G$ can
map to $S^1$ sending vertices to points other than the distinguished one.

It would be good to have more information on how the representation $\rho$ determines the covering $Z$,
the cohomology class $[\Re \alpha ]$, the map $\varphi$ and its factorization through the graph $G$, and the
lifting of $\varphi$ to the pluriharmonic map $\Phi$. Ideally one would like to be able to calculate everything
explicitly, whenever we are given $\rho$ by matrices.

\subsection{Lefschetz theory for the spectral foliation}

After pulling back to the covering $Z$ given by Corollary \ref{thetasplitting}, 
we have the following situation. We have an irreducible 
smooth quasiprojective variety $Z$ with smooth normal-crossings
compactification  $\overline{Z}$. The lift to $Z$ of the multivalued spectral form
$\alpha \in H^0(\overline{Z}, \Omega ^1_{\overline{Z}})$ is
a nonzero holomorphic one-form on the compactification. Let $A$ be the minimal abelian variety from which 
$\alpha$ comes (i.e. $A=Alb(\overline{Z})/B$ where $B$ is the maximal sub-abelian variety on which $\alpha$ 
vanishes). Let $Z'$ and $\overline{Z}'$ be the coverings of $Z$ and $\overline{Z}$
defined by fiber product with the universal covering $V\rightarrow A'$. 

Let $\Ff$ be the foliation defined by $\alpha$
on $Z'$. We have a map $Z'\rightarrow V$, and $\alpha$ corresponds to a linear form on $V$,
which we denote $g:V\rightarrow \cc$. For $t\in \cc $ we get a linear subspace $P_t\subset V$
and $g^{-1}(t)\subset \overline{Z}'$ is its pullback to $\overline{Z}'$. Let
$$
D':=\overline{Z}' - Z'
$$
be the covering of the complementary divisor.

\begin{lemma}
\label{commsg}
In the above situation,
the coverings $\overline{Z}'$ and $Z'$ are connected, and
$\pi _1(Z')$ contains the commutator subgroup of $\pi _1(Z)$. 
\end{lemma}
{\em Proof:}
The covering $Z'\rightarrow Z$ is an infinite Galois covering whose Galois group is the quotient
of $\pi _1(Z)$
given by the surjective morphism
$\pi _1(Z)\rightarrow \pi _1(A)$. In particular, $Z'$ is connected and the Galois group of $Z'/Z$ is abelian, hence 
$\pi _1(Z')$ contains the commutator subgroup of $\pi _1(Z)$.
\eop

The following is the basic Lefschetz theorem which we need. It is a generalization to the quasiprojective
case of the statement of \cite{lefschetz} which was done for 
the projective case of the present project. 

\begin{theorem}
\label{lef-th}
Suppose that the complex dimension of the image of $Z$ in $A$ is at least two. Then 
if $t \in \cc$ is a general value (i.e. not one of a countable number of singular values),
the inverse images $g^{-1}(t)$ are connected and their fundamental groups surject onto $\pi _1(Z')$.
\end{theorem}
{\em Proof:} We treat first a simplified special case which is then easily seen to give the general case.

{\em Special case:} Suppose that the following additional conditions hold: that $dim (Z)=2$,
the morphism from $\overline{Z}$ onto its image $\overline{Z}_A\subset A$ is a resolution of singularities,
that there is a Weil divisor $D_A\subset \overline{Z}_A$ containing the singular set of
$\overline{Z}_A$, such that $D$ is the inverse image of $D_A$ in
$\overline{Z}$, and that the map $Z\rightarrow \overline{Z_A}-D_A$ is an isomorphism. In particular we may 
consider this last isomorphism as an equality and view $Z$ as an open subset of the smooth points of
$\overline{Z_A}\subset A$. 

Now as explained in \cite{lefschetz}, in this situation the only zeroes of $\alpha$ on 
$Z$ lie over points in $A$. We can also look at zeros of the restriction of $\alpha$ to $D$. Since $D$ is a curve,
there are two cases: a zero of $\alpha |_D$ is either a point, or else it is a union of components of $D$. However
if $D_i$ is a component of $D$ such that $\alpha |_{D_i}=0$  then the albanese variety of $D_i$
would be an abelian subvariety
on which $\alpha$ vanishes. By the construction of $B$ and $A$, this means that $D_i$ maps to a point in $A$.
We conclude that there is a finite set of points $p_j\in \overline{Z}_A$ containing the zeros $\alpha |_Z$
and the images of the zeros of $\alpha |_D$. 

We can now apply the same argument as in \cite{lefschetz}. We don't give the details of 
the choice of nested neighborhoods and so forth necessary to make things
rigorous: this follows the same principles as in \cite{lefschetz}.

Fix a collection of relatively compact neighborhoods of these points
$\overline{N}_A(p_j)\subset \overline{Z}_A$. 

Away from these neighborhoods we can choose vector fields of bounded length on 
$\overline{Z}$, preserving $Z$ i.e. tangent to $D$, whose images by $\alpha$ are 
constant (real and imaginary fields on $\cc$). 

Following these vector fields outside of the given neighborhoods, 
we can lift a given homotopy in $\cc$ in a way which takes 
us only within a bounded region of
the noncompact $\overline{Z}'$. This allows us to solve the Lefschetz problem outside of the 
neighborhoods, see \cite{lefschetz}.

Thus we can restrict our attention to one of the neighborhoods $\overline{N}_A=\overline{N}_A(p_j)$. 
Let $\overline{N}_1$ denote its inverse image in $\overline{Z}$,
and let $N_1=\overline{N}_1\cap Z = \overline{N}_1-\overline{N}_1\cap D$ be its intersection with $Z$.
This might be disconnected because the image $\overline{Z}_A$ is not necessarily normal. Let
$N$ be a connected component of $N_1$ and $\overline{N}$ the closure of $N$ in $\overline{N}_1$.

The inverse image of of $\overline{N}$ in $\overline{Z}'$ 
is a disjoint union of isomorphic
copies of the neighborhood $\overline{N}$ itself, and each copy denoted $\overline{N}'$
is relatively
compact in $\overline{Z}'$ (i.e. the covering $Z'/Z$ does not induce any
nontrivial coverings of the neighborhoods).
Since we are in a bounded region of $\overline{Z}'$, only a finite number of these come into play.
We just have to solve the Lefschetz problem in one of the connected neighborhoods $N'\subset \overline{N}'$.

Notice that the map $\overline{N}'\rightarrow \overline{N}_A(p_j)$ is a resolution of singularities of one of the
irreducible components of the target. Any positive dimensional exceptional fiber is contained in the interior
of $\overline{N}'$ in other words does not go near the boundary; this is because the exceptional components map
to the point $p_j\in \overline{N}_A(p_j)$.

Here the situation becomes classical. We can look at integration of $\alpha$ as defining
a holomorphic map $g$ to a disc which is a fibration over the punctured disk. 
In fact this comes from a function $g_A:\overline{N}_A(p_j)\rightarrow \Delta$
which we normalize so that $g(p_j)=0$. Over the punctured disc $\Delta ^{\ast}$ the
map of triples 
$$
(\overline{N}', D\cap \overline{N}', N')_{\Delta ^{\ast}}\rightarrow \Delta^{\ast}
$$
is topologically locally trivial (we arrange the boundary of the neighborhood in a good way so that this holds).

Note furthermore that $D$ intersects the fiber $g^{-1}(0)$  only in points lying over $p_j$. 
This is because any component of $D$ upon which the differential $dg = \alpha$ restricts to zero, maps to 
an isolated point in $A$ (which in the present case must be $p_j$).
Given our hypothesis that $D_A$ contains the singular points of $\overline{Z}_A$, this means in particular
that the fiber $g_A^{-1}(0)$ intersects the singular set of $\overline{N}_A(p_j)$ only at $p_j$.

The Lefschetz question is now localized as follows:
given $p$ and $q$ in one fiber $g^{-1}(t)$ of the map $g:N'\rightarrow \Delta$ with $t\neq 0$, 
and a path $\gamma$ from $p$ to $q$ inside the total space of the neighborhood $N'$,
is $\gamma$ homotopic to a path which stays in the fiber?  

The main tool for doing this is the observation that there is at least one component of the central fiber
$g^{-1}(0)$ which isn't contained in $D$, and isn't a multiple fiber. Indeed, 
the fiber $g_A^{-1}(0)$ contains a positive dimensional component in each irreducible component of
$\overline{N}_A(p_j)$, by  purity of the dimension of a subset defined by one equation. As noted above,
away from $p_j$ this corresponds to a curve in our open subset $N'\cong N$. As $g^{-1}(0)$ is the
inverse image in $\overline{N}'$ of $g_A^{-1}(0)\subset \overline{N}_A(p_j)$, this means that $g^{-1}(0)$ 
contains a component which isn't contained in $D$ as claimed.  

Using this fact we can treat the local Lefschetz question.
Let $\xi$ be a path starting from $p$ and going around a non-multiple component of
the central fiber. 
Let $N^{\ast}$ denote the complement of $g^{-1}(0)$ in $N'$. Then $Path (N'; p,q)$
is obtained from $Path (N^{\ast }; p,q)$ by adjoining at least the relation that 
$\xi$ becomes trivial. Now the fibration 
$$
g: N^{\ast } \rightarrow \Delta ^{\ast}
$$
yields an exact sequence
$$
1 \rightarrow Path (g^{-1}(t); p,q) \rightarrow Path (N^{\ast }; p,q) \rightarrow \zz \rightarrow 0,
$$
and multiplication by $\xi$ via the action of $\pi _1(N^{\ast},p)$ on 
$Path (N^{\ast }; p,q)$ projects to translation by $1$ in $\zz$ (this is exactly the statement that 
$\xi$ is a loop around a non-multiple component; if it had gone around a component of multiplicity $r$ then
it would have projected to translation by $r$). 
We can now prove the Lefschetz property. The path $\gamma$ can be moved away from the central fiber
so it comes from an element of $Path (N^{\ast }; p,q)$; multiplying by an appropriate power of
$\xi$ yields a path in $Path (g^{-1}(t); p,q)$, and since $\xi$ is homotopically trivial in 
$N'$ this means that $\gamma$ was homotopic in $N'$ to an element of $Path (g^{-1}(t); p,q)$.

Adorning the argument with the niceties of \cite{lefschetz}, we are done in the special case $dim (Z)=2$. 

{\em General case:} We now return to the general situation of the statement of the theorem and show how to
obtain the conclusion, once we know the special case of dimension $2$ which has been treated already.

Suppose $x',y'\in g^{-1}(t)\subset Z'$ are two points, and suppose $\gamma$ is a path from $x'$ to $y'$ in $Z'$.
We want to show that $\gamma$ is homotopic to a path inside $g^{-1}(t)$. Let $x$ and $y$ denote the images of
$x'$ and $y'$ in $Z$. Choose a complete intersection of hyperplane sections $\overline{Y}\subset \overline{Z}$,
such that $x,y\in Y:= \overline{Y}\cap Z$, with $dim (Y)=2$. Choose $\overline{Y}$ 
of sufficiently high multidegree and
general subject to these conditions (although in the case
$dim (Z)=2$ we just have $\overline{Y}=\overline{Z}$).
In particular, the image of $Y$  in $A$ has dimension $2$ because by hypothesis the image of $Z$ has dimension 
$\geq 2$.  The map $\pi _1(\overline{Y})\rightarrow \pi _1(A)$ is still surjective so the covering $Y':= Y\times _ZZ'$
is connected as in Lemma \ref{commsg}. Also, the fundamental group of $Y$ surjects onto that of $Z$,
so $\pi _1(Y')\rightarrow \pi _1(Z')$ is surjective. Now $x',y'\in Y'$, and
$\gamma$ is homotopic to a path in $Y'$.
Thus, it suffices to prove that $g^{-1}(t)\cap Y'$ is connected and has fundamental
group surjecting to $\pi _1(Y')$, for then we will be able to make $\gamma$ homotopic to a path from $x'$ to $y'$
in  $g^{-1}(t)\cap Y'$. 

For the proof for $Y'$, we can replace $Y$ by a smaller open subset
without loss of generality, even if $x'$ or $y'$ is no longer in this smaller open subset. This is because 
connectedness and surjectivity of the map on fundamental groups are properties which are stable under passing between 
Zariski open subsets and supersets. Finally, we can also change the compactification $\overline{Y}$ without
changing the statement. Now by going to a small enough open subset $Y$ and choosing a normal crossings compactification
$\overline{Y}$, we can get to the special case which was treated in the first part of the proof. This completes
the proof in the general case. 
\eop

We now apply the above result to the spectral data for our pluriharmonic map.
Look first over $Z$. The differential of the pluriharmonic map is $\pm \alpha$ for a well-defined one-form $\alpha$. 
The restriction of the harmonic
map to any of the leaves $g^{-1}(t)$ of the resulting foliation, is constant. Therefore if
$\gamma \in \pi _1(Z)$ is a path which has a representative whose lift to the universal
covering stays within a leaf, then the action of $\gamma$ fixes the image point of the leaf.

\begin{lemma}
\label{dimone}
The dimension of the image of $Z$ in $A$ is one.
\end{lemma}
{\em Proof:}
Suppose to the contrary that the dimension of the 
image of $Z$ in $A$ is at least two. Then by the previous theorem \ref{lef-th},
for a generic value of $t$ we have the fiber $g^{-1}(t)\subset Z'$ which is connected and
has fundamental group surjecting to $\pi _1(Z')$. Consider the pullback 
$$
\widetilde{Z}_t:= g^{-1}(t)\times _{Z'}\widetilde{Z}.
$$
It is a union of possibly infinitely many connected components denoted $\widetilde{Z}_{t,k}$. 

Let
$\Phi: \widetilde{Z}\rightarrow \Tt$ be the pluriharmonic map. The differential of $\Phi$ is constant
along the components $\widetilde{Z}_{t,k}$ by construction of $g$ as the integral of $\alpha = \pm d\Phi$.
Choose a component $\widetilde{Z}_{t,k}$, which therefore maps to a single point $q\in \Tt$
under $\Phi$. If $x\in \widetilde{Z}_{t,k}$
is a point, and if $\gamma \in \pi _1(Z')$ then the translate is again in the same 
component: 
$$
\gamma \cdot x \in \widetilde{Z}_{t,k}
$$
In particular, $\Phi (\gamma \cdot x) = q = \Phi (x)$. The equivariance property of $\Phi$ tells us that
$\Phi (\gamma \cdot x) = \rho (\gamma )\Phi (x)$, so we get $\rho (\gamma )q = q$.
In other words, under the assumption that the dimension of the 
image of $Z$ in $A$ is at least two, applying Theorem \ref{lef-th} we conclude that 
$\pi _1(Z')$ has a fixed point in $\Tt$ under the action $\rho$.

By Lemma \ref{commsg}, the 
commutator subgroup 
the commutator subgroup $\Upsilon \subset \pi _1(Z)$ is contained in $\pi _1(Z')$. 
Thus, the assumption that the dimension of the 
image of $Z$ in $A$ is at least two would imply that the commutator subgroup $\Upsilon$ acts with a
fixed point on $\Tt$. By Lemma \ref{normalsubgroup} this would mean that $\Upsilon$ maps to the center of 
$SL(2,K)$. But that would imply that the representation into $PSL(2,K)$ factors
through the abelianization $\pi _1(X,x)^{\rm ab}$ and hence to a commutative subgroup scheme,
contradicting Zariski-denseness of $\rho$.  

We conclude that the dimension of the image $Z_A$ cannot be $\geq 2$.
Since $\alpha$ is nontrivial, we get $dim (Z_A)=1$. 
\eop

Putting all of the above discussion together we obtain the following theorem, see
\cite{Beauville} \cite{Campana94} \cite{CullerShalen} \cite{Delzant2} \cite{Eyssidieux}
\cite{Gromov} \cite{GromovSchoen} \cite{JostZuo2} \cite{Katzarkov} \cite{Klingler} 
\cite{Kollar} \cite{NapierRamachandran3} \cite{Siu} \cite{Zuo}. 

\begin{theorem}
\label{mainfactorization}
Suppose 
$X$ is a quasiprojective variety; $K$ is a complete local field with finite residue
field; $\rho : \pi _1(X,x)\rightarrow SL(2,K)$ is a representation whose monodromy transformations
at infinity are quasi-unipotent; the image of $\rho$ is Zariski-dense; and the image is not contained in any
compact subgroup. Then $\rho$ factors through a DM-curve.
\end{theorem}
{\em Proof:}
Considering the projected representation $\widetilde{\rho}$ into $PSL(2,K)$, 
we are in the situation of Hypothesis \ref{unboundedZD}. 
Apply all of the above sequence of results. We get a covering $Z\rightarrow X$ and a pluriharmonic map on $X$ which pulls
back to one on $Z$. 
Let 
$$
Z\rightarrow Y \rightarrow A
$$
be the Zariski-Stein factorization of the partial Albanese map $Z\rightarrow A$ discussed above.  
By Lemma \ref{dimone}, $Y$ is a curve. It is normal because
$Z$ is, so it is a smooth curve. The fibers of $Z\rightarrow Y$ are contained in the leaves of the foliation defined
by $\Re \alpha$, and the pluriharmonic map is constant on these leaves, by Corollary \ref{thetasplitting}. 
Thus, the pluriharmonic map is constant on the fibers of $Z\rightarrow Y$. 

The image of the fundamental group of the fiber is a normal subgroup $\Upsilon 
\subset \pi _1(Z,z)$. This acts on $\Tt$ with a nonempty set of fixed points, indeed
any image of a fiber by the pluriharmonic map is a point fixed by $\Upsilon$. So applying
Lemma \ref{normalsubgroup}, we get that $\Upsilon$ maps to the identity element of $PSL(2,K)$.
Therefore $\rho_Z$ projectively factors through the map $Z\rightarrow Y$. By 
\ref{projfactorsDMfactors}, $\rho _Z$ factors through a DM-curve. 
Then 
by Lemma \ref{descendfactorization} we get a factorization
over $X$. 
\eop

In the remainder of the paper, we will consider a representation $\rho : \pi _1(X,x)
\rightarrow SL(2,\cc )$ which doesn't factor
through a DM-curve. The above theorem (applied to various fields $K$) will allow us to prove that
$\rho$ is rigid and integral, and this will eventually give rise to a variation of Hodge structure
corresponding to a factorization through a Shimura modular stack.

\section{Rigidity}
\label{sec-rigidity}

Suppose $X$ is a smooth quasiprojective variety. Fix a smooth compactification $X\subset \overline{X}$
such that the complementary divisor $D$ has smooth irreducible components
$D_i$ meeting at normal crossing points. 
For each component $D_i$ of the divisor at infinity, let $\gamma _i$ denote a loop going out to a point near
$D_i$ then once around $D_i$ and back to $x$ along the same path. 

We start with some general notations. Suppose $G$ is a reductive group. Suppose we have chosen 
closed $ad(G)$-invariant subsets $C_i\subset G$ indexed by the components $D_i$ of the
divisor at infinity. The {\em representation space}
$$
R(X,x,G,\{ C_i\})
$$
is defined as the closed subset of $Hom (\pi _1(X,x), G)$ of representations $\rho $ such that $\rho (\gamma _i)\in C_i$.
The condition that $C_i$ be conjugation-invariant means that this condition is well-defined independently of
the choices of $\gamma _i$. 

We can then define the {\em moduli stack} and {\em coarse moduli space} denoted respectively
$$
{\mathcal M} (X,G,\{ C_i \} ) \;\; \mbox{and} \;\; M(X,G,\{ C_i \} ).
$$
These are the quotients of $R(X,x,G,\{ C_i\})$ by the conjugation action of $G$, first in the sense of algebraic
stacks and then in the sense of a universal categorical quotient of affine schemes. These are the ``Betti'' moduli
spaces \cite{LubotskyMagid}, there should also be ``de Rham'' and ``Dolbeault'' versions
\cite{BalajiBiswasNagaraj} \cite{Biswas} \cite{Biquard} \cite{BodenYokogawa} 
\cite{InabaIwasakiSaito} \cite{Konno} \cite{Nakajima} \cite{Nitsure} \cite{Moduli} \cite{SteerWren} 
\cite{Thaddeus} \cite{Yokogawa}. 

The main case we will consider is when $C_i$ are the closures of conjuacy classes in $G$. Suppose 
$\rho : \pi _1(X,x) \rightarrow G$ is a representation. Let $C_i(\rho )$ be the closure of the conjugacy class of
the image $\rho (\gamma _i)$. We say that a Zariski-dense representation 
$\rho$ is {\em rigid} if it represents an isolated point
in the moduli space $M(X,G,\{ C_i (\rho )\})$, and {\em nonrigid} otherwise. 

\subsection{Conjugacy classes in $SL(2)$}

Specialize now to the case $G=SL(2)$. In this case there are basically two types of conjugacy classes:
the unipotent matrices, and the semisimple matrices with distinct eigenvalues.
One has to add in the semisimple matrices with only one eigenvalue, which in this case means just $\pm 1$.
The notion of ``unipotent'' should also be extended to include $-1$ times a unipotent matrix which we sometimes call
``projectively unipotent''. 

\begin{lemma}
\label{conjclasses}
Suppose $C_i$ is the closure of a conjugacy class in $SL(2,K)$ with $K$ an algebraically closed field of characteristic
different from $2$. Then there are five cases:
\newline
(a)\, $C_i = \{ 1\}$; \newline
(b)\, $C_i = \{ -1\}$; \newline
(c)\, $C_i$ is the closure of the set of matrices conjugate to
$\left( \begin{array}{cc} 1 & 1 \\ 0 & 1 \end{array} \right)$; \newline
(d)\, $C_i$ is the closure of the set of matrices conjugate to
$\left( \begin{array}{cc} -1 & 1 \\ 0 & -1 \end{array} \right)$; or \newline
(e)\, there is $a\neq \pm 1$ in $K$ such that $C_i$ is the set of matrices
conjugate to 
$\left( \begin{array}{cc} a & 0 \\ 0 & a^{-1} \end{array} \right)$.
\end{lemma}
\eop

Note that the conjugacy classes 
(a) and (b) are contained in the closures of the conjugacy classes (c) and (d), 
whereas the conjugacy class in (e) 
is closed. The equations for cases (c), (d) and (e) are respectively 
$$
Tr (A)=2; \;\; Tr (A) = -2; \;\; Tr (A) = a + a^{-1}.
$$
The equations for (a) and (b) are not obtained by traces but rather directly by the equations $A=1$ or $A=-1$.
In any case one sees explicitly the equations to be imposed on $R(X,x,SL(2))$ to cut out
$R(X,x,SL(2), \{ C_i\} )$.

A closed subset $C_i$ is said to be {\em quasi-unipotent} if there is $n$ such that $A^n$ is unipotent for 
any $A\in C_i$. The closure of a conjugacy class is quasi-unipotent if and only if the conjugacy class is that of a
quasi-unipotent matrix. 
A conjugacy-class-closure which is quasi-unipotent 
is contained in one of the cases (a), (b), (c), (d), or (e) with $a$ a root of unity. 
Thus if the $C_i$ are quasi-unipotent then the equations defining $R(X,x,SL(2), \{ C_i \} )$
are defined over any ring extension of $\zz$ which contains the roots of unity in question. 

Special to the case of $SL(2)$ is the fact that a quasi-unipotent conjugacy class consists either
of matrices which are either of finite order (case (e)) or else projectively unipotent, i.e.
a central multiple of a unipotent matrix (cases (c) and (d)).

\subsection{Zariski density}

We need to keep track of the Zariski density hypothesis as we reduce and localize. 
First recall the connected subgroups of $SL(2,k)$.

\begin{lemma}
\label{sl2subs}
Suppose $k$ is \revision an algebraically closed field. 
The connected algebraic subgroups of $SL(2,k)$ are the following:
\newline
---the trivial subgroup;
\newline
---a torus, conjugate to the subgroup of diagonal matrices;
\newline
---a unipotent subgroup conjugate to the group of strictly upper triangular matrices;
\newline
---a solvable subgroup conjugate to the group of upper triangular matrices; or
\newline
---the whole group.

Furthermore, if $H\subset SL(2,k)$ is a strict algebraic subgroup which is not finite, then 
either $H/H^o$ is of order $\leq 2$, or else $H^o$ is conjugate to the unipotent subgroup, $H$ is contained in
the associated solvable subgroup, and 
$H/H^o$ is a finite cyclic subgroup of the torus. 
\end{lemma}
{\em Proof:}
The first part comes from the Levi decomposition and the fact that the only reductive groups
smaller than $SL(2)$ are tori. For the last statement, note that $H$ normalizes its connected component 
$H^o$. If $H^o$ is the diagonal torus, its normalizer is the semidirect product of $H^o$ with the group of
transpositions which has order two. If $H^o$ is the solvable (but not unipotent or toric) 
subgroup then it is equal to its own normalizer.
The remaining case is where $H^o$ is the unipotent subgroup; its normalizer is the 
solvable subgroup and $H/ H^o$ is contained in the torus (quotient of the solvable subgroup by the unipotent one).
\eop

\begin{corollary}
\revision
Still assuming $k=\overline{k}$, suppose $H\subset SL(2,k)$ is a strict algebraic subgroup containing an element whose eigenvalues are elements of 
infinite order in $k^{\ast}$. Then $H$ is either a torus, the solvable subgroup, or the normalizer of a torus.
\end{corollary}
{\em Proof:}
The condition on eigenvalues means that $H$ must contain a torus. Thus $H^o$ is either a torus or
the solvable subgroup. In the first case, $H$ is either the torus or its normalizer (which contains the
torus with index $2$). In the second case, $H=H^o$ since the solvable group is its own normalizer.
\eop

We will consider a subgroup $H$ containing an element whose eigenvalues are of infinite order.
We would like a good algebraic way of distinguishing between the cases $H=SL(2,k)$, and
the other three possibilities mentionned in the corollary.  For this, notice that for the connected subgroups
mentionned in Lemma \ref{sl2subs}, the traces of elements detect only the diagonal part and this is
commutative. Thus for example, $Tr (ABAB) = Tr (A^2 B^2)$ for these subgroups whereas this is not generically
true. The other possible subgroup, the normalizer of a torus, is
an extension of order two, which we can get around by looking at the squares of matrices. This leads to the
following Katz-style criterion \cite{KatzESDE}, also somewhat similar to Corollary 1.2.2 of \cite{CullerShalen}. 
We thank M. Larsen for comments about it.

\begin{lemma}
\label{zdcriterion}
Suppose $k$ is a field and  $\Gamma \subset SL(2,k)$ is a subgroup. Suppose that $\Gamma$ contains an
element whose eigenvalues are of infinite order in $k^{\ast}$. Then $\Gamma $ is Zariski-dense if and
only if there exist elements $\alpha , \beta \in \Gamma$ such that 
$$
Tr (\alpha ^2 \beta ^2 \alpha ^2 \beta ^2) - Tr (\alpha ^4 \beta ^4) \neq 0.
$$
\end{lemma}
{\em Proof:}
\revision
We reduce to the case $k=\overline{k}$ so that there is no problem with twisted forms of tori. 
If $\Gamma \subset SL(2,k)$ is not Zariski-dense in $SL(2,\overline{k})$
then the Zariski closure of $\Gamma$ over
$\overline{k}$ is a subgroup $G\subset SL(2,\overline{k})$ preserved by $Gal (\overline{k}/k)$,
in particular it comes from a subgroup over $k$. Thus $\Gamma$ is not Zariski-dense over $k$.
The other direction is easy, so Zariski-density is preserved by going from $k$ to $\overline{k}$.
The criterion given in the statement of the lemma doesn't depend on $k$ so it is preserved too.
Thus we may assume $k=\overline{k}$. 

Now we get to the proof. If $H\subset SL(2,k)$ is a torus, the normalizer of a torus, or
a solvable subgroup, and if $\alpha , \beta \in H$ then 
$$
Tr (\alpha ^2 \beta ^2 \alpha ^2 \beta ^2) - Tr (\alpha ^4 \beta ^4) = 0.
$$
In view of Lemma \ref{sl2subs}, this shows that if $\alpha , \beta \in \Gamma$ exist as in the statement,
then $\Gamma$ is Zariski-dense. 

Suppose now that $\Gamma$ is Zariski-dense. It is easy to see that $\Gamma \times \Gamma$ is 
Zariski-dense in $SL(2,k)\times SL(2,k)$. On the other hand, the set of pairs of matrices
$(A,B)\in  SL(2,k)\times SL(2,k)$ where the identity $Tr (A^2B^2A^2B^2) = Tr (A^4B^4)$
holds is a Zariski-closed subset.
In order to prove that $\Gamma \times \Gamma$ goes outside of this subset it suffices to exhibit
any pair of matrices $(A,B)\in  SL(2)\times SL(2)$ which is not in the subset. Furthermore we can do this
after going to the algebraic closure. 

The subset of matrices in $SL(2,\overline{k})$ which are squares is Zariski-dense, so again its product is Zariski dense
in the product. Thus, it suffices to exhibit a pair of matrices 
$(A,B)\in  SL(2,k)\times SL(2,k)$ such that $Tr (ABAB) \neq Tr (A^2B^2)$. For this, set
$$
A:= \left( \begin{array}{cc} 1 & 1 \\ 0 & 1 \end{array} \right) 
$$
$$
B:= \left( \begin{array}{cc} 1 & 0 \\ 1 & 1 \end{array} \right) 
$$
with $Tr(ABAB) = 7 \neq 6 = Tr(A^2B^2)$ in any characteristic.
\eop

\begin{lemma}
\label{rigidequiv}
Suppose $\rho$ is a Zariski-dense representation. Then $\rho$ is rigid in the above sense, i.e. 
it represents an isolated point
in the moduli space $M(X,SL(2), \{ \overline{C}_i(\rho ) \} )$ ,
if and only if  there is no
non-isotrivial family of representations all having the same conjugacy classes at infinity, going through
$\rho$. 
\end{lemma}
{\em Proof:}
We have a map 
$$
R(X,SL(2), \{ \overline{C}_i(\rho ) \} )\rightarrow M(X,SL(2), \{ \overline{C}_i(\rho ) \} ).
$$
As is well-known, the points of the moduli space $M$ represent $S$-equivalence classes of representations in 
$R$. In particular, if $\rho$ is an irreducible representation as is the case here, any representation
$S$-equivalent to $\rho$ is actually isomorphic to $\rho$. Thus the fiber 
of the map $R\rightarrow M$ over the point $[\rho ]$ consists only of conjugates of $\rho$, in other words
the fiber is set-theoretically an orbit of $SL(2)$. The point $[\rho ]\in M$ is an isolated point if
and only if the fiber $SL(2)\cdot \rho$ is a connected component of $R$. 

If this is the case, and if $\{ \rho _t\}$ is a continuous family of representations with $\rho _0=\rho$
and all having the same conjugacy classes at infinity, then we get a path in 
$R(X,SL(2), \{ \overline{C}_i(\rho ) \} )$ going through $\rho$. Since the orbit 
$SL(2)\cdot \rho$ is a connected component, this means that all of the $\rho _t$ are contained in the orbit
so the family is isotrivial. This shows one direction of the lemma.

In the other direction, suppose that $[\rho ]$ is not an isolated point of $M$. Then the orbit 
$SL(2)\cdot \rho$ is not a connected component of $R$. However, it is a closed subset. Therefore there is a path
$\{ \rho _t\}$ in $R$ which starts at $\rho _0$ in the orbit, and goes out of the orbit. Up to changing
the basis of the underlying space we may assume $\rho _0 = \rho$. The family is nonisotrivial. We just have
to see that the conjugacy classes at infinity stay the same. The monodromy transformations at infinity
$\rho _t(\gamma _i )$ are families in $\overline{C}_i(\rho )$. However, by definition 
$\overline{C}_i(\rho )$ is the closure of the conjugacy class $C_i(\rho )$
of $\rho (\gamma _i)$, in particular $C_i(\rho )$ (which is a locally closed subset of the group) is 
an open subset of $\overline{C}_i(\rho )$. On the other hand, $\rho _0(\gamma _i)=\rho (\gamma _i)$ is
in this open subset. So, for $t$ close to $0$, we still have $\rho _t(\gamma _i)\in C_i(\rho )$
which is the condition we needed. 
\eop

\begin{lemma}
\label{rigidinvariant}
If $V$ is a rigid local system and $\sigma : \cc \rightarrow \cc$ is an automorphism of $\cc$
then $V^{\sigma}$ is also rigid. 
\end{lemma}
{\em Proof:}
Standard.
\eop

\subsection{Nonfactorization implies rigidity}

To motivate the proof of the next main theorem, we note a corollary of Lemma \ref{finiteDMfactors}.

\begin{corollary}
\label{nonfactorsopen}
The subset of points in the space of $SL(2,\cc )$ representations of $\pi _1(X,x)$
corresponding to representations which don't
factor through a hyperbolic DM curve, is open. 
\end{corollary}
{\em Proof:}
By Lemma \ref{finiteDMfactors} there are only a finite number of maps from $X$ to hyperbolic DM-curves.
For each of these maps, the set of representations which factor is given by a closed subset.
\eop

\begin{theorem}
\label{nonfactorsrigid}
Suppose $\rho : \pi _1(X,x)\rightarrow SL(2,\cc )$ is a Zariski-dense
representation with quasi-unipotent monodromy at infinity. If $\rho$  doesn't factor through a map 
from $X$ to a DM-curve, then it is rigid in the
above sense.
\end{theorem}
{\em Proof:}
Suppose that we have a representation $\rho$ which doesn't factor through a DM-curve, and which is not rigid.
We derive a contradiction. Note that $\rho$ doesn't projectively factor through an orbicurve 
(Lemma \ref{projfactorsDMfactors}).

We first give the short if somewhat heuristic version of the argument. From 
Corollary \ref{nonfactorsopen} there is an open subscheme of the moduli scheme corresponding to
representations which don't factor through DM-curves. We can take an open subset to obtain an
irreducible positive-dimensional open subscheme whose representations don't factor. This open set
may be considered as having finite type over $\zz$. Then there is a
curve defined over a finite field mapping into here. The trace functions of the monodromy elements are not
all constant, in particular there is a point on the curve where at least one trace function has a pole.
Taking the completion at this point we obtain a representation 
$$
\hat{\rho} : \pi _1(X,x)\rightarrow SL(2,\ff _q((t)))
$$
whose image is not contained in a compact subgroup. Using the criterion of Lemma \ref{zdcriterion} we can
assume that $\hat{\rho}$ is Zariski dense.  This is an elementary example of the more general problem 
of independence of $\ell$ for monodromy groups such as considered in \cite{Larsen}.  
By Theorem \ref{mainfactorization} the representation $\hat{\rho}$ would have to factor through a DM curve,
contradicting our choice of open set.

To get a more accurate version of this proof, we fill in explicitly some of the equations for the 
properties preserved when going to the completion at a finite point. The argument uses what is known as 
a ``spread'' \cite{GreenGriffiths}. 

Start by enumerating the possible factorization maps through hyperbolic orbicurves. Denote these
by $f_i : X \rightarrow Y_i$ for a finite set of indices $i\in I$. Let $\varphi _i : \pi _1(X,x) \rightarrow \Gamma _i$
be the corresponding maps on fundamental groups, and let $K_i$ be the kernel of $\varphi _i$. 
The fact that $\rho$ doesn't projectively factor through $f_i$ translates by saying that $\rho (K_i)$ is not
contained in the center of $SL(2,\cc )$. 

Note that $K_i$ is a normal subgroup of $\pi _1(X,x)$ so the Zariski closure of $\rho (K_i)$ is a normal subgroup
of the full image of $\rho$ which is by assumption $SL(2,\cc )$. The only normal subgroups are the identity, the
center, and the group itself. Therefore $\rho |_ {K_i}$ is Zariski-dense in $SL(2,\cc )$. In particular there
exists an element $k_i \in K_i$ such that 
$$
Tr (\rho (k_i)) \neq \pm 2.
$$
The nonfactorization hypothesis translates into the existence of such a $k_i$ for each $i\in I$. 

The hypothesis that $\rho $ is nonrigid means that there is a curve of representations passing through
$\rho$ and projecting to a nontrivial curve in the moduli space of representations. Explicitly write down this
curve as follows. We may assume it is affine, and even (applying Noether normalization after going to an 
open set if necessary) that it is finite over
${\mathbb A}^1$ with a coordinate ring of the form
$$
A = \frac{ \cc [y] }{(y^m + a_{m-1}y^{m-1} + \ldots + a_0)}.
$$
Our path of representations is a representation
$$
\rho _A : \pi _1(X,x)\rightarrow SL(2,A)
$$
and there is a point $p:A\rightarrow \cc$ such that $\rho _A \otimes _{A,p}\cc = \rho$.

Recall that Procesi's theorem says that the traces of the images of group elements give an embedding of the
moduli space  \cite{Procesi} \cite{CullerShalen}. Therefore the hypothesis that our curve is nontrivial 
in the moduli space means that there is an element $\gamma \in \pi _1(X,x)$ such that $Tr (\rho _A(\gamma ))$
is not contained in the ground field $\cc$. Concretely, the elements of $A$ are written as
sums $u_{m-1}y^{m-1} + \ldots + u_0$ so the non-constancy of this element means that
one of its coordinates $u_j$ is nonzero for $0<j\leq m-1$. Write this condition as
$u_j (Tr (\rho _A(\gamma ))) \neq 0$.

Similarly, the nonfactorization hypothesis says
$Tr (\rho _A(k_i))\neq \pm 2$ for all $i\in I$. 

We can specialize our curve of representations so that it is defined over a finite extension field $L$ of $\qq$
(in other words we specialize the coefficients $a_i$ of the equations for $A$ as well as the coefficients
of the matrix coefficients of $\rho _A$). We can do this and still keep the nonfactorization conditions that
$Tr (\rho _A(k_i))\neq \pm 2$ and the nonrigidity condition that 
$u_j (Tr (\rho _A(\gamma ))) \neq 0$.

Now our curve of representations can also be defined over a finite extension ring $R$ of $\zz$
(in general $R$ will be the localization of the ring of integers $\Oo _L$ at a finite number of 
primes). Then we can find a prime ideal ${\mathfrak p}\subset R$ 
such that the nonfactorization and nonrigidity conditions
hold modulo ${\mathfrak p}$. We obtain a curve of representations 
$\rho _B : \pi _1(X,x) \rightarrow SL(2,B)$ 
such that 
$$
B = \frac{ {\mathbb F}_q [y] }{(y^m + b_{m-1}y^{m-1} + \ldots + b_0)},
$$
such that $Tr (\rho _B(k_i))\neq \pm 2$ and $u_j (Tr (\rho _B(\gamma ))) \neq 0$.
We may assume that the characteristic of the finite field ${\mathbb F}_q$ is different from $2$.

Note that $Spec(B)$ is a curve over the finite field
${\mathbb F}_q$ and we may view the element $Tr (\rho _B(\gamma ))$ as a nonconstant 
function on this curve.
Complete at a point at infinity on this curve, choosing a point with the property
that the function $Tr (\rho _B(\gamma ))$ has a pole.  The completed field, which may be seen as the completion
of $B$ with respect to the nonarchimedean norm given by the order of pole at the point, will be denoted by 
$\widehat{B}$.
It is a complete local field of the type envisioned in Hypothesis \ref{unboundedZD} above. 
The ring of integers $\Oo _{\widehat{B}}$
is a discrete valuation ring with finite residue field $\ff _q$, so 
$$
\Oo _{\widehat{B}} \cong \ff _q [[t]], \;\;\; 
\widehat{B}\cong \ff _q ((t)).
$$
We have a representation
$$
\rho _{\widehat{B}}: \pi _1(X,x)\rightarrow SL(2,\widehat{B} ) \cong SL(2,\ff _q ((t))\, ).
$$
The fact that $Tr (\rho _{\widehat{B}}(k_i))\neq \pm 2$ means that $\rho _{\widehat{B}}$ doesn't projectively
factor
through the orbicurve $f_i:X\rightarrow Y_i$. This holds for all $i\in I$, so 
$\rho _{\widehat{B}}$ doesn't factor through an orbicurve. On the other hand, the fact that
$Tr (\rho _B(\gamma ))$ has a pole at the point where we took the completion means that 
$Tr (\rho _{\widehat{B}}(\gamma ))$ is not in $\Oo _{\widehat{B}}$. 

We have $\rho _{\widehat{B}}(\gamma )\in SL(2,\widehat{B})$ a matrix with eigenvalues 
$\beta ^{\pm 1}$ for $\beta \in \overline{\widehat{B}}$, therefore
$\rho _{\widehat{B}}(\gamma ^n)$ has eigenvalues $\beta ^{\pm n}$.
In particular if $Tr(\rho _{\widehat{B}}(\gamma ))$ has a pole of order 
$\geq 1$ then $Tr (\rho _{\widehat{B}}(\gamma ^n))$ has a pole of order $\geq n$. 
This shows that 
$\rho _{\widehat{B}}$ doesn't have image in any compact subgroup of $SL(2,\widehat{B})$.

This situation leads to a contradiction with Theorem \ref{mainfactorization}.  Before getting there,
we still have to verify the hypothesis
that $\rho _{\widehat{B}}$ have Zariski dense image. This is done much as above, using
the explicit criterion \ref{zdcriterion}. 
Note that the image of $\rho$ was not finite, so it contained an element of infinite order.
By Lemma \ref{zdcriterion}, there exist $\alpha , \beta \in \pi _1(X,x)$ such that
$$
Tr (\rho (\alpha ^2 \beta ^2 \alpha ^2 \beta ^2 )) - Tr (\rho (\alpha ^4 \beta ^4)) \neq 0.
$$
In our process of reducing modulo a prime, we may choose a prime which doesn't divide this nonzero
quantity. Then generalizing along a curve and completing we still obtain a representation 
$\rho _{\widehat{B}}$ with the property that 
$$
Tr (\rho _{\widehat{B}}(\alpha ^2 \beta ^2 \alpha ^2 \beta ^2 )) - 
Tr (\rho _{\widehat{B}}(\alpha ^4 \beta ^4)) \neq 0.
$$
On the other hand, our construction of $\rho _{\widehat{B}}$ was such that its image contains an element
of infinite order. Thus, Lemma \ref{zdcriterion} again applies (choose our prime to have characteristic 
different from $2$). We conclude that the image of $\rho _{\widehat{B}}$ is Zariski dense. 

We can now apply Theorem \ref{mainfactorization} to obtain a contradiction, completing the proof of the theorem.
\eop

{\em Remark}: 
The above proof constitutes a new proof of the factorization result of \cite{ubiquity}
in the projective case. Delzant seems to have found a similar proof, as he has indicated to us in
correspondence. Many of the details---fundamentally elementary---going into the proof, 
stem from wanting to refer only to trees of
finite type. The theory of Gromov and Schoen should in principle allow us to prove the result
more directly for a complex curve of representations going to infinity, but this would involve considerations
of non-locally compact trees which we have preferred to avoid. 

Alternatively, one could adapt to the 
quasiprojective (and quasi-unipotent monodromy) case the archimedean argument of \cite{ubiquity}.
This type of thing is currently being done by T. Mochizuki in a much wider context \cite{TMochizuki} \cite{TMochizuki2}. 
For example he
obtains the theorem that any representation can be deformed to a complex variation of Hodge structure. 
We felt that in the context of the present paper, it would be more interesting to deduce Theorem 
\ref{nonfactorsrigid}  from the Gromov-Schoen theory of harmonic
maps to trees, since this provides a unified approach to the two basic results
(rigidity as above, and integrality in the next section). We cannot entirely avoid refering to the
archimedean case: in \S \ref{sec-cvhs} below we state the theorem that rigid representations are
complex variations of Hodge structure, which requires the theory of harmonic maps to symmetric spaces.

\section{Integrality}
\label{sec-integrality}

Suppose $\rho : \pi _1(X,x) \rightarrow SL(2,\cc )$ is a representation.  
Following Bass \cite{Bass} we say that $\rho $ is {\em integral} if for every $\gamma \in \pi _1(X,x)$ the trace
$Tr (\rho (\gamma ))$ is an algebraic integer. Bass proved the following fundamental characterization.

\begin{lemma}
A Zariski-dense representation $\rho$ has integral traces if and only if
there exists a number field $L$ (i.e. a finite extension of $\qq$) with its ring of algebraic
integers $\Oo _L\subset L$, an embedding $\eta : L\rightarrow \cc$, 
a rank two projective $\Oo _L$-module $V_{\Oo _L}$, 
and a representation 
$$
\rho _{\Oo _L} : \pi _1(X,x) \rightarrow SL(V_{\Oo _L})
$$
such that
$$
\rho _{\Oo _L} \otimes _{\Oo _L,\eta }\cc \cong \rho .
$$
\end{lemma}
{\em Proof:} See \cite{Bass}.
\eop

\begin{lemma}
\label{integralfiniteindex}
Suppose $V$ is a local system, and suppose $Y\rightarrow X$ is a finite etale covering such that
$V|_Y$ is integral. Then $V$ is integral.
\end{lemma}
{\em Proof:}
A local system is integral if and only if for any element $\gamma \in \pi _1(X,x)$, all of the eigenvalues
of $\rho _V(\gamma )$ are algebraic integers. Under the hypotheses of the lemma, if $\gamma \in \pi _1(X,x)$
then there is an $n$ such that $\gamma ^n$ comes from $\pi _1(Y,y)$ (where $y$ is a lift of the basepoint
to $Y$). Then the eigenvalues of  $\rho _V(\gamma ^n)$ (which are the $n$-th powers of 
the eigenvalues of $\rho _V(\gamma )$) are algebraic integers. This imples that the eigenvalues
of $\rho _V(\gamma )$ are algebraic integers (indeed an equation of the form $z^n = \alpha$ 
is an integral equation for $z$, if $\alpha$ is an algebraic integer). 
\eop

The main result about integrality is the following. It is pretty well-known, see the work of Jost and Zuo
\cite{JostZuo2}, Eyssidieux \cite{Eyssidieux}, Katzarkov \cite{Katzarkov}, Klingler \cite{Klingler},
Napier-Ramachandran \cite{NapierRamachandran3}  and others
which are based on the theory of Gromov-Schoen \cite{GromovSchoen}. 

\begin{theorem}
\label{nonfactorsintegral}
If $\rho : \pi _1(X,x)\rightarrow SL(2,\cc )$ is 
a Zariski dense representation with quasi-unipotent monodromy at infinity, and $\rho$ doesn't
factor through a DM-curve, then $\rho$ is integral. 
\end{theorem}
{\em Proof:}
Suppose $\rho$ doesn't factor. By Theorem \ref{nonfactorsrigid} $\rho$ is rigid. This implies that it can
be defined over an algebraic number field $L$ (which we may assume is finite over $\qq$).
Assume that it is not integral. Then the trace of some monodromy element is not an
algebraic integer, and in particular it is not integral at some prime of $L$.
Taking the completion at this prime, we get a Zariski dense representation into 
$SL(2,K)$ where $K$ is a complete local field, and the representation doesn't go into a 
compact subgroup because the trace of some element is not an integer of $K$.
The residue field of $K$ is finite because the integers of $L$ are of finite type over $\zz$.
The monodromy transformations at infinity are the same as for the original $\rho$, so they are quasi-unipotent.  
Therefore Theorem \ref{mainfactorization} applies: the representation into $SL(2,K)$ factors through a DM-curve.
By invariance of factorization under change of base field, Corollary \ref{invariancefieldextension}, 
the original representation
factors through a DM-curve contradicting our hypothesis. Therefore $\rho$ is integral.
\eop

\subsection{Hypergeometric cases}
\label{sec-hypergeometric}

There is not a perfect dichotomy between the properties of factorization through DM-curves and rigidity. 
There are some cases where these overlap, corresponding to rigid local systems on DM-curves.
For Zariski-dense representations on orbicurves,
rigidity can be determined by a dimension count which depends only on the local monodromy data. 
Katz's algorithm \cite{Katz} in principle gives a way to determine which rigid monodromy data can come from
representations. Furthermore, when the data do come from representations, his algorithm provides an explicit construction
of the representation as a motivic one. As a corollary, rigid representations on DM-curves are automatically
integral. 

In the case of local systems of rank two, the only rigid cases correspond to hypergeometric equations; their
motivic expression is basically the classical integral expression of hypergeometric functions \cite{DeligneMostow} \cite{KatzESDE} \cite{Katz}.
In particular, the result that rigidity implies integrality can be understood explicitly in this case.

Suppose $Y$ is a smooth compact curve with points $y_1,\ldots , y_k$ and basepoint $y$ distinct
from the marked points. Suppose $\rho : \pi _1(Y-\{ y_i\}, y) \rightarrow SL(r,\cc )$ is a
Zariksi-dense representation, which is rigid among those with fixed conjugacy classes at the punctures.

\begin{lemma}
\label{genuszero}
In the above situation, $Y\cong \pp ^1$.
\end{lemma}
{\em Proof:}
A dimension count shows that
the only potentially problematic case which needs to be discussed is that of representations on an elliptic
curve, with prescribed monodromy around a single puncture. This yields the equation 
$$
aba^{-1} b^{-1} \in C
$$
for a fixed conjugacy class $C$. The space of solutions is divided by the conjugation action of $SL(2)$. 
The dimension of the space of solutions, divided by the conjugation action, is therefore at least 
equal to 
$$
2{\rm dim} (SL(2)) - {\rm dim} (SL(2)) - ({\rm codim} _{SL(2)}C) = dim (C).
$$
In particular, as soon as $C$ is positive-dimensional, the representation cannot be rigid. The only zero-dimensional
conjugacy classes are $\{ 1 \}$ and $\{ -1\}$. In the first case the image is abelian; and in the second case
it is an extension of an abelian group by $\{\pm 1\}$. Neither of these can be Zariski-dense subgroups of
$SL(2)$. Therefore this case cannot arise, so $Y$ has to be rational.
\eop

Now we proceed with a dimension count. Note that the naive dimension count actually gives the dimension of the
moduli space at a Zariski-dense point, in this case \cite{Katz}. Let $C_1,\ldots , C_k$ denote the conjugacy 
classes of the representation at the points $y_i$. Here we are including orbifold points because the space of
representations on the orbifold is the same as the space of representations on the Zariski open subset with prescribed
conjugacy classes (in the case of orbifold points, the conjugacy classes will be of finite order).

Recall that the space of representations is the space of solutions of $a_1\cdots a_k = 1$ with $a_i\in C_i$, 
the whole up to a global conjugation (and on Zariski-dense points the conjugation action is faithful).
The virtual dimension of the moduli space is therefore 
$$
v := {\rm dim} (C_1) + \ldots + {\rm dim} (C_k) - 2{\rm dim} (SL(2)).
$$
Rigidity happens for Zariski-dense representations exactly when $v=0$. 

Recall that the dimensions of the conjugacy classes for the cases given by Lemma \ref{conjclasses} are as follows:
for (a) and (b), ${\rm dim} (C)=0$; whereas for (c), (d) and (e), ${\rm dim} (C)=2$. Together with the above we conclude
the following statement. 

\begin{lemma}
If $\rho$ is a rigid rank two representation on a DM-curve, then there are exactly three orbifold or singular points
with conjugacy classes different from $\{ 1\}$ or $\{ -1\}$.
\end{lemma}
\eop

The points with conjugacy classes $\{ -1\}$ can be chosen freely; they can also combine with the three main points
(resulting in an interchange between (c) and (d), or taking the negative of the eigenvalue  in case (e)). 
These operations correspond to tensoring with rank one local systems on $Y$ (in fact, the Prym local systems of
hypergeometric curves over $Y$). From now on we ignore these shifts, and use the term ``unipotent'' for cases (c) and (d). 

If all three conjugacy classes are of type (e), then we are in the case of the classical hypergeometric equation.
It is well  known that for eigenvalues which are roots of unity, these hypergeometric local systems exist, and are
motivic (hence integral). 

We give explicit calculations for the cases where one or more of the conjugacy classes is unipotent, showing that
the monodromy representations are integral. Note by \cite{Katz} that if a rigid representation exists,
it is uniquely determined by the monodromy data; thus it suffices to write down a single integral representation for
each given collection of monodromy data. In a similar way, if there exists a reducible representation with given 
conjugacy classes, then there will not exist an irreducible one. In this case we must multiply one of the matrices by $-1$.
Also for these calculations, notice that a matrix with trace $\pm 2$ which is not $\pm$ the identity, must be
in conjugacy classes (c) or (d) (i.e. nontrivially unipotent). 

For three unipotent matrices, one or three must have eigenvalues $-1$ otherwise an irreducible solution will not 
exist. In this case we have the solution 
$$
\left( \begin{array}{cc} 1 & 1 \\ 0 & 1 \end{array} \right) 
\left( \begin{array}{cc} 1 & 0 \\ -4 & 1 \end{array} \right)
= 
\left( \begin{array}{cc} -3 & 1 \\ -4 & 1 \end{array} \right) .
$$
For  two unipotent matrices and one matrix of type (e) with eigenvalue $\alpha$, note that 
$$
\left( \begin{array}{cc} 1 & 1 \\ 0 & 1 \end{array} \right) 
\left( \begin{array}{cc} 1 & 0 \\ b & 1 \end{array} \right)
= 
\left( \begin{array}{cc} 1+b & 1 \\ b & 1 \end{array} \right) .
$$
If $b= \alpha + \alpha ^{-1} -2$ then the trace on the right will be $\alpha + \alpha ^{-1}$ so we get a solution.
For any root of unity $\alpha \neq \pm 1$ the coefficient $b$ is an algebraic integer so this gives an integral solution.
Note also that there is a solution for $-\alpha$, so by tensoring with a rank one system we can treat 
unipotent conjugacy classes
of type (d) also. 

Finally, for the case of one unipotent matrix and two matrices of type (e) with eigenvalues $\alpha$ and $\beta$,
look at
$$
\left( \begin{array}{cc} \alpha & 1 \\ 0 & \alpha ^{-1} \end{array} \right) 
\left( \begin{array}{cc} \beta  & 0 \\ x & \beta ^{-1} \end{array} \right)
= 
\left( \begin{array}{cc} \alpha \beta + x & \beta ^{-1} \\ x\alpha ^{-1} & (\alpha \beta )^{-1} \end{array} \right) .
$$
Up to tensoring by a rank one Prym system (corresponding to changing the sign of one of $\alpha$ or $\beta$)
we can assume that we want the trace on the right to be equal to $2$. Notice that in any case the matrix on the
right is not the identity so if its trace is $2$ it will be of type (c). This gives the equation
$$
x = 2 - (\alpha \beta  + (\alpha \beta )^{-1}).
$$
If $\alpha$ and $\beta$ are roots of unity and $\alpha \beta \neq 1$ then this has an integral solution.
If $\alpha \beta = 1$ then there is a non-irreducible representation so there can be no irreducible solution
(however if one of the matrices is multiplied by $-1$ then there is a solution). 

We have concluded the proof of the following result.

\begin{lemma}
\label{DMrigidintegral}
If $\rho$ is a Zariski dense rank two rigid representation on a DM-curve, then it is integral. 
\end{lemma}
\eop

{\em Remark:} Our main result \ref{rigidcase} below implies that in all of the above cases, the representation
is motivic. Thus we are spared the trouble of writing down an explicit motivic presentation for the cases
of one or more unipotent conjugacy classes. Such presentations certainly exist in the literature, 
see  \cite{DettweilerReiter} \cite{KatzESDE} \cite{Katz}. Darmon discusses motivic properties of rigid local systems on $\pp ^1$ in 
an arithmetic setting \cite{Darmon}.

\section{Variations of Hodge structure} 
\label{sec-cvhs}
We recall here the classical definitions concerning complex, real and integral variations of Hodge structure.
In this section, let $V$ be a complex local system over $X$. We think of it as a $C^{\infty}$ vector
bundle with flat connexion $\nabla$. A structure of {\em $\cc$VHS of weight $w$} on $V$ is the data
of a decomposition of $C^{\infty}$ vector bundles $V =\bigoplus _{p+q=w}V^{p,q}$ satisfying the 
{\em Griffiths identities} 
$$
\nabla (V^{p,q}) \subset A^{1,0}(X; V^{p,q}) \oplus 
A^{1,0}(X; V^{p-1,q+1}) \oplus 
A^{0,1}(X; V^{p,q}) \oplus 
A^{0,1}(X; V^{p+1,q-1}).
$$
A {\em real variation of Hodge structure} arises when $V$ has a real structure $V_{\mathbb R}$ (i.e. a real local system
whose tensor product with ${\mathbb C}$ is $V$), such that $V^{p,q}=\overline{V^{q,p}}$.
Finally an {\em integral variation of  Hodge structure} is a real variation of Hodge structure together with an 
integral local system $V_{\mathbb Z}\subset V_{\mathbb R}$ whose extension by scalars to ${\mathbb R}$ is $V_{\mathbb R}$.

If $V$ is a $\cc$VHS of weight $w$, an {\em hermitian polarization} of $V$
is a sesquilinear $\nabla$-invariant form $\Phi ( \cdot ,\cdot )$
on $V$, such that $\Phi$ is $(-1)^1$-sesquisymmetric
$$
\Phi (u,v) =(-1)^w \overline{\Phi (v,u)};
$$
the different Hodge subspaces are orthogonal
$$
\Phi (u,v)=0\;\;\; \mbox{for}\;\;\; u\in V^{p,q}, \;\; v\in V^{r,s}, \;\; (p,q)\neq (r,s);
$$
and such that the positivity condition 
$$
(\sqrt{-1})^{p-q}\Phi (u,u) > 0 ,\;\;\; u\in V^{p,q}
$$
holds. 

A {\em real polarization} on a real variation $(V,V_{\rr})$ is an hermitian polarization $\Phi$ such that
$\Phi (u,v)\in \rr$ whenever $u,v\in V_{\rr}$.  
An {\em integral polarization} on 
an integral variation $(V,V_{\zz})$ is an hermitian polarization $\Phi$ such that
$\Phi (u,v)\in \zz$ whenever $u,v\in V_{\zz}$.  

In the weight $1$ case which interests us here, our conventions for the polarization
say that $\Phi$ should be a hermitian antisymmetric form, in particular $\sqrt{-1}\Phi$ is an indefinite unitary 
form with $\sqrt{-1}\Phi (u,u)$ positive
for $u\in V^{1,0}$ and negative for $u\in V^{0,1}$. If $\Phi$ is a real or integral polarization,
then it induces an alternating form on $V_{\rr}$ or $V_{\zz}$ in particular $\Phi (u,u)=0$ for $u\in V_{\rr}$.
The hermitian antisymmetric form $\Phi$ on the complexification is 
uniquely determined by $\Phi |_{V_{\rr}\otimes V_{\rr}}$
and the positivity condition holds on the Hodge subspaces which don't intersect $V_{\rr}$.

Recall \cite{DeligneTravauxDeShimura} \cite{DeligneCorvallis}
that, if we make the convention that $V^{p,q}$ is to be nonzero only when $p,q\geq 0$
then an integral polarizable weight one variation of Hodge structure is the same thing as an algebraic smooth family of
abelian varieties over $X$. 

\begin{theorem}
\label{rigidvhs}
Suppose $V$ is a rank two local system with quasi-unipotent monodromy at infinity. If $V$  is rigid  as
defined in \S \ref{sec-rigidity}
then it underlies a complex variation of Hodge structure of weight one. 
\end{theorem}
{\em Proof:}
This is Corollary 1 of \cite{actesToulouse}, see also \cite{JostZuo1}, \cite{TMochizuki2}.
Since we are talking about rank two local systems, there are two possibilities for the Hodge types:
either the local system is unitary in which case one can assume that there is a single Hodge type,
or else there are two adjacent Hodge types with Hodge  numbers $h^{1,0} = h^{0,1} = 1$. 
\eop

\begin{corollary}
\label{nonfactorsvhs}
Suppose $V$ is a rank two local system with quasi-unipotent monodromy at infinity, whose monodromy representation
is Zariski dense. If $V$ doesn't factor 
through a DM-curve, then for any automorphism $\sigma$ of $\cc$, the conjugate local system $V^{\sigma}$ underlies
a complex variation of Hodge structure of weight one. 
\end{corollary}
{\em Proof:}
By Theorem \ref{nonfactorsrigid}, $V$ is rigid. 
Therefore any conjugate is also rigid, and by Theorem \ref{rigidvhs} the conjugates
all underly complex variations of Hodge structure.
\eop

When a representation is integral and all of the Galois conjugates come from complex variations of
Hodge structure, then adding them together we get a $\zz$-variation of Hodge structure. 

\begin{lemma}
\label{integralvhsZvhs}
If $V_{\Oo _L}$ is a local system of projective $\Oo _L$-modules for an algebraic number field $L$,
such that for every embedding $\eta : L\rightarrow \cc$ the resulting complex local system
$V^{\eta}$ underlies a complex variation of Hodge structure, then $V_{\Oo _L}^{\oplus 2}$ considered as a local system 
of $\zz$-modules underlies a $\zz$-variation of Hodge structure. If each factor can be associated to a 
variation of weight one, then the $\zz$-variation of Hodge structure can be assumed to have weight one and
in particular it comes from a family of abelian varieties.
\end{lemma}
{\em Proof:} See \cite{hbls}.
\eop

As a corollary we get a first statement which approaches our classification statement:

\begin{corollary}
\label{classif-first}
Suppose $\rho : \pi _1(X,x)\rightarrow SL(2,\cc )$ is a representation with quasi-unipotent monodromy
at infinity, which doesn't factor through
a map to a DM-curve. Then there is a family of abelian varieties
over $X$ such that $\rho$ is
a direct factor of the underlying complex monodromy representation.
\end{corollary}
{\em Proof:}
By Corollary \ref{nonfactorsvhs} we are in the situation of Lemma \ref{integralvhsZvhs},
so $\rho$ is a complex direct factor in the monodromy of a weight one $\zz$-variation of Hodge structure;
this is the monodromy representation of the corresponding family of abelian varieties. 
\eop

The task in the remainder of the paper will be to analyse more carefully the kind of variation of Hodge structure
which can be made to occur here, and show that it corresponds to a factorization of the representation through a
Shimura modular stack.

\section{Polydisk Shimura modular stacks} 
\label{sec-shimura}

Moduli for families of abelian varieties is obtained through the theory of {\em Shimura varieties}. 
In general, a Shimura variety will have a tautological representation of $\pi _1$ into some algebraic group.
Since we are interested in rank $2$ representations, we consider the special case where this algebraic group is some form of 
$SL(2)$. The classical case which has been the most extensively studied is that of {\em Hilbert modular varieties},
corresponding to $SL(2,F)$ where $F$ is a totally real field \cite{Darmon}
\cite{Freitag} \cite{Goren}  \cite{HirzebruchVandeVen} \cite{McMullen}  \cite{Rapoport} \cite{vanderGeer}. 

It turns out that the Hilbert modular case is very slightly overly restrictive
for our purposes. For a rigid and integral representation thanks to a lemma explained to us by M. Larsen, we can assume that the
representation goes into $SL(2,L)$ where $L$ is a totally imaginary quadratic extension of a totally real field $F$, see
\S \ref{sec-improvements} below. 
However, the universal example which we consider in this section shows that in general, the traces will not necessarily lie in $F$,
so it is not always possible to reduce to $SL(2,F)$. We thus fall into the more general realm of the theory of
Shimura varieties, more specifically unitary Shimura varieties \cite{Larsen2}
\cite{Pappas} \cite{RapoportZink} which have recently come to be known as ``Shimura varieties of PEL type'' \cite{Wedhorn}. 
Our case is midway between the Hilbert modular case and the general unitary case, since the fact that the group is $SL(2)$ means
that the universal covering will be a product of disks. For this reason, and for apparent lack of any standard terminology,
we propose to call these ``polydisk Shimura varieties''. As with the case of DM-curves, it is more convenient for getting
good factorization statements, to ignore level structure and consider moduli stacks. 

Thus, we will review a version of Shimura's construction which we call {\em polydisk Shimura modular DM-stacks}. 
Aside from the restriction to the group $SL(2)$, this is standard material
\cite{DeligneTravauxDeShimura}
\cite{DeligneCorvallis}  \cite{Kottwitz} \cite{Larsen} 
\cite{Pappas} \cite{RapoportZink} \cite{Wedhorn},
but it seems
like a good idea to have a discussion which corresponds to our specific situation. 

Let $L$ be a totally imaginary extension of a totally real field $F$. Let $\Oo _L$ be the ring of algebraic integers
in $L$. Let $P$ be a projective $\Oo_L$-module of rank two, and put $P_L:= P\otimes _{\Oo_L}L$.

Suppose $\sigma : L \rightarrow \cc$ is an embedding. Let $P^{\sigma}:= P_L \otimes _L\cc$ be the 
tensor product using $\sigma$ as structural map. 

Let $\iota : L \rightarrow L$ be the complex conjugation map. It is independent of the embedding
(indeed it is the generator for the order $2$ group $Gal (L/F)$), so for any embedding $\sigma$ we have
$$
\sigma (\iota x) = \overline{\sigma x}.
$$
Let $\overline{\sigma} = \sigma \circ \iota$ denote the embedding composed with complex conjugation.

Now we assume that we are given an $\iota$-hermitian antisymmetric form $\Phi$ on $P_L$. This means 
$$
\Phi : P_L \times P_L \rightarrow L,
$$
such that $\Phi (u,v) = -\iota \Phi (v,u)$, with $\Phi$ separately additive in each variable and such that
$\Phi (a u, v) = a\Phi (u,v)$ thus also $\Phi (u,av) = (\iota a) \Phi (u,v)$. The collection 
$(L,P,\Phi )$ is basically what is known as a {\em PEL datum} \cite{Wedhorn}. 

Let $\Uu (P_L, \Phi)$ denote the group of $L$-linear transformations of $P_L$ which preserve $\Phi$,
and let $\Uu (P,\Phi )$ denote the subgroup of those transformations which also preserve $P\subset P_L$.
Let $\mcS \Uu (P_L, \Phi)$  and $\mcS \Uu (P,\Phi )$ denote the subgroups of transformations of determinant one. 

If $\sigma : L \rightarrow \cc$ is an embedding then we obtain a complex hermitian antisymmetric form 
on $P^{\sigma} := P_L\otimes _L\cc$
defined by
the condition that
$$
\Phi ^{\sigma}(u^{\sigma}, v^{\sigma}) = \sigma \Phi (u,v)
$$
where $u^{\sigma}$ and $v^{\sigma}$ denote the images of $u$ and $v$ in $P^{\sigma}$.

A {\em Hodge structure on $(P_L,\Phi )$} is the data $H$ for each embedding $\sigma$ of a Hodge decomposition 
$$
P^{\sigma} = H^{1,0}_{\sigma} \oplus H^{0,1}_{\sigma}, 
$$
such that this decomposition is polarized by the form $\Phi ^{\sigma}$
(recall that this means that $\Phi^{\sigma}$ is positive imaginary on $H^{1,0}_{\sigma}$ and
negative imaginary on $H^{0,1}_{\sigma}$, and that these two spaces are orthogonal), and such that
the decompositions on $P^{\sigma}$ and $P^{\overline{\sigma}}$ are complex conjugate.
Notate this complex conjugation condition as follows. We have an antilinear isomorphism 
$$
\overline{(\cdot )}: P^{\sigma} \leftrightarrow P^{\overline{\sigma}}, \;\;\; \overline{u^{\sigma}}:= u^{\overline{\sigma}}
$$
(recall that $\overline{\sigma}:= \sigma \iota$), and using this operation we require that
$$
H^{p,q}_{\overline{\sigma}} = \overline{H^{q,p}_{\sigma}}.
$$
To see why such structures exist, we need the following lemma.

\begin{lemma}
\label{polconj}
Suppose 
$$
P^{\sigma} = \bigoplus
H^{p,q}_{\sigma}
$$ 
is a Hodge
decomposition polarized by $\Phi ^{\sigma}$. Then setting 
$$
H^{p,q}_{\overline{\sigma}}:= \overline{H^{q,p}_{\sigma}}
$$
we obtain a Hodge decomposition 
$$
P^{\overline{\sigma}} = \bigoplus
H^{p,q}_{\overline{\sigma}}
$$ 
for $P^{\overline{\sigma}}$ polarized by $\Phi ^{\overline{\sigma}}$.
\end{lemma}
{\em Proof:}
Suppose $u,v\in P$ so that $u^{\sigma},v^{\sigma}\in P^{\sigma}$.
Then 
$$
\Phi ^{\overline{\sigma}}(\overline{u^{\sigma}}, \overline{v^{\sigma}}) = 
\overline{\sigma}\Phi (u,v) = \overline{\sigma \Phi (u,v)} = \overline{\Phi ^{\sigma}(u^{\sigma},v^{\sigma})}.
$$
This shows that in general for $x,y\in P^{\sigma}$ we have
$$
\Phi ^{\overline{\sigma}}(\overline{x},\overline{y}) =\overline{\Phi ^{\sigma}(x,y)}.
$$
Now if $x\in H^{p,q}_{\sigma}$ then $\overline{x}\in H^{q,p}_{\overline{\sigma}}$
and to get the polarization condition for $H^{q,p}_{\overline{\sigma}}$ we would like to show that
$$
(\sqrt{-1})^{q-p}\Phi ^{\overline{\sigma}}(\overline{x},\overline{x}) > 0.
$$
The quantity on the left is equal by the previous formula to
$$
(\sqrt{-1})^{q-p}\overline{\Phi ^{\sigma}(x,x)} = \overline{(\sqrt{-1})^{p-q}\Phi ^{\sigma}(x,x)}.
$$
The polarization condition for $H^{p,q}_{\sigma}$ says that 
$$
(\sqrt{-1})^{p-q}\Phi ^{\sigma}(x,x) > 0
$$
which yields the desired polarization condition for $H^{q,p}_{\overline{\sigma}}$.
The condition about orthogonality is proved similarly. 
\eop

Our forms $\Phi ^{\sigma}$ are hermitian antisymmetric, giving hermitian symmetric 
forms denoted $\sqrt{-1}\Phi^{\sigma}$. 

Say that $\sigma$ is {\em positive}, {\em mixed} or {\em negative} respectively
if $\sqrt{-1}\Phi ^{\sigma}$ is a positive definite, indefinite, or negative definite unitary form. 
Note that $\overline{\sigma}$ is respectively negative, mixed or positive
when $\sigma$ is positive, mixed or negative (a special case of the above lemma).

If $\sigma$ is positive (resp. negative) then we are forced to set $H^{1,0}_{\sigma} :=P^{\sigma}$ 
(resp. $H^{0,1}_{\sigma} :=P^{\sigma}$) so in these cases the Hodge decompositions are uniquely determined.
In the mixed case, each Hodge subspace has dimension $1$ and the space of such decompositions for a given 
$\sigma$ is the unit disc. The complex conjugacy condition determines the Hodge structure on 
$P^{\sigma\iota}=P^{\overline{\sigma}}$
once it is given for $P^{\sigma}$, and the complex conjugate Hodge structure on $P^{\sigma\iota}$
is polarized by $\Phi ^{\sigma\iota}$ if and only if the original one was polarized by $\Phi ^{\sigma}$,
by Lemma \ref{polconj}.
Therefore the space of Hodge decompositions of $(P_L,\Phi )$ is
isomorphic to a product of discs, with the number of factors equal to half the number of mixed embeddings for $\Phi$.

If $X$ is a smooth analytic variety, then a {\em VHS of type $(P,\Phi )$ over $X$} is
a local system $V$ of $\Oo _L$ modules whose stalks are isomorphic to $P$, together with an antihermitian
form $\Phi _V$ on $V$ such that for any $x\in X$ the stalk $(V_x, \Phi _{V,x})$ is isomorphic to 
$(P,\Phi )$ (and in particular $V$ is a local system of projective rank two $\Oo _L$-modules), together
with a family of Hodge decompositions of the $(V_x, \Phi _{V,x})$ in the above sense which for each $V^{\sigma}$
forms a complex variation of Hodge structure polarized by $\Phi ^{\sigma}_V$. This notion can be extended
in an obvious way to the case where the base $X$ is a smooth analytic Deligne-Mumford stack. 

Let $\Dd (P, \Phi )$ denote the period domain parametrizing Hodge decompositions on $(P, \Phi )$ as above.
The group $\Uu (P, \Phi)$ acts on $\Dd (V_L, \Phi )$. Put
$$
\Hhh ^{\rm an}(P,\Phi ):= \Dd (P, \Phi ) /  \Uu (P, \Phi).
$$

\begin{theorem}
\label{shimstackan}
The analytic stack $\Hhh ^{\rm an}(P,\Phi )$ is a smooth analytic Deligne-Mumford stack
(i.e. an orbifold). This means in particular that the stabilizers are finite.
It admits a tautological VHS of type $(P,\Phi )$ denoted $\Vv$.  This represents the functor of
VHS's of type $(P,\Phi )$  in the following sense: if $X$ is any smooth analytic variety (or analytic
DM-stack) and $(V,\Phi _V)$ is a VHS of type $(P,\Phi )$ over $X$, then there is a map
$f:X\rightarrow \Hhh ^{\rm an}(P,\Phi )$ and an isomorphism 
$$
\alpha : (V,\Phi _V)\cong f^{\ast} (\Vv,\Phi _\Vv).
$$
Furthermore, given two maps $f,g: X\rightarrow \Hhh ^{\rm an}(P,\Phi )$, the
set of isomorphisms between $f$ and $g$ maps isomorphically to the set of 
isomorphisms between  $f^{\ast} (\Vv,\Phi _\Vv)$ and
$g^{\ast} (\Vv,\Phi _\Vv)$.
\end{theorem}

\begin{theorem}
\label{shimstackalg}
There exists an algebraic Deligne-Mumford stack $\Hhh (P,\Phi )$ and an isomorphism between its associated
analytic space and the complex DM-stack $\Hhh ^{\rm an}(P,\Phi )$ defined above. The algebraic stack
represents the functor of VHS's of type $(P,\Phi )$ over algebraic varieties (with the same statement as
before which we don't copy down a second time). In particular, if $X$ is any smooth algebraic DM-stack,
then any analytic map $X^{\rm an}\rightarrow \Hhh ^{\rm an}(P,\Phi )$ comes from an algebraic map
(and the analytic and algebraic isomorphisms between algebraic maps are the same). 
\end{theorem}

This is a classical result of Baily-Borel type. In the Hilbert modular case see Rapoport \cite[Theorem 1.20]{Rapoport}. 
Rapoport attributes the materiel of his \S 1 to Deligne. 
For Shimura varieties of PEL type, Wedhorn \cite{Wedhorn} refers to Kottwitz \cite{Kottwitz}
who states a moduli problem in terms of abelian varieties with appropriate structure, 
in our case that would mean purely imaginary multiplication by $L$.
 
The consequence about unicity
of the algebraic structure seems even now a bit surprising. It basically comes from the distance-decreasing 
property for variations of Hodge structure, which implies that any analytic map is 
distance decreasing for the Poincar\'e metric. 

The stacks $\Hhh (P,\Phi )$ appearing in Theorem \ref{shimstackalg} will be called {\em polydisk Shimura modular DM-stacks}. Note that
there is a tautological representation
$$
\pi _1(\Hhh (P,\Phi )) \rightarrow \Uu (P, \Phi) \hookrightarrow SL(2,L)
$$
which gives one tautological representation $\pi _1(\Hhh (P,\Phi )) \rightarrow \rightarrow SL(2,\cc )$ for
each embedding $\sigma : L \hookrightarrow \cc$. 

A necessary remark is that our terminology is meant to suggest that these are the standard examples of Shimura varieties whose
universal covering is a polydisk. However, we don't claim that {\em every} Shimura variety of PEL type whose universal covering is
a polydisk, is of the above form (or else Hilbert modular). It would be good to clarify this. 

Ideally one should also discuss the notion of {\em level structure} here, of 
a type generically denoted by $N$ for the purposes of the present paragraph. This allows one to get rid of the 
stack structure by looking at finite etale
coverings $\Hhh (P,\Phi ,N)\rightarrow \Hhh (P,\Phi )$ by covering stacks which are actually quasiprojective varieties.
One would then need to define the notion of {\em VHS of type $(P,\Phi )$ with level structure} to correspond
to morphisms into the polydisk Shimura modular varieties with level structure $\Hhh (P,\Phi ,N)$. 
For a given map $X\rightarrow 
\Hhh (P,\Phi )$ one could try to describe the explicit finite etale covering of $X$ which would lift to a map
into $\Hhh (P,\Phi ,N)$. Such a discussion would have the benefit of giving a factorization result in the world
of varieties rather than in the world of stacks. However, choices would have to be made and the factorization
would only take place over a covering of $X$. Just as we have amply seen in the case of curves, the formulation
in terms of stacks gives a much more natural and cleaner factorization statement.

\section{Improvements for a rigid integral local system}
\label{sec-improvements}

Let $\Upsilon := \pi _1(X,x)$. 
Suppose $\rho : \Upsilon \rightarrow SL(2,\cc )$ is a rigid, integral Zariski-dense representation. 
As was pointed out by M. Larsen (see \cite{hbls} Lemma 4.8), we can assume after simultaneous conjugation that $\rho$ is 
defined over a totally imaginary quadratic extension
$L$ of a totally real field $F$.  In other words $\rho$ comes from 
$\rho _L : \Upsilon \rightarrow SL(2,L)$ by extension of scalars via an embedding $\eta : \rightarrow \cc$. 

Furthermore, by Bass-Serre theory \cite{Bass} the fact that $\rho$ is integral means that
there is a rank two projective $\Oo _L$-module $P$ with $P\otimes _{\Oo _L}L\cong L^2$ such that $\rho _L$ comes from 
a representation $\rho _P: \Upsilon \rightarrow SL(P)$. 

We adopt the same notations as in the previous section, for example $\iota$ denotes the complex conjugation on $L$. 
If $\sigma : L\rightarrow \cc$ is any embedding then we obtain a complex vector space $P^{\sigma} := P\otimes _{\Oo _L}\cc$
and a representation $\rho ^{\sigma}:\Upsilon \rightarrow SL(P^{\sigma})$ as in the previous section. 
For the given embedding $\sigma = \eta$ this gives back the original representation. If $\rho$ was rigid 
then all of the $\rho ^{\sigma}$ are rigid too. Thus, they all correspond to complex variations of Hodge structure.

This is almost the data we need in order to define a Shimura modular VHS and get a map from $X$ into 
a polydisk Shimura modular stack $\Hhh (P,\Phi )$. But we have to discuss the question of the polarization.

\subsection{Uniqueness of the polarization}

Let $V$ be a vector space of dimension two over $L$, with an action of $\Upsilon$. 
Define the {\em conjugate representation} $V^{\iota}$ as follows. The underlying additive group is the same
as that of $V$, but the scalar multiplication operation operation $\mu$ is defined by
$$
\mu _{V^{\iota }} (u,x):= \mu _{V}(\iota u, x).
$$
The action of $\Upsilon$ is the same as on $V$, and it is again by $L$-linear endomorphisms.

A sesquilinear form $\Phi$ on $V$ is a map $\Phi : V\times V \rightarrow L$ which is separately additive
in each variable and which satisfies
$$
\Phi (ax, y) = a\Phi (x,y),
$$
and 
$$
\Phi (x, ay) = (\iota a)\Phi (x,y).
$$
It is {\em symmetric} (resp. {\em antisymmetric}) if 
$$
\Phi (y,x) = \iota \Phi (x,y) \;\; \left( \mbox{resp.} \; \Phi (y,x) = -\iota \Phi (x,y) \right) .
$$
We say that $\Phi$ is {\em $\Upsilon$-invariant} if $\Phi (\gamma x, \gamma y) = \Phi (x,y)$ for any
$\gamma \in \Upsilon$ and $x,y\in V$.

Define the {\em tensor corresponding to $\Phi$} denoted 
$$
{\bf t}\Phi \in V^{\ast} \otimes _L ( V^{\iota})^{\ast},
$$
to be the unique element such that 
$$
{\bf t}\Phi \cdot x\otimes y = \Phi (x,y)
$$
for any $x\in V$ and $y\in \iota V$ (the latter being the same as saying $y\in V$); and
where the product is the natural contraction 
$$
\left( V^{\ast} \otimes _L ( V^{\iota})^{\ast} \right) \times 
\left( V \otimes _L (V^{\iota}) \right) \rightarrow L. 
$$
There is an action of $\Upsilon$ on $V^{\ast} \otimes _L ( V^{\iota})^{\ast}$ and
$\Phi$ is $\Upsilon$-invariant if and only if ${\bf t}\Phi$ is fixed by this action.

Define the $\iota$-antilinear involution $C: V^{\ast} \otimes _L ( V^{\iota})^{\ast} \rightarrow 
V^{\ast} \otimes _L (V^{\iota})^{\ast}$ with $C^2=1$ by
$$
C (e \otimes f) = - f \otimes e.
$$
If ${\bf t}\Phi = e\otimes f$ this means that $\Phi (x,y) = e(x)f(y)$.
Define $\Phi ^C$ by ${\bf t}(\Phi ^C) = C({\bf t}\Phi )$. Then 
$$
{\bf t}\Phi^C = - f\otimes e
$$
so $\Phi ^C(x,y) = - f(x)e(y) = -\Phi (y,x)$. By linearity this extends to all $\Phi$ so
$$
\Phi ^C(x,y) = - \Phi (y,x).
$$
In particular, $\Phi ^C = \Phi$ or equivalently, ${\bf t}\Phi$ is fixed by $C$, if and only if
$\Phi$ is antisymmetric. 

The involution $C$ defines an $F$-structure on the $L$-vector space $V^{\ast} \otimes _L (V^{\iota})^{\ast}$ 
such that the $F$-valued points (i.e. those fixed by $C$) are exactly the tensors corresponding to
antisymmetric sesquilinear forms. 

\begin{lemma}
If $V$ is an irreducible representation of $\Upsilon$ then the dimension of the $F$-vector space
of points of $V^{\ast} \otimes _L (V^{\iota})^{\ast}$ 
fixed by $\Upsilon$ and $C$ is at most one.
\end{lemma}
{\em Proof:}
The $F$-vector space is a reduction to $F$ of the $L$-vector space of points fixed only by $\Upsilon$.
Thus it suffices to prove that the space of points fixed by $\Upsilon$ has dimension $\leq 1$ over $L$.
But this space is the same as the space of $\Upsilon$-equivariant $L$-linear maps from $V^{\iota}$ to $V^{\ast}$,
and since both are irreducible representations of $\Upsilon$, the dimension of the space of maps is at most $1$. 
\eop

Fix an embedding $\sigma : L\hookrightarrow \cc$ inducing $F\hookrightarrow \rr$.
Everything in the above discussion can be tensored over $F$ with $\rr$.
The field $L$ becomes $\cc$ and $V\otimes _F\rr = V\otimes _L\cc$ is a two-dimensional $\cc$-vector space.
We have 
$$ 
( V^{\iota})\otimes _F\rr = \overline{(V\otimes _F\rr )},
$$ 
and
$$
\left( V^{\ast} \otimes _L ( V^{\iota})^{\ast} \right) \otimes _F\rr = 
$$
$$
(V\otimes _F\rr )^{\ast} \otimes _{\cc} \overline{(V\otimes _F\rr )}^{\ast}.
$$
The involution $C$ becomes the involution with a similar formula on this space, whose
fixed points are again the tensors corresponding to antisymmetric sesquilinear forms. 

The $\rr$-dimension of the space of points in  
$\left( V^{\ast} \otimes _L (\iota V)^{\ast} \right) \otimes _F\rr $
fixed by $\Upsilon$ and $C$ is equal to the $F$-dimension of the space of fixed points before tensoring.
By the above lemma, it is $0$ or $1$. 

We conclude the following facts from this discussion. 

\begin{proposition}
\label{thescholium}
If for some $\sigma$ the representation $V\otimes _{L,\sigma}\cc$ admits a
$\Upsilon$-invariant antisymmetric sesquilinear form $\Phi _0$ then the space of such forms is
a one-dimensional real vector space. If $\Phi _0$ exists then  $V$ admits a $\Upsilon$-invariant antisymmetric 
$\iota$-sesquilinear form
$\Phi$ and the space of such forms is again a one-dimensional $F$-vector space. For any $\Phi$  there is a
real constant $\lambda$ such that $\Phi _0 = \lambda \Phi ^{\sigma}$. 
\end{proposition}
\eop

\subsection{Creating a polarized VHS}

We go back to the situation considered at the start of this section. We have a rank two projective $\Oo _L$-module $P$
with a representation
$$
\rho _P:\pi _1(X,x)\rightarrow SL(P)
$$
yielding 
$$
\rho ^{\sigma}:\pi _1(X,x)\rightarrow SL(P^{\sigma})\cong SL(2,\cc )
$$
for any embedding $\sigma : L\rightarrow \cc$.  
Our original representation $\rho$ is
equal to $\rho ^{\eta}$.

Let $V:=P\otimes _{\Oo _L}L$ be the two dimensional $L$-vector space on which
$\rho _L$ acts. In the notation of the previous subsection $V\otimes _{L,\sigma}\cc = P^{\sigma}$. Let $A_L$ be the
local system of $L$-vector spaces of rank two corresponding to $\rho _L$ and for any $\sigma$ 
let $A_{\sigma} = A_L\otimes _{L,\sigma}\cc$
be the local system of complex vector spaces corresponding to the representation $\rho ^{\sigma}$. 

\begin{lemma}
\label{choicevhs}
We can choose for each $\sigma$ a structure of $\cc$VHS of weight one on the local system $A_{\sigma}$ such that
$$
A^{1,0}_{\overline{\sigma}} = \overline{A^{0,1}_{\sigma}}, 
\;\;\;
A^{0,1}_{\overline{\sigma}} = \overline{A^{1,0}_{\sigma}}.
$$
Furthermore, these complex variations have hermitian polarizations.
\end{lemma}
{\em Proof:}
By hypothesis $\rho$ is rigid, which implies that each $\rho ^{\sigma}$ is rigid. Thus
each $\rho ^{\sigma}$ underlies an hermitian polarized $\cc$VHS (Theorem \ref{rigidvhs}).
The fact that we are in rank two means that we can choose these variations to have weight one. 
If $\sigma$ is mixed then
the structure has to have weights $(1,0)$ and $(0,1)$ so it is unique. In this case $\overline{\sigma}$ is
also mixed, and since the complex conjugate variation is one possible choice for $A_{\overline{\sigma}}$ it is
the only choice and we get the desired formula. 

On the other hand, if $\sigma$ is positive
or negative, then we may choose the structure as we like. In this case $\overline{\sigma}=\sigma 
\circ \iota$ has the opposite sign. Thus we may assume that the complex conjugate choice is made for 
$\overline{\sigma}$, which again gives the formula. 
\eop

Fix a choice of variations $A_{\sigma}$ as in the lemma. 
The fact that $A_{\sigma}$ has an hermitian polarization implies that
there exists a $\pi _1(X,x)$-invariant antisymmetric sesquilinear form on $V\otimes _{L,\sigma}\cc$. 
By Proposition \ref{thescholium}, there also exists a $\pi _1(X,x)$-invariant antisymmetric $\iota$-sesquilinear
form $\Phi $ on $V$, and it is unique up to multiplication by scalars in $F$. 

The remaining problem is to get the right signs for the polarizations on the various pieces. 
For any embedding $\sigma : L \rightarrow \cc$, the image $\Phi ^{\sigma}$ is 
a real multiple of the polarization form of the variation of Hodge structure $A_{\sigma}$ chosen in Lemma
\ref{choicevhs} above. 
Thus there is a sign function depending on $\Phi$ and denoted 
$$
g(\Phi , \sigma ) \in \{ \pm 1\}
$$ 
such that
$g(\Phi , \sigma )\Phi ^{\sigma }$ polarizes the complex variation of Hodge structure 
$A_{\sigma }$.

\begin{lemma}
\label{fixsigns}
There is a $\pi _1(X,x)$-invariant antisymmetric sesquilinear form $\Phi$ on $V$ such that
$g(\Phi , \sigma )=1$ for all embeddings $\sigma : L\hookrightarrow \cc$, in other words
$\Phi ^{\sigma }$ polarizes the $\cc$VHS $V^{\sigma}$ for each mixed embedding $\sigma$.
\end{lemma}
{\em Proof:}
Fix one choice of $\Phi$ which might not have the right signs. 
We claim that for each $\sigma$, $g(\Phi , \sigma \circ \iota ) = g(\Phi , \sigma )$.  Recall that we have chosen the
variations of Hodge structure $A_{\sigma}$ in Lemma \ref{choicevhs} so that $A_{\overline{\sigma}}=\overline{A_{\sigma}}$.
By Lemma \ref{polconj} the form $\Phi^{\sigma}$ polarizes $A_{\sigma}$ if and only if $\Phi ^{\overline{\sigma}}$
polarizes $A_{\overline{\sigma}}$. The same is true if we replace $\Phi$ by $-\Phi$. Thus $g(\overline{\sigma})=g(\sigma )$
as claimed. 

For each embedding $\zeta : F\rightarrow \rr$ there are exactly two embeddings $\sigma , \overline{\sigma}:L\rightarrow \cc$
which induce $\zeta$ on $F$. The claim in the previous paragraph says $g(\Phi , \sigma )$ only depends on $\zeta = \sigma |_F$.
Abusing notation we write $g(\Phi , \zeta )$. 

Next let ${\bf sgn}(r)\in \{ \pm 1\}$ denote the sign of a nonzero real number.
If $\lambda \in F ^{\ast}$ then 
$$
g(\lambda \Phi , \zeta ) = {\bf sgn}(\zeta \lambda ) g(\Phi ,\zeta ).
$$
Therefore, if we can find $\lambda \in F^{\ast}$ such that ${\bf sgn}(\zeta \lambda ) = g(\Phi , \zeta )$ for 
all $\zeta : F\rightarrow \cc$, then the form $\lambda \Phi$ will be a  polarization form with the right signs for the lemma. 

Now recall that
$$
F\otimes _{\qq}\rr = \prod _{\zeta : F\rightarrow \rr} \rr .
$$
Since $F$ is dense in $F\otimes _{\qq}\rr$, there is a point $\lambda \in F$ which lies in the open quadrant of 
$\prod _{\zeta : F\rightarrow \rr} \rr $ determined by the specification ${\bf sgn}(\zeta \lambda ) = g(\Phi , \zeta )$.
Thus $\lambda\Phi$ answers the question. 
\eop

Fix a form $\Phi$ as in the lemma. If $\sigma$ is a non-mixed embedding then
the $\cc$VHS $V^{\sigma}$ is unitary, so it can be chosen as equal to a single Hodge subspace
of Hodge type either $(1,0)$ or $(0,1)$. We can choose which one so that $V^{\sigma}$ with this
Hodge structure, is polarized by $\sqrt{-1}\Phi ^{\sigma }$. With this condition the choice
is unique, and the complex conjugate embedding corresponds to the opposite choice of Hodge type
(Lemma \ref{polconj}).
Thus, if $\Phi$ is fixed as in the lemma, there is unique choice of collection of complex variations
of Hodge structure for the $V^{\sigma}$ such that the are all polarized by $\sqrt{-1}\Phi ^{\sigma }$.

We can clear denominators in $\Phi$ by multiplying it by an appropriate integer. Thus we may assume that
$\Phi : P\times P \rightarrow \Oo _L$ is an integral polarization form. We have now obtained the following proposition.

\begin{proposition}
\label{improvement}
Suppose $\rho : \pi _1(X,x) \rightarrow SL(2,\cc )$ is a rigid, integral Zariski-dense representation. 
Then there exists:  a totally imaginary extension $L$ of a totally real field $F$; a projective rank two $\Oo _L$-module
$P$ together with an $\iota$-sesquilinear antisymmetric form $\Phi : P\times P \rightarrow \Oo _L$; a representation 
$\rho _P : \pi _1(X,x)\rightarrow \mcS \Uu (P,\Phi )$; and a variation of Hodge structure $A = \{ A_{\sigma}\}$ of type $(P,\Phi )$
polarized by $\Phi$ and whose underlying representation is $\rho _P$; such that for one of the embeddings $\eta : L\rightarrow
\cc$ we have $\rho = \rho ^{\eta}$. 
\end{proposition}

By Theorem \ref{shimstackalg} the VHS of type $(P,\Phi )$ comes from a map $X\rightarrow \Hhh (P,\Phi )$ to the Shimura modular stack
corresponding to $(L,P,\Phi )$. 

\section{Classification}
\label{sec-classif}

We sum up our main classification results.  
Suppose $X$ is a quasiprojective variety with normal crossings compactification $\overline{X}$ and
$\overline{X}-X = D$. Fix a representation $\rho :\pi _1(X,x) \rightarrow SL(2,\cc )$.
Suppose that $\rho$ is quasi-unipotent at infinity, i.e. the monodromy elements $\rho (\gamma _i)$
are quasi-unipotent for loops $\gamma _i$ going around the components $D_i$ of the divisor at infinity.
We also suppose that the image of $\rho$ is Zariski-dense in $SL(2,\cc )$.

Recall that we say $\rho$ is {\em rigid} if there are no nontrivial global
deformations which fix the conjugacy
classes of the $\rho (\gamma _i)$, see Lemma \ref{rigidequiv}, and $\rho$ is {\em non-rigid} otherwise. 

\begin{theorem}
\label{nonrigidcase}
Suppose with the above hypotheses, that $\rho$ is {\em nonrigid}. 
Then there exists a DM-curve $Y$, 
a rank two representation $\rho _Y$ of $\pi _1(Y)$, and a map 
$$
f:X\rightarrow Y
$$
such that $\rho \cong f^{\ast} \rho _Y$. 
\end{theorem}
{\em Proof:}
This the contrapositive of Theorem \ref{nonfactorsrigid}.
\eop

\begin{theorem}
\label{rigidcase}
Suppose with the above hypotheses, that $\rho$ is {\em rigid}. 
Then $\rho$ is also integral; and there exists a totally imaginary field
$L$ with an embedding $\eta : L\rightarrow \cc$, 
a rank two projective $\Oo _L$-module $P$ with an antihermitian polarization form $\Phi$,
and a map to the polydisk Shimura modular stack
$$
f:X\rightarrow \Hhh (P,\Phi )
$$
such that if $\Vv^{\eta }$ is the complex local system on $\Hhh (P,\Phi )$
corresponding to the tautological 
local system tensored with $\cc$ over the embedding $\eta$, then $\rho$ is conjugate
to the monodromy representation of 
$f^{\ast }\Vv ^{\eta }$.
\end{theorem}
{\em Proof:}
By Lemma \ref{DMrigidintegral}, $\rho$ is integral. Proposition \ref{improvement} gives 
data $(L,P,\Phi )$ and a 
variation of Hodge structure $V$ of type $(P,\Phi )$ over $X$.
If $\rho _P : \pi _1(X,x)\rightarrow \mcS \Uu (P,\Phi )$ is the representation underlying this 
variation of Hodge structure then from \ref{improvement}, for one of the embeddings $\eta : L\rightarrow
\cc$ we have $\rho = \rho _P^{\eta}$. 

Now Theorems \ref{shimstackan} and \ref{shimstackalg}
say that the variation $V$ on $X$ comes from a map $f:X\rightarrow \Hhh (P,\Phi )$
in the sense that the pullback of the tautological variation $f^{\ast}(\Vv )$ is the variation $V$ on $X$. 
The tautological rank two complex local systems over $\Hhh (P,\Phi )$ all come from the tautological
local system of $\Oo _L$-modules of type $P$ underlying $\Vv$, via the embeddings $L\hookrightarrow \cc$. Thus
the representation $\rho \cong \rho _P^{\eta}$ is the pullback of the tautological rank two complex
local system $\Vv ^{\eta}$ on  $\Hhh (P,\Phi )$.
\eop

{\em Proof of Theorem \ref{classification-global}:} Suppose $\rho$ is a Zariski dense representation
which is quasi-unipotent at infinity. It is either rigid, in which case Theorem \ref{rigidcase}
says that it comes from a polydisk Shimura modular stack; or else nonrigid, in which case Theorem \ref{nonrigidcase}
says that it comes from a DM-curve. 
\eop

There are some cases where a representation on a DM-curve $Y$ can be rigid, in particular it comes from a map
from $Y$ to a Shimura modular stack. These cases were discussed in \S \ref{sec-hypergeometric}. In these
cases, Theorem \ref{rigidcase} applies. 

\begin{corollary}
\label{rigidmotivic}
Any rigid local system of rank two, quasiunipotent at infinity, with Zariski dense monodromy, is motivic.
\end{corollary}
{\em Proof:}
The tautological local systems on Shimura modular stacks are motivic.
\eop

\section{Questions}
\label{sec-questions}

\subsection{Further development}

The most obvious question is: can 
we give a similar classification, or at least a start, for rank $3$ representations? 

Aside from this, 
several questions are left open concerning the case of local systems of rank $2$. 
For example, what happens in the non-quasi-unipotent case? And, what about local systems with smaller 
monodromy groups, notably a solvable monodromy group? This question is addressed by Delzant in \cite{Delzant} 
\cite{Delzant2}.   

For a given DM-curve we would like to have a
good understanding of the irreducible components of the moduli spaces of rank $2$ local systems. 

In the present treatment we have restricted to the
case where the original $X$ was a smooth quasiprojective variety. What about the case where $X$ is a general DM-stack,
possibly with singularities?

\subsection{Calculating the factorization}

The transcendental nature of our methods leads naturally to the question of whether one can effectively
(and preferably, easily) calculate the map to a Shimura modular variety corresponding to a given rank $2$
local system. Concretely, $V$ might be given in any number of ways: 
\newline
(A)\, by representation matrices corresponding
to the generators of the fundamental group (in a Zariski-type calculation of the fundamental group of $X$);
\newline
(B)\, as a vector bundle (probably the trivial bundle, at least if we allow going to a neighborhood in $X$)
plus an explicitly given regular-singular connection;
\newline
(C)\, as a stable parabolic Higgs bundle with vanishing parabolic Chern classes, which might in many cases take on 
one of several special forms:
\newline
(C1)\, a stable vector bundle with $\theta =0$ over a projective variety; 
\newline
(C1')\, a stable parabolic vector bundle with $\theta = 0$;
\newline
(C2)\, a stable system of Hodge bundles, which in the rank two case and excepting (C1) above means
a pair of line bundles and a map 
$$
\theta : L \rightarrow M \otimes \Omega ^1_X
$$
\newline
(C2')\, a stable system of Hodge bundles with parabolic structure and/or logarithmic $\theta$. 

For each kind of concrete example of $V$ there could arise the question of how to explicitly construct
the map from $X$ (or a finite covering space $X'$ of $X$) to a Shimura modular variety plus an isomorphism
between our given object and the pullback of the universal object. 

While it might be possible to prove using some general principles that there exists an algorithm
which solves this problem in principle, we would actually like to have a good method which works 
on specific examples and which is related to geometry \cite{AllcockCarlsonToledo}
\cite{AmorosBauer} \cite{Dimca} \cite{DimcaNemethi} \cite{FriedmanMorgan} \cite{KapovichMillson} \cite{RobbTeicher} \cite{Zuo}.

There are now known some specific examples. For example D. Panov constructs rank two unitary representations
corresponding to polyhedral K\"ahler structures on some surfaces in his thesis \cite{Panov}.

Any algebraic solution of the Painlev\'e VI equation will correspond to a rank two local system on a surface.
Such solutions have been constructed by Hitchin \cite{Hitchin}, Dubrovin \cite{Dubrovin}, Boalch \cite{Boalch}, 
Ben Hamed-Gavrilov \cite{BenHamedGavrilov}.
In many cases the solutions were shown to have geometric origin already. 

More generally there are various potential sources of other examples 
\cite{BogomolovKatzarkov} \cite{Campana04} \cite{DeligneMostow} \cite{DimcaNemethi} \cite{Moishezon}
\cite{Nori} \cite{SommeseVandeVen} \cite{Toledo}, in particular line or hyperplane
arrangements as was the case in Panov's construction \cite{AllcockCarlsonToledo} 
\cite{CohenOrlik} \cite{FalkYuzvinsky} \cite{MoishezonTeicher} \cite{PapadimaSuciu}. These might in the future lead to
constructions of rank two local systems for which the classification results would apply. 

To what extent is the factorization map uniquely determined by the representation? This is easy 
for factorization through an orbicurve, but might be an interesting question for factorization through
a Shimura modular stack. What kind of Shimura modular stack is associated to a given local system? 
For example what is its dimension? Are there details specific to the theory of Hilbert modular varieties,
which could be extended to the polydisk Shimura case?
On the other side, how do we tell given $(P,\Phi )$ how many mixed embeddings there are?

The problem of calculating explicitly the factorization is also interesting when $X$ is a curve. 
In \S \ref{sec-geom} we have started to look at the real geometry of the pluriharmonic map this case. This was very preliminary
but it looks
interesting and should be further pursued. How does the representation $\rho : \pi _1\rightarrow SL(2,K)$
determine the geometry of the map from $\widetilde{X}$ to the tree?

The analytic proof that nonrigid representations factor, which we have sidestepped here, also leads to 
an effectivity question. 
The space of representations whose Higgs fields are nilpotent is
compact---can we give an explicit bound for the size of these representations?

\end{document}